\documentclass{ucbthesis}

\usepackage{amsmath,amssymb,amsthm,hyperref,enumerate,stmaryrd}

\usepackage{tikz}
\usepgflibrary[arrows]
\usetikzlibrary{arrows,matrix}

\newtheorem{conj}{Conjecture}[section]
\newtheorem{prop}{Proposition}[section]
\newtheorem{thrm}{Theorem}[section]
\newtheorem{lem}{Lemma}[section]
\newtheorem{cor}{Corollary}[section]

\theoremstyle{remark}
\newtheorem{rmk}{Remark}[section]

\theoremstyle{definition}
\newtheorem{defn}{Definition}[section]
\newtheorem{ex}{Example}[section]
\newtheorem*{ques}{Question}

\DeclareMathOperator{\Man}{\mathbf{Man}}
\DeclareMathOperator{\Symp}{\mathbf{Symp}}

\DeclareMathOperator{\dom}{dom}
\DeclareMathOperator{\im}{im}

\newcommand{\R}{\mathbb R}
\newcommand{\F}{\mathcal F}
\newcommand{\X}{\mathcal X}
\renewcommand{\O}{\mathcal O}
\newcommand{\g}{\mathfrak g}
\newcommand{\h}{\mathfrak h}
\newcommand{\wt}{\widetilde}
\newcommand{\red}{\twoheadrightarrow}
\newcommand{\cored}{\rightarrowtail}

\begin{document}


\title{Double Groupoids, Orbifolds, and the Symplectic Category}
\author{Santiago Valencia Ca\~nez}
\degreesemester{Spring}
\degreeyear{2011}
\degree{Doctor of Philosophy}
\chair{Professor Alan Weinstein}
\othermembers{Professor Peter Teichner \\
  Professor Robert Littlejohn}
\numberofmembers{3}
\field{Mathematics}
\campus{Berkeley}




\maketitle
\copyrightpage

\begin{abstract}
Motivated by an attempt to better understand the notion of a symplectic stack, we introduce the notion of a \emph{symplectic hopfoid}, which should be thought of as the analog of a groupoid in the so-called symplectic category. After reviewing some foundational material on canonical relations and this category, we show that symplectic hopfoids provide a characterization of symplectic double groupoids in these terms. Then, we show how such structures may be used to produce examples of symplectic orbifolds, and conjecture that all symplectic orbifolds arise via a similar construction. The symplectic structures on the orbifolds produced arise naturally from the use of canonical relations.

The characterization of symplectic double groupoids mentioned above is made possible by an observation which provides various ways of realizing the core of a symplectic double groupoid as a symplectic quotient of the total space, and includes as a special case a result of Zakrzewski concerning Hopf algebra objects in the symplectic category. This point of view also leads to a new proof that the core of a symplectic double groupoid itself inherits the structure of a symplectic groupoid. Similar constructions work more generally for any double Lie groupoid---producing what we call a \emph{Lie hopfoid}---and we describe the sense in which a version of the ``cotangent functor'' relates such hopfoid structures.
\end{abstract}

\begin{frontmatter}

\begin{dedication}
\null\vfil
\begin{center}
\emph{to Adelina and my parents}\\\vspace{12pt}
\end{center}
\vfil\null
\end{dedication}

\tableofcontents

\begin{acknowledgements}
First and foremost, I wish to thank my advisor, Alan Weinstein, for all the guidance and support he has provided me over the years. His love of mathematics and geometric intuition are second to none, and it has been an honor working under his supervision. Thanks for listening to all my wild ideas, and for sharing with me some of your own. I can only hope that one day I will be as effective a teacher as you have been.

Thanks to my mathematical brothers---Aaron, Hanh, Beno\^it, and Sobhan---for the many great times we have shared together. Our group meetings have always been a source of inspiration, and I have learned a great deal from you all. Thanks to the visitors our group has welcomed at various points over the last few years: Camilo, Leandro, Ludovic, Sean, Chiara, and Mitsuhiro. Thanks also to all the students I have had the pleasure to teach as a graduate student for helping to keep me motivated.

Thanks to my wife Adelina for all the love and support she has shared with me these past eight years; thanks for keeping me grounded and calm, and I look forward to many more years of you doing the same. Thanks to my mom and dad, who provided me with a nurturing and loving environment growing up, and a special thanks to my dad who has acted as both mother and father to me for so long. Thanks also to my friends---even if the majority of you are from the midwest---for helping to keep me sane.

Finally, thanks to television and the internet for being awesome, and to Apple for making products which are both incredibly useful and incredibly entertaining.
\end{acknowledgements}

\end{frontmatter}

\pagestyle{headings}


\chapter{Introduction}

Suppose that $\X$ is a differentiable stack \cite{BX} presented by a Lie groupoid $G \rightrightarrows M$.\footnote{We recall the basics of groupoids and stacks in Appendix~\ref{appendix}.} The tangent stack $T\X$ \cite{H} is then easy to construct: one simply applies the tangent functor to the groupoid $G \rightrightarrows M$ to obtain the tangent groupoid $TG \rightrightarrows TM$, and $T\X$ is the differentiable stack presented by this groupoid. Two natural questions then arise:
\begin{itemize}
\item Is there a similar construction of the cotangent stack $T^*\X$, and
\item What precisely is the cotangent stack, and in what sense is it symplectic?
\end{itemize}
This thesis arose out of an attempt to address these questions. While we do not give complete answers, we introduce a new type of structure---that of a \emph{symplectic hopfoid}---which seems to capture some of the data necessary in defining the notion of a ``symplectic stack''.

The notion of the cotangent stack of a stack is a subtle one in general. On the one hand, such a stack is usually not needed in practice since one can define the notion of a differential form on a differentiable stack (see \cite{BX} and \cite{LM}) via sheaf-theoretic means and without reference to a ``cotangent stack''. However, it is conceivable that for other purposes---quantization for example---having such a stack could well prove useful. In the algebraic category, definitions of cotangent stacks of algebraic stacks can be found in various place throughout the literature (see for example \cite{BD}), but these definitions do not easily carry over to the differential-geometric setting. Similarly, the notion of a symplectic structure on a differentiable stack also requires care. One aim of this thesis is to provide evidence that such structures may have nice descriptions in terms of the so-called ``symplectic category'', which we soon define.

Naively, one may attempt to construct such a cotangent stack as in the case of tangent stacks by applying the cotangent functor $T^*$ to the groupoid $G \rightrightarrows M$, thereby producing a groupoid $T^*G \rightrightarrows T^*M$. Of course, one immediately runs into the problem that such a ``cotangent functor'' does not actually exist, owing to the fact that in general a smooth map $G \to M$ does not induce a smooth map $T^*G \to T^*M$; indeed, such a map exists only when $G \to M$ is \'etale, which suggests that this construction might work in the case where $\X$ is actually an \emph{orbifold}. Even if we instead try to view $T^*$ as a contravariant functor, there is still no well-defined map $T^*M \to T^*G$ in general.

One way to try to get around this problem is to allow a more general type of morphism between cotangent bundles. Taking the point of view that specifying a map is the same as specifying its graph, we are led to the idea that general relations (or in the smooth setting, smooth relations) should be thought of as ``generalized maps''. In this setting, there is indeed a version of the ``cotangent functor'', where $T^*$ assigns to a smooth map $G \to M$ a certain submanifold of $T^*G \times T^*M$. In fact, the submanifolds produced have an additional property with respect to the standard symplectic structures on cotangent bundles: they are actually \emph{lagrangian submanifolds} of $\overline{T^*G} \times T^*M$, where the bar denotes the same manifold but with the sign of the symplectic form switched. Such a lagrangian submanifold will be called a \emph{canonical relation}, and can be thought of as a ``generalized symplectic map'' from $T^*G$ to $T^*M$.

The following picture thus emerges. We construct a ``symplectic category'' whose objects are symplectic manifolds and whose morphisms are given by canonical relations; in particular, symplectic maps (such as symplectomorphisms) between symplectic manifolds always produce such canonical relations via their graphs. Thus, we view the ``symplectic category'' as an enlargement of a more standard category of symplectic manifolds whose morphisms are symplectic maps. This is a point of view which has its roots in the works of various people: H\"ormander, Duistermaat, Guillemin, Sternberg, and Weinstein to name a few. As is well-known, this picture is incomplete since, while it is easy to define the composition of relations set-theoretically, the resulting objects in the smooth setting may no longer be smooth themselves. We will discuss such issues in the sequel, but for the time being we naively press ahead as if such a composition were well-defined.

Returning to our attempted construction of cotangent stacks, enlarging the target category of $T^*$ to include not simply smooth maps but more generally smooth (canonical) relations gives a well-defined operation, and applying this functor to the groupoid $G \rightrightarrows M$ produces an object which we denote by $T^*G \rightrightarrows T^*M$. Understanding this resulting structure in the context of the symplectic category was the main motivation for this thesis. The naive hope is that $T^*G \rightrightarrows T^*M$ is itself a ``groupoid'' in the symplectic category, and would then be a ``presentation'' of the ``cotangent stack'' $T^*\X$; this idea was first proposed in a comment by Weinstein in \cite{W5}. We will attempt to make sense of this statement, and will see that when $\X$ is an orbifold (i.e. a proper, \'etale differentiable stack), this idea works well, but the general case is not so clear. 

Now, the resulting structure on $T^*G \rightrightarrows T^*M$, which we will call a \emph{symplectic hopfoid} structure, exploits the fact that the cotangent bundle of a groupoid actually carries additional structure, namely that of a \emph{symplectic double groupoid} \cite{M}. We will see that the relation between this structure and that on $T^*G \rightrightarrows T^*M$ generalizes to any symplectic double groupoid, obtaining a characterization of symplectic double groupoids in terms of the symplectic category. Thus, in a sense, the ``groupoid'' objects in the symplectic category are really the symplectic double groupoids. 

In general, it is not clear what the objects $T^*G \rightrightarrows T^*M$ have to do with cotangent stacks, but we will see evidence that there should be some kind of relation between these two concepts. First, as mentioned above, in the case where $\X$ is actually an orbifold, we will see that $T^*G \rightrightarrows T^*M$ indeed encodes the cotangent orbifold $T^*\X$. Second, this construction generalizes to the symplectic hopfoids arising from general symplectic double groupoids in that, under certain assumptions, these structures do encode the data necessary to define a \emph{symplectic orbifold}. Finding the correct relation in general between symplectic structures on arbitrary differentiable stacks and objects in the symplectic category is an interesting problem which we hope to return to in the future.

We now provide a brief outline of this thesis, summarizing the main results and introducing the various structures which we will use throughout.

\section{Groupoids and Stacks}
Lie groupoids arise as generalizations of both Lie groups and of manifolds themselves, and may also be thought of as ``generalized'' smooth equivalence relations. The shortest possible definition of a (set-theoretic) groupoid is the following: a \emph{groupoid} is a small category in which every morphism is invertible. Then, a \emph{Lie groupoid} would be such an object equipped with compatible smooth structures. This definition, however, is not the one which one usually works with, and instead one usually spells out this definition in more concrete terms, as we will do.

As a motivating example, consider the smooth action of a Lie group $G$ on a  smooth manifold $M$. Then the quotient $M/G$, while always naturally a topological space, need not be smooth nor even Hausdorff unless the action is free and proper. Rather than attempting to work with this ill-behaved quotient, we may instead work with a certain Lie groupoid built from this data---the so-called action groupoid. This object encodes the quotient $M/G$ but also ``remembers'' how this quotient is built, and can be viewed as a smooth model for the quotient.

It is easy to see that different such groupoids can all be viewed as describing the same quotient in the sense above. To get around this, one defines a certain kind of equivalence relation among groupoids---namely, \emph{Morita equivalence}---and works with equivalence classes instead; this leads to the notion of a \emph{differentiable stack}. Stacks were first considered in the algebraic setting by Grothendieck, and have a general definition in terms of so-called ``categories fibered in groupoids''. The foundations of this point of view in the smooth setting were laid out in \cite{BX}, and have proved to be extremely useful. However, we will adopt the point of view that such structures can be described solely via Morita equivalence of Lie groupoids, and will avoid the full power of the category-theoretic definition. In particular, Weinstein's work \cite{W5} on the notion of a volume form on a differentiable stack is carried out in these terms.

Apart from attempts to model ill-behaved quotients in a smooth way, Lie groupoids also arise naturally in Poisson geometry as the ``global'' objects giving rise to Poisson structures. In a broad sense, the relation here is similar to that between Lie groups and Lie algebras. Making this precise leads to the notion of a \emph{symplectic groupoid}; indeed, such structures first arose from such considerations as part of a program to quantize Poisson manifolds. The rough idea is to first lift the Poisson structure to the corresponding symplectic groupoid, then quantize that symplectic groupoid via some quantization scheme compatible with the groupoid structure, and finally to push the result back down to the Poisson structure. Symplectic groupoids will play a key role in the results of this thesis.

We review some basics of the theory of Lie and symplectic groupoids in Appendix~\ref{appendix}, which should be read first if necessary. In particular, we give an explicit description of the standard symplectic groupoid structure on the cotangent bundle of a groupoid---an example which will be key in the what follows. We also give a definition of ``differentiable stack'' which will be suitable for our purposes, and review some constructions involving orbifolds. All of the material in this appendix is well-known, and we include it in order to provide the necessary background.

\section{The Symplectic Category}
The idea of viewing canonical relations as the morphisms of a category in which to do symplectic geometry arose in \cite{W3} and \cite{W4}, with the linear case in particular having been explored by Guillemin and Sternberg in \cite{GS2}. The main motivation came from H\"ormander's work \cite{Hor} on Fourier integral operators (see also \cite{D}), where he shows how to quantize canonical relations (with extra data) to produce linear maps between Hilbert spaces. This point of view led Weinstein to formulate the so-called ``symplectic creed'': Everything is a lagrangian submanifold.

The idea behind this statement is that it should be a goal to phrase constructions in symplectic geometry, as much as possible, in terms of lagrangian submanifolds and canonical relations. Apart from not actually being a true category, the so-called symplectic category indeed allows one to express various constructions in symplectic geometry much more naturally and in terms of purely ``categorical'' concepts. In particular, the notions of a symplectic groupoid and of a Hamiltonian action have particularly nice characterizations in this language, and the results of this thesis should hopefully provide further support for this point of view.

We review the basics of the symplectic category in Chapter~\ref{chap:symp}, beginning with some general remarks about smooth relations. Most of this is well-known, but some material---in particular that concerning a description of Hamiltonian actions of symplectic groupoids in terms of the symplectic category---has, to the author's knowledge, never explicitly appeared in print. In addition, the description of symplectic reduction in this setting will be a key construction in the rest of the thesis.

We also present some new ideas concerning the notion of a \emph{simplicial object} in the symplectic category. The main trouble in attempting to rigorously define the notion of a groupoid object in this category is a lack of a fiber product of canonical relations in general; the use of simplicial objects is one possible way of dealing with this. Indeed, we will see that there are many natural examples of such simplicial structures.

\section{Double Groupoids}
Double (Lie) groupoids\footnote{The definition of a double Lie groupoid and related constructions are recalled in Chapter~\ref{chap:dbl-grpds}.} arise in various ways in differential geometry, and can be viewed as a special case of so-called ``double categories'' (roughly, a ``category in the category of categories''), mimicking the first definition of a groupoid given above. As in that case, it is better to spell out the entire data in full and work with that. Much of the theory of double Lie groupoids was developed in \cite{M2} and \cite{M3}.

Symplectic double groupoids, in particular, first arose in attempts \cite{LW} to integrate Poisson-Lie groups and later Lie bialgebroids. More recently, a relation between symplectic double groupoids and ``symplectic 2-groupoids'' has been described in \cite{MT}. A key construction will be that of the \emph{core} of a symplectic double groupoid, which is a groupoid encoding some of the data of the full double groupoid structure.

In Chapter~\ref{chap:dbl-grpds}, the main chapter of this thesis, we present a characterization of symplectic double groupoids in terms of the symplectic category, made possible by various ways in which we can realize the core of a symplectic double groupoid as a symplectic quotient of the total space. We call the resulting object in the symplectic category a \emph{symplectic hopfoid} (see Definition~\ref{defn:symp-hopf}), which should be viewed as the analog of a groupoid object. The main result is:

{
\renewcommand{\thethrm}{\ref{thrm:main-cor}}
\begin{thrm}
There is a $1$-$1$ correspondence (up to symplectomorphism) between symplectic double groupoids and pairs of dual symplectic hopfoids.
\end{thrm}
\addtocounter{thrm}{-1}
}

Without assuming the presence of any symplectic structure, we also derive as a byproduct a characterization of double Lie groupoids in terms of so-called \emph{Lie hopfoids}, which can be thought of as the ``groupoid'' objects in the category of smooth manifolds and smooth relations. A nice relation between such structures and symplectic hopfoids is given by:

{
\renewcommand{\thethrm}{\ref{thrm:lie-hopf}}
\begin{thrm}
Suppose that $D$ is a double Lie groupoid with core $C$ and endow $T^*D$ with the induced symplectic double groupoid structure of Theorem \ref{thrm:ctdbl} with core $T^*C$. Then the induced symplectic hopfoid structure on $T^*D \rightrightarrows T^*C$ is simply the one obtained by applying the cotangent functor to the Lie hopfoid $D \rightrightarrows C$.
\end{thrm}
\addtocounter{thrm}{-1}
}

In the course of deriving the characterization mentioned above, we also present a new proof of Mackenzie's result that the core of a symplectic double groupoid naturally inherits a symplectic groupoid structure. Indeed, once we have these various structures phrased in ``categorical'' and diagrammatic terms, this construction falls out quite naturally.

\section{Symplectic Orbifolds}
Orbifolds arise naturally in topology and geometry as a certain type of ``singular space''. The original definition, as for manifolds, was phrased in terms of gluing local models, where the local models are quotients of vector spaces by linear actions of finite groups. The ``stacky'' point of view of orbifolds is much more recent, and has many benefits over the original theory; for example, the notion of a morphism between orbifolds becomes much simpler to consider from the point of view of stacks (see \cite{L} and \cite{IM}).

As opposed to general differentiable stacks, orbifolds have well-defined notions of cotangent stacks and of symplectic structures. The key to these nice descriptions rests in the fact that orbifolds have finite isotropy, which leads to one being able to avoid having to deal with the full machinery of stacks.

In Chapter~\ref{chap:stacks} we show how the symplectic hopfoid structures of the previous chapter induce certain Lie groupoid structures, equipped with additional structure coming from the use of canonical relations. The main results are the following:\footnote{We use the notation introduced in Chapter~\ref{chap:dbl-grpds} for the structure maps of a symplectic double groupoid.}

{
\renewcommand{\thethrm}{\ref{thrm:ind}}
\begin{thrm}
Let $S$ be a symplectic double groupoid with core $C$, and suppose that the restriction $\wt \ell_P|_C: C \to P^*$ is transverse to $1^{P^*} M \subseteq P^*$. Then the symplectic hopfoid structure on $S \rightrightarrows C$ gives rise to a Lie groupoid structure on
\[
D := \wt r_P^{-1}(1^{P^*}M) \cap \wt \ell_P^{-1}(1^{P^*}M) \rightrightarrows Y
\]
where $Y$ is the preimage of $1^{P^*}M$ under the restriction $\wt\ell_P|_C$.
\end{thrm}
\addtocounter{thrm}{-1}
}

{
\renewcommand{\thethrm}{\ref{thrm:symp-orb}}
\begin{thrm}
Let $S \rightrightarrows C$ be the symplectic hopfoid corresponding to a symplectic double groupoid $S$, and let $D \rightrightarrows Y$ be its induced groupoid. Then $Y$ is a coisotropic submanifold of $C$ whose characteristic leaves are the orbits of the induced groupoid. Moreover, the restriction of the symplectic form on $C$ to $Y$ has the property that its pullback under the target and source of the induced groupoid agree.
\end{thrm}
\addtocounter{thrm}{-1}
}

The structure described by the theorem above is, in the orbifold case, precisely that of a symplectic orbifold; a construction of cotangent orbifolds will be a special case. We also examine how and when this construction generalizes to simplicial objects in the symplectic category. The goal of this point of view, which we leave as a conjecture in this thesis, will be to show that all symplectic orbifolds arise in this manner, thus providing a new tool for their study.

We also note here the curious observation of Example~\ref{ex:inertia-induced}, where, given a Lie groupoid $G \rightrightarrows M$, we construct its \emph{inertia groupoid} $G \ltimes_M LG \rightrightarrows LG$ via the methods developed in this thesis; understanding why we obtain this result should eventually lead to a deeper understanding of our constructions.

\section{Further Directions}
Finally, in Chapter~\ref{chap:further} we consider further directions. This chapter is mainly speculative, and touches upon some of the original motivations which produced this thesis.

In particular, the main results of Chapter~\ref{chap:stacks} quoted above suggest that symplectic hopfoids encode ``symplectic-like'' structures on differentiable stacks, with ``symplectic-like'' actually meaning symplectic in the orbifold case. We summarize some questions which would have to be answered to make sense of this in general, such as what the correct notion of ``morphism'' between symplectic hopfoids should be. We end with some comments about attempting to quantize such structures, and suggest that the corresponding ``quantum'' objects should be \emph{quantum groupoids}\footnote{See \cite{X2} for a discussion of quantum groupoids.}, or more generally \emph{Hopf algebroids}.

It remains to be seen precisely why the constructions of this thesis---in particular, those of Chapter~\ref{chap:dbl-grpds}---work as they do: that is, is there a deeper ``categorical'' interpretation of the structures produced? We can speculate that the answer should come from some general properties of double categories and possibly monoidal categories, but will not discuss this issue any further in this thesis. This also remains an important avenue of further consideration. It would also be nice to know whether our constructions have any relation to the ``symplectic 2-groupoids'' of \cite{MT}, and whether our constructions may be useful in providing new insight into the integration problem for Poisson groupoids.

\section*{Notation and Conventions}
All manifolds considered in this thesis will be smooth and, except for possibly the total spaces of Lie groupoids, Hausdorff. By \emph{submanifold} we mean embedded submanifold, and by \emph{immersed submanifold} we simply mean the image of an immersion. The category of smooth manifolds and smooth maps will be denoted by $\Man$. In addition, all Lie group actions will be smooth.

For a smooth map $f: M \to N$, $df: TM \to TN$ will denote the induced map on tangent bundles; $df_p^*: T^*_{f(p)}N \to T^*_pM$ is then the induced map on cotangent spaces. We also denote the pairing of vectors and covectors by angled brackets: $\langle \cdot,\cdot \rangle$. For a foliation $\F$ on a manifold $\F$, $T\F$ denotes the associated involutive tangent distribution, $N^*\F$ denotes the conormal bundle of $\F$:
\[
N^*\F := \{(p,\xi) \in T^*M\ |\ \xi \text{ annihilates } T_p\F\},
\]
and $M/\F$ denotes the leaf space of $\F$. Also, we will occasionally have the need to identify the cotangent bundle $T^*G$ of a Lie group $G$ with $G \times \g^*$, where $\g$ is the Lie algebra of $G$, using \emph{left translations}; by this we mean the identification
\begin{align*}
T^*G &\cong G \times \g^* \\
(g,\xi) &\mapsto (g,(dL_g)_1^*\xi) \\
(g,(dL_{g^{-1}})_g^*\eta) &\mapsfrom (g,\eta)
\end{align*}
where $L_h: G \to G$ is left multiplication by $h \in G$ and $1 \in G$ is the identity element.

Given a symplectic manifold $(S,\omega)$, $\overline S$ will denote the symplectic manifold $(S,-\omega)$, which we call the \emph{dual} of $S$. For a subspace $U \subset V$ of a symplectic vector space $V$, $U^\perp$ will denote its symplectic orthogonal; when $U$ is coisotropic in $V$, $U/U^\perp$ will be called the \emph{reduction} of $V$ by $U$---it is naturally a symplectic vector space itself. For a $2$-form $\omega$ on a manifold $X$, $\ker\omega$ denotes the (possibly singular) vector bundle:
\[
\ker\omega = \{(p,v) \in TX\ |\ \omega_p(v,u) = 0 \text{ for all } u \in T_pX\}.
\]
When $\omega$ is closed and $\ker\omega$ is regular---meaning that all fibers have the same dimension---$(X,\omega)$ will be called \emph{presymplectic}; in this case, $\ker\omega$ is an involutive distribution on $X$.
\chapter{The Symplectic Category}\label{chap:symp}

In this chapter, we provide the necessary background on the symplectic category. In particular, the characterizations of symplectic reduction and of symplectic groupoids given here will be used throughout this thesis. We also discuss some new ideas concerning fiber products and simplicial objects in the symplectic category.

\section{Motivation}
As a first attempt to construct a category whose objects are symplectic manifolds, one may think to allow only symplectomorphisms as morphisms. This gives a perfectly nice category but is too restrictive since, for example, there will be no morphisms between symplectic manifolds of different dimensions. Even allowing general symplectic maps (i.e. maps $f: (X,\omega_X) \to (Y,\omega_Y)$ such that $f^*\omega_Y = \omega_X$) does not give enough morphisms since this condition requires that $f$ be an immersion.

Rather, it turns out that a good notion of ``morphism'' is provided by the notion of a \emph{canonical relation}. We point out here two motivations for this. First, it is well-known that a diffeomorphism $f: X \to Y$ between symplectic manifolds $X$ and $Y$ is a symplectomorphism if and only if its graph is a lagrangian submanifold of $\overline X \times Y$. Hence, it makes sense to consider arbitrary lagrangian submanifolds of this product as a kind of ``generalized'' symplectomorphisms.

The second motivation comes from a quantization perspective. Under the ``quantization dictionary''---describing a ``functor'' $Q$ from a category of symplectic manifolds to a category of vector spaces\footnote{One usually requires that these actually be Hilbert spaces.}---symplectic manifolds should quantize to vector spaces and lagrangian submanifolds to vectors:
\[
Q: X \mapsto Q(X),\ L \subset X \mapsto v \in Q(X).
\]
Furthermore, under this correspondence duals should quantize to duals and products to tensor products:
\[
Q: \overline{X} \mapsto Q(X)^*,\ X \times Y \mapsto Q(X) \otimes Q(Y).
\]
Then considering a linear map $Q(X) \to Q(Y)$ as an element of $Q(X)^* \otimes Q(Y)$, we see that the object which should quantize to a linear map is a lagrangian submanifold of $\overline{X} \times Y$; this is what we will call a \emph{canonical relation}.

\section{Symplectic Linear Algebra}
Here we collect a few basic facts from symplectic linear algebra which will be useful. Throughout, we assume that $(U,\omega_U)$, $(V,\omega_V)$, and $(W,\omega_W)$ are symplectic vector spaces. As usual, a bar written over a symplectic vector space indicates the the sign of the symplectic form is switched, and $^\perp$ indicates the symplectic orthogonal of a subspace.

\begin{prop}\label{prop:dom}
Suppose that $L \subset \overline{V} \otimes W$ is a lagrangian subspace. Then the image $pr_V(L)$ of $L$ under the projection $pr_V: \overline{V} \oplus W \to \overline{V}$ to the first factor is a coisotropic subspace of $\overline{V}$, and hence of $V$ as well.
\end{prop}

\begin{proof}
Let $v \in pr_V(L)^\perp$. For any $u + w \in L$, we have
\[
-\omega_V \oplus \omega_W(v+0,u+w) = -\omega_V(v,u) + \omega_W(0,w) = 0
\]
since $u \in pr_V(L)$ and $v$ is in the symplectic orthogonal to this space. Thus since $L$ is lagrangian, $v + 0 \in L$. Hence $v \in pr_V(L)$ so $pr_L(V)$ is coisotropic as claimed.
\end{proof}

By exchanging the roles of $V$ and $W$ above and exchanging the components of $L$, we see that the same result holds for the image of $L$ under the projection to the second factor.

\begin{prop}\label{prop:comp}
Suppose that $L \subset \overline{U} \times V$ and $L' \subset \overline{V} \times W$ are lagrangian subspaces such that $L \oplus L'$ and $U \oplus \Delta_V \oplus W$ are transverse subspaces of $\overline{U} \oplus V \oplus \overline{V} \oplus W$, where $\Delta_V \subset V \oplus \overline{V}$ is the diagonal. Then the image $L' \circ L$ of
\[
L \times_V L' := (L \oplus L') \cap (U \oplus \Delta_V \oplus W)
\]
under the natural projection from $\overline{U} \oplus V \oplus \overline{V} \oplus W$ to $\overline{U} \oplus W$ is lagrangian.
\end{prop}

\begin{proof}
We first claim that the image $L' \circ L$ of the above projection is isotropic in $\overline{U} \oplus W$. Indeed, suppose that $(u,w), (u',w') \in L \circ L'$. Then there are $v,v' \in V$ such that
\[
(u,v) \in L, (v,w) \in L' \text{ and } (u',v') \in L, (v',w') \in L'.
\]
Since $L$ and $L'$ are lagrangian (and hence isotropic), we have
\[
0 = -\omega_U \oplus \omega_V((u,v),(u',v')) = -\omega_U(u,u') + \omega_V(v,v'),
\]
and
\[
0 = -\omega_V \oplus \omega_W((v,w),(v',w')) = -\omega_V(v,v') + \omega_W(w,w').
\]
Thus
\[
-\omega_U \oplus \omega_W((u,w),(u',w')) = -\omega_U(u,u') + \omega_W(w,w') = -\omega_V(v,v') + \omega_V(v,v') = 0,
\]
showing that $L' \circ L$ is isotropic as claimed.

Next, a simple dimension count using transversality shows that $L \times_V L'$ has half the dimension of $\overline U \oplus W$. And finally, if $(0,w,w,0) \in L \times_V L'$, then it is to see that $(0,w,w,0)$ annihilates $U \oplus \Delta_V \oplus W$ and $L \times L'$ with respect to the symplectic form
\[
-\omega_U \oplus \omega_V \oplus (-\omega_V) \oplus \omega_W,
\]
and hence annihilates all of $\overline{U} \oplus V \oplus \overline{V} \oplus W$ by transversality. This implies that $w = 0$ by the non-degeneracy of the above symplectic structure, so we conclude that the projection of $L \times_V L'$ onto $L' \circ L$ is injective. Then $L' \circ L$ has half the dimension of $\overline{U} \oplus W$ as well and is thus lagrangian.
\end{proof}

\begin{prop}\label{prop:lin-factor}
Suppose that $L \subset \overline{V} \oplus W$ is a lagrangian subspace, so that $C := pr_V(L)$ and $Y := pr_W(L)$ are coisotropic subspaces of $V$ and $W$ respectively. Then the symplectic vector spaces $C/C^\perp$ and $Y/Y^\perp$ are naturally symplectomorphic.
\end{prop}

\begin{proof}
We define a map $T: C \to Y/Y^\perp$ as follows: for $v \in C$, choose $w \in Y$ such that $(v,w) \in L$ and set $T(c) := [w]$. To see that this is well-defined, suppose that $(v,w), (v,w') \in L$ and let $z \in Y$. Choose $x \in C$ such that $(x,z) \in L$. Then
\[
0 = -\omega_V \oplus \omega_W((v,w),(x,z)) = -\omega_V(v,x) + \omega_W(w,z),
\]
and similarly $\omega_V(v,x) = \omega_W(w',z)$. Thus
\[
\omega_W(w-w',z) = \omega_W(w,z) - \omega_W(w',z) = \omega_V(v,x) - \omega_V(v,x) = 0.
\]
Thus $w-w' \in Y^\perp$, so $[w] = [w']$ and hence $T$ is well-defined; it is clearly linear.

The equation
\[
-\omega_V(v,v') + \omega_W(w,w') = 0 \text{ for } (v,w), (v,',w') \in L
\]
easily implies that $C^\perp = \ker T$, so we get an induced isomorphism
\[
T: C/C^\perp \to Y/Y^\perp.
\]
It is then straightforward to check that this is a symplectomorphism using the fact that $L \subset \overline{V} \oplus W$ is lagrangian.,
\end{proof}

\section{Smooth Relations}
Before moving on to canonical relations, we first consider the general theory of smooth relations.

\begin{defn}
A \emph{smooth relation} $R$ from a smooth manifold $M$ to a smooth manifold $N$ is a closed submanifold of the product $M \times N$. We will use the notation $R: M \to N$ to mean that $R$ is a smooth relation from $M$ to $N$. We will suggestively use the notation $R: m \mapsto n$ to mean that $(m,n) \in R$, and will think of relations as partially-defined, multi-valued functions. The \emph{transpose} of a smooth relation $R: M \to N$ is the smooth relation $R^t: N \to M$ defined by $(n,m) \in R^t$ if $(m,n) \in R$.
\end{defn}

\begin{rmk}
To make compositions easier to read, it may be better to draw arrows in the opposite direction and say that $R: N \leftarrow M$ is a smooth relation \emph{to} $N$ \emph{from} $M$; this is the approach used by Weinstein in \cite{W1}. We will stick with the notation above.
\end{rmk}

Given a smooth relation $R: M \to N$, the \emph{domain} of $R$ is
\[
\dom R := \{m \in M\ |\ \text{there exists $n \in N$ such that $(m,n) \in R$}\} \subseteq M.
\]
In other words, this is the domain of $R$ when we think of $R$ as a partially-defined, multi-valued function. Given a subset $U \subseteq M$, its \emph{image} under the smooth relation $R: M \to N$ is the set
\[
R(U) := \{n \in N\ |\ \text{there exists $m \in M$ such that $(m,n) \in R$}\} \subseteq N.
\]
In particular, we can speak of the image $R(m)$ of a point $m \in M$ and of the image $\im R := R(M)$ of $R$. The domain of $R$ can then be written as $\dom R = R^t(N) = \im R^t$.

We may now attempt to define a composition of smooth relations simply by using the usual composition of relations: given a smooth relation $R: M \to N$ and a smooth relation $R': N \to Q$, the composition $R' \circ R: M \to Q$ is defined to be
\[
R' \circ R := \{ (m,q) \in M \times Q\ |\ \text{there exists } n \in N \text{ such that } (m,n) \in R \text{ and } (n,q) \in R'\}.
\]
This is the same as taking the intersection of $R \times R'$ and $M \times \Delta_N \times Q$ in $M \times N \times N \times Q$, where $\Delta_N$ denotes the diagonal in $N \times N$, and projecting to $M \times Q$. However, we immediately run into the problem that the above composition need no longer be a smooth closed submanifold of $M \times Q$, either because the intersection of $R \times R'$ and $M \times \Delta_N \times Q$ is not smooth or because the projection to $M \times Q$ is ill-behaved, or both. To fix this, we introduce the following notions:

\begin{defn}
A pair $(R,R')$ of smooth relations $R: M \to N$ and $R': N \to Q$ is \emph{transversal} if the submanifolds $R \times R'$ and $M \times \Delta_N \times Q$ intersect transversally. The pair $(R,R')$ is \emph{strongly transversal} if it is transversal and in addition the projection of
\[
(R \times R') \cap (M \times \Delta_N \times Q)
\]
to $M \times Q$ is a proper embedding.
\end{defn}

As a consequence, for a strongly transversal pair $(R,R')$, the composition $R' \circ R$ is indeed a smooth relation from $M$ to $Q$.

\begin{defn}
A relation $R: M \to N$ is said to be:
\begin{itemize}
\item \emph{surjective} if for any $n \in N$ there exists $m \in M$ such that $(m,n) \in R$,
\item \emph{injective} if whenever $(m,n), (m',n) \in R$ we have $m=m'$,
\item \emph{cosurjective} if for any $m \in M$ there exists $n \in N$ such that $(m,n) \in R$,
\item \emph{coinjective} if whenever $(m,n), (m,n') \in R$ we have $n=n'$.
\end{itemize}
\end{defn}

Note that, if we think of relations as partially defined multi-valued functions, then cosurjective means ``everywhere defined'' and coinjective means ``single-valued''. Note also that $R$ is cosurjective if and only if $R^t$ is surjective and $R$ is coinjective if and only if $R^t$ is injective.

\begin{defn}
A smooth relation $R: M \to N$ is said to be a \emph{surmersion} if it is surjective and coinjective, the projection of $R$ to $M$ is a proper embedding, and the projection of $R$ to $N$ is a submersion; it is a \emph{cosurmersion} if $R^t: N \to M$ is a surmersion.
\end{defn}

It is a straightforward check to see that $R$ is a surmersion if and only if $R \circ R^t = id$ and hence a cosurmersion if and only if $R^t \circ R = id$. It is also straightforward to check that a pair $(R,R')$ is always strongly transversal if either $R$ is a surmersion or $R'$ a cosurmersion.

\section{Canonical Relations}
\begin{defn}
A \emph{canonical relation} $L: P \to Q$ from a symplectic manifold $P$ to a symplectic manifold $Q$ is a smooth relation which is lagrangian as a submanifold of $\overline{P} \times Q$.
\end{defn}

\begin{rmk}
The term ``lagrangian relation'' would perhaps be a better choice of words, and would fit in better with the terminology in \cite{WW}. Other sources use ``symplectic relation'' instead. The term ``canonical relation'' is motivated by Hamiltonian mechanics, where symplectomorphisms were classically called ``canonical transformations''.
\end{rmk}

\begin{ex}
As mentioned before, the graph of a symplectomorphism $f: P \to Q$ is a canonical relation $P \to Q$, which by abuse of notation we will also denote by $f$. In particular, given any symplectic manifold $P$, the graph of the identity map will be the canonical relation $id: P \to P$ given by the diagonal in $\overline P \times P$. More generally, the graph of a symplectic \'etale map is a canonical relation.
\end{ex}

\begin{ex}
For any manifold $M$, the \emph{Schwartz transform} on $T^*M$ is the canonical relation
\[
s: T^*M \to \overline{T^*M},\ (p,\xi) \mapsto (p,-\xi)
\]
given by multiplication by $-1$ in the fibers. Alternatively, it is the lagrangian submanifold of $T^*M \times T^*M \cong T^*(M \times M)$ given by the conormal bundle to the diagonal of $M \times M$.
\end{ex}

\begin{ex}
For any symplectic manifold $S$, a canonical relation $pt \to S$ or $S \to pt$ is nothing but a closed lagrangian submanifold of $S$.
\end{ex}

Here is a basic fact, which follows from Proposition \ref{prop:dom}:

\begin{prop}
Suppose that $L: P \to Q$ is a canonical relation with the projection of $L$ to $P$ of constant rank. Then the domain of $L$ is a coisotropic submanifold of $P$.
\end{prop}

The same result also shows---replacing $L$ by $L^t$---that the image of a canonical relation is coisotropic when the projection to $Q$ is a submersion.

We also have as a consequence of Proposition \ref{prop:comp} and of the previous discussion on composing smooth relations:

\begin{prop}
If $L: X \to Y$ and $L': Y \to Z$ are canonical relations with $(L,L')$ strongly transversal, then $L' \circ L$ is a canonical relation.
\end{prop}

In other words, the only obstacle to the composition of canonical relations being well-defined comes from smoothness issues and not from the requirement that the resulting submanifold be lagrangian. Later, we will discuss a way of getting around these smoothness complications.

\begin{rmk}
The composition of canonical relations is well-defined under weaker assumptions than strong transversality. In particular, assuming only transversality, Proposition \ref{prop:comp} implies that the composition will be an \emph{immersed} lagrangian submanifold; this also holds under a weaker \emph{clean intersection} hypothesis. In this thesis, however, we will only need to consider strong transversality.
\end{rmk}

\begin{rmk}
Note that, for general smooth relations $R: X \to Y$ and $R': Y \to Z$, transversality of $(R,R')$ does not necessarily imply that $R' \circ R$ is immersed---this is a key difference between smooth and canonical relations.
\end{rmk}

We note here that when $U \subset X$ is a lagrangian submanifold viewed as a canonical relation $pt \to X$ and the pair $(U,L)$ is strongly transversal, then the image $L(U)$ is just the composition
\begin{center}
\begin{tikzpicture}[>=angle 90]
	\node (1) at (0,1) {$pt$};
	\node (2) at (2,1) {$X$};
	\node (3) at (4,1) {$Y$};

	\tikzset{font=\scriptsize};
	\draw[->] (1) to node [above] {$U$} (2);
	\draw[->] (2) to node [above] {$L$} (3);
\end{tikzpicture}
\end{center}
and is a lagrangian submanifold of $Y$. In this way, a canonical relation $L: X \to Y$ induces a map from (a subset of) the set of lagrangian submanifolds of $X$ to the set of lagrangian submanifolds of $Y$.

\begin{defn}
A canonical relation $L: X \to Y$ is said to be a \emph{reduction} if, as a smooth relation, it is a surmersion; it is a \emph{coreduction} if it is a cosurmersion. We use $L: X \red Y$ to denote that $L$ is a reduction, and $L: X \cored Y$ to denote that $L$ is a coreduction.
\end{defn}

The use of the term ``reduction'' is motivated by the following example.

\begin{ex}(Symplectic Reduction)
Let $(M,\omega)$ be a symplectic manifold and $C$ a coisotropic submanifold. The distribution on $C$ given by $\ker\omega \subset TC$ is called the \emph{characteristic distribution} of $C$. It follows from $\omega$ being closed that $\ker\omega$ is integrable; the induced foliation $C^\perp$ on $C$ will be called the \emph{characteristic foliation} of $C$. If the leaf space $C/C^\perp$ is smooth and the projection $C \to C/C^\perp$ is a submersion, then $C/C^\perp$ is naturally symplectic and the relation
\[
red: M \red C/C^\perp
\]
assigning to an element of $C$ the leaf which contains it is a canonical relation which is a reduction in the sense above. This process will be called \emph{symplectic reduction}.

Symplectic reduction via Hamiltonian actions of Lie groups is a special case of the construction above. Indeed, suppose that $G$ acts properly and freely on $(P,\omega)$ in a Hamiltonian way with equivariant momentum map $\mu: P \to \g^*$. Then $C := \mu^{-1}(0)$ is a coisotropic submanifold of $P$ and the orbits of the induced $G$ action on $\mu^{-1}(0)$ are precisely the leaves of the characteristic foliation. Hence the reduction $C/C^\perp$ in the above sense is the quotient $\mu^{-1}(0)/G$.
\end{ex}

\begin{ex}
We also note two more well-known examples of symplectic reduction.\footnote{See, for example, \cite{BW}.} Suppose that $X$ is a smooth manifold with $Y \subseteq X$ a submanifold. Then the restricted cotangent bundle $T^*X|_Y$ is a coisotropic submanifold of $T^*X$ whose reduction is symplectomorhpic to $T^*Y$. Thus we obtain a reduction $T^*X \twoheadrightarrow T^*Y$.

As a generalization, suppose now that $\F$ is a regular foliation on $Y \subseteq X$ with smooth, Hausdorff leaf space $Y/\F$. Then the conormal bundle $N^*\F$ is a coisotropic submanifold of $T^*X$ (this is in fact equivalent to the distribution $T\F$ being integrable) and its reduction is canonically symplectomorphic to $T^*(Y/\F)$, giving rise to a reduction relation $T^*X \twoheadrightarrow T^*(Y/\F)$. The previous example is the case where $\F$ is the zero-dimensional foliation given by the points of $Y$.
\end{ex}

Any reduction can be expressed as the composition of a symplectic reduction as above with (the graph of) a symplectic \'etale map. More generally, it follows from Proposition~\ref{prop:lin-factor} that any canonical relation $\Lambda: X \to Y$ can be factored (modulo constant rank issues) into the composition of a reduction, followed by a symplectic \'etale map, followed by a coreduction:
\begin{equation}\label{factor}
\begin{tikzpicture}[>=angle 90,baseline=(current  bounding  box.center)]
	\node (U1) at (0,1) {$X$};
	\node (U2) at (3,1) {$Y$};
	\node (L1) at (0,-1) {$X_{\dom\Lambda}$};
	\node (L2) at (3,-1) {$Y_{\im\Lambda}$};

	\tikzset{font=\scriptsize};
	\draw[->] (U1) to node [above] {$\Lambda$} (U2);
	\draw[->>] (U1) to node [left] {$red$} (L1);
	\draw[>->] (L2) to node [right] {$cored$} (U2);
	\draw[->] (L1) to node [above] {$\Lambda$} (L2);
\end{tikzpicture}
\end{equation}
where $X_{\dom\Lambda}$ denotes the reduction of $X$ by the coisotropic $\dom L$ and $Y_{\im\Lambda}$ the reduction of $Y$ by the coisotropic $\im\Lambda$.

\section{The Symplectic Category, Version I}
We are now ready to introduce the category we will be working in. We start with a preliminary definition, which we call a category even though, as explained above, compositions are not always defined.

\begin{defn}
The \emph{symplectic category} is the category $\Symp$ whose objects are symplectic manifolds and whose morphisms are canonical relations.
\end{defn}

Apart from the existence of compositions, it is easy to see that this indeed satisfies the other conditions of being a category. In particular identity morphisms are provided by diagonals and the associativity of composition follows from the usual identification
\[
X \times (Y \times Z) = (X \times Y) \times Z.
\]

To illustrate a use of the symplectic category, we prove the following:

\begin{prop}\label{prop:act}
Let $G$ be a Lie group acting smoothly on a symplectic manifold $P$ and let $\mu: P \to \g^*$ be a smooth, equivariant map with respect to the given action on $P$ and the coadjoint action on $\g^*$. Then the action of $G$ on $P$ is Hamiltonian with moment map $\mu$ if and only if the relation $T^*G \times P \to P$ given by
\[
((g,\mu(p)),p) \mapsto gp,
\]
where we have identified $T^*G$ with $G \times \g^*$ using left translations, is a canonical relation.
\end{prop}

\begin{proof}
First, it is simple to see that the submanifold $L := \{((g,\mu(p)),p,gp)\}$ has half the dimension of the ambient space $\overline{T^*G} \times \overline{P} \times P$ and that a tangent vector to $L$ looks like
\[
(X, d\mu_p V, V, dp_g X + dg_p V) \text{ where } X \in T_gG, V \in T_pP
\]
under the identification $T^*G \cong G \times \g^*$, and where (abusing notation) $p: G \to P$ and $g: P \to P$ respectively denote the maps
\[
h \mapsto hp \text{ and } q \mapsto gq.
\]
Then, the symplectic form on $\overline{T^*G} \times \overline{P} \times P$ acts as follows:
\begin{align*}
&-\omega_{T^*G} \oplus -\omega_P \oplus \omega_P((X, d\mu_p V, V, dp_g X + dg_p V),(X', d\mu_p V', V', dp_g X' + dg_p V')) \\
&\quad = -\omega_{T^*G}((X, d\mu_p V),(X',d\mu_p V')) - \omega_P(V,V') + \omega_P(dp_g X, dp_g X') \\
&\quad \quad + \omega_P(dp_g X, dg_p V') + \omega_P(dg_p V, dp_g X') + \omega_P(dg_p V, dg_p V') \\
&\quad = -\langle (dL_{g^{-1}})_gX, d\mu_p V' \rangle + \langle (dL_{g^{-1}})_gX', d\mu_p V \rangle - \langle (dL_{g^{-1}})_gX, d\mu'_p X' \rangle - \omega_P(V,V') \\ 
&\quad \quad + \omega_P(dp_g X, dp_g) + \omega_P(dp_g X, dg_p V') + \omega_P(dg_p V, dp_g X') + \omega_P(dg_p V, dg_p V'),
\end{align*}
where $L_{g^{-1}}: G \to G$ is left multiplication by $g^{-1}$ and $\mu': G \to \g^*$ is the map $h \mapsto h\cdot\mu(p)$ induced by the coadjoint action of $G$ on $\g^*$.

Now, if the action is Hamiltonian, then it is in particular symplectic so
\[
\omega_P(dg_p V, dg_p V') = \omega_P(V,V').
\]
Thus the above reduces to
\begin{align*}
&[-\langle (dL_{g^{-1}})_gX, d\mu_p V' \rangle + \omega_P(dp_g X, dg_p V')] + [\langle (dL_{g^{-1}})_gX', d\mu_p V \rangle - \omega_P(dp_g X', dg_p V)] \\
& \qquad\qquad\qquad\qquad\quad + [-\langle (dL_{g^{-1}})_gX, d\mu'_p X' \rangle + \omega_p(dp_g X, dp_g X')],
\end{align*}
and each of these vanishes by the moment map condition, where for the third term we use the fact that $\mu' = \mu \circ p$, which follows from the equivariance of $\mu$. Hence $L$ is lagrangian and so gives a canonical relation.

Conversely, suppose that $L$ is lagrangian. Then the above computation should produce zero, and setting $X' = 0$ and $V = 0$ gives the requirement that
\[
-\langle (dL_{g^{-1}})_gX, d\mu_p V' \rangle + \omega_P(dp_g X, dg_p V') = 0.
\]
This is precisely the moment map condition, and hence we conclude that the $G$-action is Hamiltonian with moment map $\mu$.
\end{proof}

The category $\Symp$ thus defined has additional rich structure. In particular, it is a \emph{monoidal category}, where the tensor operation is given by Cartesian product and the unit is given by the symplectic manifold consisting of a single point $pt$. Moreover, $\Symp$ is \emph{symmetric monoidal} and \emph{rigid}, where the dualizing operation is given by $X \mapsto \overline X$ on objects and $L \mapsto L^t$ on morphisms. In addition, if we allow the empty set as a symplectic manifold, then it is simple to check that $\emptyset$ is both an initial and terminal object in this category, and that the categorical product of symplectic manifolds $X$ and $Y$ is the disjoint union $X \sqcup Y$.

\section{The Symplectic Category, Version II}
Now we describe a method for getting around the lack of a well-defined composition of canonical relations in general, by using what are in a sense ``formal'' canonical relations:

\begin{defn}
A \emph{generalized lagrangian correspondence} from $X_0$ to $X_n$ is a chain
\begin{center}
\begin{tikzpicture}[>=angle 90]
	\node (1) at (0,1) {$X_0$};
	\node (2) at (2,1) {$X_1$};
	\node (3) at (4,1) {$\cdots$};
	\node (4) at (6,1) {$X_{n-1}$};
	\node (5) at (8,1) {$X_n$};

	\tikzset{font=\scriptsize};
	\draw[->] (1) to node [above] {$L_0$} (2);
	\draw[->] (2) to node [above] {$L_1$} (3);
	\draw[->] (3) to node [above] {$L_{n-2}$} (4);
	\draw[->] (4) to node [above] {$L_{n-1}$} (5);
\end{tikzpicture}
\end{center}
of canonical relations between intermediate symplectic manifolds. The composition of generalized lagrangian correspondences is given simply by concatenation. The \emph{Wehrheim-Woodward (symplectic) category} is the category whose morphisms are equivalence classes of generalized lagrangian correspondence under the equivalence relation generated by the requirement that
\begin{center}
\begin{tikzpicture}[>=angle 90]
	\node (1) at (0,1) {$X$};
	\node (2) at (2,1) {$Y$};
	\node (3) at (4,1) {$Z$};

	\tikzset{font=\scriptsize};
	\draw[->] (1) to node [above] {$L_0$} (2);
	\draw[->] (2) to node [above] {$L_1$} (3);
\end{tikzpicture}
\end{center}
be identified with $L_1 \circ L_0: X \to Z$ whenever the pair $(L_0,L_1)$ is strongly transversal.
\end{defn}

\begin{rmk}
Note that the same construction makes sense for general smooth relations.
\end{rmk}

The Wehrheim-Woodward category is an honest category, used in \cite{WW} to define \emph{quilted Floer homology}. The following result of Weinstein simplifies the types of correspondences which need to be considered when working with this category:

\begin{thrm}[Weinstein, \cite{W1}]
Any generalized lagrangian correspondence
\[
X_0 \to \cdots \to X_n
\]
is equivalent to a two-term lagrangian correspondence $X_0 \cored Y \red X_n$, where the first relation can be taken to be a coreduction and the second a reduction.
\end{thrm}

All compositions that we consider in this thesis will actually be strongly transversal, so that we need not use the full language of generalized lagrangian correspondences. To be precise, we will technically work in the Wehrheim-Woodward category but will only consider single-term correspondences. It would be interesting to know how and if our results generalize to the full Wehrheim-Woodward category.

\section{The Cotangent Functor}

We define a functor $T^*: \Man \to \Symp$, called the \emph{cotangent functor}, as follows. First, $T^*$ assigns to a smooth manifold its cotangent bundle. To a smooth map $f: M \to N$, $T^*$ assigns the canonical relation $T^*f: T^*M \to T^*N$ given by
\[
T^*f: (p,df_p^*\xi) \mapsto (f(p),\xi).
\]
This is nothing but the composition $T^*M \to \overline{T^*M} \to T^*N$ of the Schwartz transform of $T^*M$ followed by the canonical relation given by the conormal bundle to the graph of $f$ in $M \times N$. We call $T^*f$ the \emph{cotangent lift} of $f$. It is a simple check to see that pairs of cotangent lifts are always strongly transversal and that $T^*$ really is then a functor: i.e. $T^*(f \circ g) = T^*f \circ T^*g$ and $T^*(id) = id$. Note also that the same construction makes sense even when $f$ is only a smooth relation. 

\begin{ex}
When $\phi: M \to N$ is a diffeomorphism, $T^*\phi: T^*M \to T^*N$ is precisely the graph of the lifted symplectomorphism.
\end{ex}

The following is easy to verify.

\begin{prop}
The cotangent lift of $f$ is a reduction if and only if $f$ is a surmersion; the cotangent lift of $g$ is a coreduction if and only if $g$ is cosurmersion.
\end{prop}

Much of this thesis was motivated by the process of applying this functor to various structures arising in $\Man$ and studying the resulting structures in $\Symp$; in particular, in Chapter~\ref{chap:dbl-grpds} we describe the type of structure arising from applying $T^*$ to a Lie groupoid. For now, let us compute the cotangent lift of a group action, which motivates the canonical relation used in Proposition \ref{prop:act}:

\begin{prop}{\label{prop:action-lift}}
Let $\tau: G \times M \to M$ be a smooth action of a Lie group $G$ on a smooth manifold $M$. Then the cotangent lift $T^*\tau$ is given by
\[
((g,\mu(p,\xi)),(p,\xi)) \mapsto g(p,\xi),
\]
where $g(p,\xi)$ denotes the lifted cotangent action of $G$ on $T^*M$, $\mu$ is its standard momentum map: $\mu(p,\xi) = dp_e^*\xi$, and we have identified $T^*G$ with $G \times \g^*$ using left translations.
\end{prop}

\begin{proof}
The cotangent lift $T^*\tau: T^*G \times T^*M \to T^*M$ is given by
\[
((g,p),d\tau_{(g,p)}^*\xi) \mapsto (gp,\xi).
\]
Abusing notation, for any $g \in G$ and $p \in P$ we denote the maps
\[
P \to P,\ q \mapsto gq \text{ and } G \to P,\ h \mapsto hp
\]
by $g$ and $p$ respectively. Then using the fact that
\[
d\tau_{(g,p)} = dp_g \circ pr_1 + dg_p \circ pr_2,
\]
we can write $T^*\tau$ as
\[
((g,dp_g^*\xi),(p,dg_p^*\xi)) \mapsto (gp,\xi),
\]
which is then
\[
((g,dp_g^*(dg^{-1}_{gp})^*\eta,(p,\eta)) \mapsto g(p,\eta)
\]
where $\eta = dg_p^*\xi$ and $g(p,\eta)$ now denotes the lifted cotangent action of $G$ on $T^*M$.

Now, letting $L_g$ denote left multiplication on $G$ by $g$, we have
\begin{align*}
(dL_g)_e^*dp_g^*(dg^{-1}_{gp})^*\eta &= d(g^{-1} \circ p \circ L_g)_e^*\eta \\
&= dp_e^*\eta.
\end{align*}
Thus identifying $T^*G$ with $G \times \g^*$ using left translations:
\[
(g,\gamma) \mapsto (g,(dL_g)_e^*\gamma),
\]
we get that $T^*\tau$ is given by the desired expression.
\end{proof}

\section{Symplectic Monoids and Comonoids}
We can now use the symplectic category to provide simple, ``categorical'' descriptions of various objects encountered in symplectic geometry; in particular, symplectic groupoids and their Hamiltonian actions.

Since $\Symp$ is monoidal, we can speak about \emph{monoid objects} in $\Symp$:

\begin{defn}
A \emph{symplectic monoid} is a monoid object in $\Symp$. Thus, a symplectic monoid is a triple $(S,m,e)$ consisting of a symplectic manifold $S$ together with canonical relations
\[
m: S \times S \to S \text{ and } e: pt \to S,
\]
called the \emph{product} and \emph{unit} respectively, so that
\begin{center}
\begin{tikzpicture}[>=angle 90]
	\node (UL) at (0,1) {$S \times S \times S$};
	\node (UR) at (3,1) {$S \times S$};
	\node (LL) at (0,-1) {$S \times S$};
	\node (LR) at (3,-1) {$P$};
	
	\tikzset{font=\scriptsize};
	\draw[->] (UL) to node [above] {$id \times m$} (UR);
	\draw[->] (UL) to node [left] {$m \times id$} (LL);
	\draw[->] (UR) to node [right] {$m$} (LR);
	\draw[->] (LL) to node [above] {$m$} (LR);
\end{tikzpicture}
\end{center}
and
\begin{center}
\begin{tikzpicture}[thick]
	\node (UL) at (0,1) {$S$};
	\node (UR) at (3,1) {$S \times S$};
	\node (URR) at (6,1) {$S$};
	\node (LR) at (3,-1) {$S$};

	\tikzset{font=\scriptsize};	
	\draw[->] (UL) to node [above] {$e \times id$} (UR);
	\draw[->] (UR) to node [right] {$m$} (LR);
	\draw[->] (URR) to node [above] {$id \times e$} (UR);
	\draw[->] (UL) to node [above] {$id$} (LR);
	\draw[->] (URR) to node [above] {$id$} (LR);
\end{tikzpicture}
\end{center}
commute. We also require that all compositions involved be strongly transversal. The first diagram says that $m$ is ``associative'' and the second says that $e$ is a ``left'' and ``right unit''. We often refer to $S$ itself as a symplectic monoid, and use subscripts in the notation for the structure morphisms if we need to be explicit.
\end{defn}

The main example of such a structure is the following:

\begin{ex}
Let $S \rightrightarrows P$ be a symplectic groupoid. Then $S$ together with the groupoid multiplication thought of as a relation $S \times S \to S$ and the canonical relation $pt \to S$ given by the image of the unit embedding $P \to S$ is a symplectic monoid.
\end{ex}

In fact, Zakrzewski gave in \cite{SZ1}, \cite{SZ2} a complete characterization of symplectic groupoids in terms of such structures. Let us recall his description. First, the base space of the groupoid is the lagrangian submanifold $E$ of $S$ giving the unit morphism
\[
e: pt \to S,\ E := e(pt).
\]
The associativity of $m$ and unit properties of $e$ together then imply that there are unique maps $\ell, r: S \to E$ such that
\[
m(\ell(s),s) \ne \emptyset \ne m(s,r(s)) \text{ for all $s$}.
\]
These maps will form the target and source maps of the sought after symplectic groupoid structure.

The above is not enough to recover the symplectic groupoid yet; in particular, the above conditions do not imply that the product $m$ must be single-valued as we would need for a groupoid product. We need one more piece of data and an additional assumption:

\begin{defn}
A \emph{*-structure} on a symplectic monoid $S$ is an anti-symplectomorphism $s: S \to S$ (equivalently a symplectomorphism $s: \overline S \to S$) such that $s^2 = id$ and the diagram
\begin{center}
\begin{tikzpicture}[>=angle 90]
	\node (UL) at (0,1) {$\overline{S} \times \overline{S}$};
	\node (U) at (3,1) {$\overline{S} \times \overline{S}$};
	\node (UR) at (6,1) {$S \times S$};
	\node (LL) at (0,-1) {$\overline{S}$};
	\node (LR) at (6,-1) {$S,$};

	\tikzset{font=\scriptsize};
	\draw[->] (UL) to node [above] {$\sigma$} (U);
	\draw[->] (U) to node [above] {$s \times s$} (UR);
	\draw[->] (UL) to node [left] {$m$} (LL);
	\draw[->] (UR) to node [right] {$m$} (LR);
	\draw[->] (LL) to node [above] {$s$} (LR);
\end{tikzpicture}
\end{center}
where $\sigma$ is the symplectomorphism exchanging components, commutes. A symplectic monoid equipped with a $*$-structure will be called a \emph{symplectic $*$-monoid}.

A $*$-structure $s$ is said to be \emph{strongly positive}\footnote{This terminology is motivated by quantization, where it is viewed as the analog of the algebraic positivity condition: ``$aa^* > 0$''.} if the diagram
\begin{center}
\begin{tikzpicture}[>=angle 90]
	\node (UL) at (0,1) {$S \times \overline{S}$};
	\node (UR) at (3,1) {$S \times S$};
	\node (LL) at (0,-1) {$pt$};
	\node (LR) at (3,-1) {$S,$};

	\tikzset{font=\scriptsize};
	\draw[->] (UL) to node [above] {$id \times s$} (UR);
	\draw[->] (LL) to node [left] {$$} (UL);
	\draw[->] (UR) to node [right] {$m$} (LR);
	\draw[->] (LL) to node [above] {$e$} (LR);
\end{tikzpicture}
\end{center}
where $pt \to S \times \overline{S}$ is the morphism given by the diagonal of $S \times \overline{S}$, commutes.
\end{defn}

\begin{thrm}[Zakrzewski, \cite{SZ1}\cite{SZ2}]
Symplectic groupoids are in $1$-$1$ correspondence with strongly positive symplectic $*$-monoids.
\end{thrm}

\begin{rmk}
There is a similar characterization of Lie groupoids as strongly positive $*$-monoids in the category of smooth relations.
\end{rmk}

It is unclear precisely what general symplectic monoids correspond to; in particular, as mentioned before, the product is then not necessarily single-valued. 

We also note that a strongly positive $*$-structure which produces out of a symplectic monoid $(S,m,e)$ a symplectic groupoid is in fact unique if it exists: the lagrangian submanifold of $S \times S$ which gives this $*$-structure must equal the composition
\begin{center}
\begin{tikzpicture}[>=angle 90]
	\node (1) at (0,1) {$pt$};
	\node (2) at (2,1) {$S$};
	\node (3) at (4,1) {$S \times S.$};

	\tikzset{font=\scriptsize};
	\draw[->] (1) to node [above] {$e$} (2);
	\draw[->] (2) to node [above] {$m^t$} (3);
\end{tikzpicture}
\end{center}
Indeed, the $*$-structure is nothing but the inverse of the corresponding symplectic groupoid. Thus, such a structure is not really extra data on a symplectic monoid, but can rather be thought of as an extra condition the monoid structure itself must satisfy. Still, we will continue to refer directly to the strongly positive $*$-structure on a symplectic groupoid to avoid having to reconstruct it from the monoid data.

\begin{ex}
Recall that the cotangent bundle of a groupoid has a natural symplectic groupoid structure. As a specific case of the previous example, let us explicitly spell out the symplectic monoid structure on $T^*G$ for a Lie group $G$. The product
\[
T^*G \times T^*G \to T^*G
\]
is obtained by applying the cotangent functor to the usual product $G \times G \to G$; explicitly, this is the relation
\[
((g,(dR_h)_g^*\xi),(h,(dL_g)_h^*\xi)) \mapsto (gh,\xi).
\]
The unit $pt \to T^*G$ is given by the lagrangian submanifold $\g^*$ of $T^*G$ and is obtained by applying $T^*$ to the inclusion $pt \to G$ of the identity element. Finally, the $*$-structure is the symplectomorphism $\overline{T^*G} \to T^*G$ given by
\[
(g,-di_g^*\xi) \mapsto (g^{-1},\xi),
\]
where $i: G \to G$ is inversion, and can be obtained as the composition $T^*G \to \overline{T^*G} \to T^*G$ of the Schwartz transform of $T^*G$ followed by the cotangent lift $T^*i$.
\end{ex}

Reversing the arrows in the above definition leads to the notion of a \emph{symplectic comonoid}; we will call the structure morphisms of a symplectic comonoid the \emph{coproduct} and \emph{counit}, and will denote them by $\Delta$ and $\varepsilon$ respectively. Similarly, one can speak of a (strongly positive) $*$-structure on a symplectic comonoid.

\begin{ex}\label{ex:std-com}
Let $M$ be a manifold. Then $T^*M$ has a natural symplectic $*$-comonoid structure, obtained by reversing the arrows in its standard symplectic groupoid structure. To be explicit, the coproduct $T^*M \to T^*M \times T^*M$ is
\[
\Delta: (p,\xi+\eta) \mapsto ((p,\xi),(p,\eta)),
\]
which is obtained as the cotangent lift of the standard diagonal map $M \to M \times M$, and the counit $\varepsilon: T^*M \to pt$ is given by the zero section and is obtained as the cotangent lift of the canonical map $M \to pt$. The $*$-structure is the Schwartz transform.
\end{ex}

Canonical relations which preserve (co)monoid structures are referred to as (co)monoid morphisms:

\begin{defn}
A \emph{monoid morphism} between symplectic monoids $P$ and $Q$ is a canonical relation $L: P \to Q$ such that
\begin{center}
\begin{tikzpicture}[>=angle 90]
	\node (UL) at (0,1) {$P \times P$};
	\node (UR) at (3,1) {$Q \times Q$};
	\node (LL) at (0,-1) {$P$};
	\node (LR) at (3,-1) {$Q$};
	
	\node at (5,0) {$\text{ and }$};
	
	\node (UL2) at (7,1) {$P$};
	\node (UR2) at (10,1) {$Q$};
	\node (LL2) at (8.5,-1) {$pt$};

	\tikzset{font=\scriptsize};
	\draw[->] (UL) to node [above] {$L \times L$} (UR);
	\draw[->] (UL) to node [left] {$m_P$} (LL);
	\draw[->] (UR) to node [right] {$m_Q$} (LR);
	\draw[->] (LL) to node [above] {$L$} (LR);
	
	\draw[->] (UL2) to node [above] {$L$} (UR2);
	\draw[->] (LL2) to node [left] {$e_P$} (UL2);
	\draw[->] (LL2) to node [right] {$e_Q$} (UR2);
\end{tikzpicture}
\end{center}
commute. When the first diagram commutes, we say that $L$ \emph{preserves products}, and when the second commutes, we say that $L$ \emph{preserves units}. Reversing the arrows in the above diagrams leads to the notion of a \emph{comonoid morphism} between symplectic comonoids, and we speak of a canonical relation \emph{preserving coproducts} and \emph{preserving counits}.
\end{defn}

\begin{prop}
The cotangent lift $T^*f: T^*X \to T^*Y$ of a map $f: X \to Y$ is a comonoid morphism with respect to the standard cotangent comonoid structures.
\end{prop}

\begin{proof}
First, the composition $\Delta_Y \circ T^*f$ looks like
\[
(p,df_p^* \xi) \mapsto (f(p),\xi) \mapsto ((f(p),\xi_1),(f(p),\xi_2))
\]
where $\xi = \xi_1 + \xi_2$. The composition $(T^*f \times T^*f) \circ \Delta_X$ looks like
\[
(p,\eta) \mapsto ((p,\eta_1),(p,\eta_2)) \mapsto ((f(p),\gamma_1),(f(p),\gamma_2))
\]
where $\eta_i = df_p^*\gamma_i$ and $\eta = \eta_1 + \eta_2 = df_p^*\gamma_1+df_p^*\gamma_2$. Hence these two compositions agree by the linearity of $df_p^*$.

Similarly, the composition $\varepsilon_Y \circ T^*f$ is
\[
(p,df_p^*\xi) \mapsto (f(p),\xi) \mapsto pt
\]
where $\xi = 0$. Then $df_p^*\xi$ is also $0$ so the composition is $\varepsilon_X$ as required.
\end{proof}

Moreover, as shown in \cite{SZ2}, the above proposition in fact completely characterizes those canonical relations between cotangent bundles which are cotangent lifts.

After considering monoids in $\Symp$, it is natural to want to consider group objects. However, we soon run into the following problem: the diagrams required of a group object $G$ in a category make use of a ``diagonal'' morphism $G \to G \times G$ and a morphism $G \to pt$, but there are no canonical choices for such morphisms in the symplectic category. Indeed, such structures should rather be thought of as coming from a comonoid structure on $G$, and from this point of view the notion of a ``group object'' gets replaced by that of a ``Hopf algebra object'':

\begin{defn}
A \emph{Hopf algebra object} in $\Symp$ consists of a symplectic manifold $S$ together with
\begin{itemize}
\item a symplectic monoid structure $(S,m,e)$,
\item a symplectic comonoid structure $(S,\Delta,\varepsilon)$, and
\item a symplectomorphism $i: S \to S$
\end{itemize}
such that the following diagrams commute, with all compositions strongly transversal:
\begin{itemize}
\item (compatibility between product and coproduct)
\begin{center}
\begin{tikzpicture}[>=angle 90]
	\node (U1) at (0,1) {$S \times S$};
	\node (U2) at (3,1) {$S$};
	\node (U3) at (6,1) {$S \times S$};
	\node (L1) at (0,-1) {$S \times S \times S \times S$};
	\node (L3) at (6,-1) {$S \times S \times S \times S$};

	\tikzset{font=\scriptsize};
	\draw[->] (U1) to node [above] {$m$} (U2);
	\draw[->] (U2) to node [above] {$\Delta$} (U3);
	\draw[->] (U1) to node [left] {$\Delta \times \Delta$} (L1);
	\draw[->] (L3) to node [right] {$m \times m$} (U3);
	\draw[->] (L1) to node [above] {$id \times \sigma \times id$} (L3);
\end{tikzpicture}
\end{center}
where $\sigma: S \times S \to S \times S$ is the symplectomorphism exchanging components,
\item (compatibilities between product and counit, between coproduct and unit, and between unit and counit respectively)
\begin{center}
\begin{tikzpicture}[>=angle 90]
	\node (U1) at (0,1) {$S \times S$};
	\node (U2) at (3,1) {$S$};
	\node (L1) at (1.5,-1) {$pt$};
	
	\node at (3,0) {$\text{,}$};
	
	\node (U3) at (4,1) {$S$};
	\node (U4) at (7,1) {$S \times S$};
	\node (L2) at (5.5,-1) {$pt$};
	
	\node at (8,0) {, and };
	
	\node (U5) at (10.5,1) {$S$};
	\node (L5) at (9,-1) {$pt$};
	\node (L6) at (12,-1) {$pt$};

	\tikzset{font=\scriptsize};
	\draw[->] (U1) to node [above] {$m$} (U2);
	\draw[->] (U1) to node [left] {$\varepsilon \times \varepsilon$} (L1);
	\draw[->] (U2) to node [right] {$\varepsilon$} (L1);
	
	\draw[->] (U3) to node [above] {$\Delta$} (U4);
	\draw[->] (L2) to node [left] {$e$} (U3);
	\draw[->] (L2) to node [right] {$e \times e$} (U4);
	
	\draw[->] (L5) to node [left] {$e$} (U5);
	\draw[->] (U5) to node [right] {$\varepsilon$} (L6);
	\draw[->] (L5) to node [above] {$id$} (L6);
\end{tikzpicture}
\end{center}
\item (antipode conditions)
\begin{center}
\begin{tikzpicture}[>=angle 90]
	\node (U1) at (0,1) {$S \times S$};
	\node (U3) at (4,1) {$S \times S$};
	\node (L1) at (0,-1) {$S$};
	\node (L2) at (2,-1) {$pt$};
	\node (L3) at (4,-1) {$S$};
	
	\node at (5.5,0) {$\text{and}$};
	
	\node (U4) at (7,1) {$S \times S$};
	\node (U6) at (11,1) {$S \times S$};
	\node (L4) at (7,-1) {$S$};
	\node (L5) at (9,-1) {$pt$};
	\node (L6) at (11,-1) {$S$};

	\tikzset{font=\scriptsize};
	\draw[->] (U1) to node [above] {$id \times i$} (U3);
	\draw[->] (L1) to node [left] {$\Delta$} (U1);
	\draw[->] (U3) to node [right] {$m$} (L3);
	\draw[->] (L1) to node [above] {$\varepsilon$} (L2);
	\draw[->] (L2) to node [above] {$e$} (L3);
	
	\draw[->] (U4) to node [above] {$i \times id$} (U6);
	\draw[->] (L4) to node [left] {$\Delta$} (U4);
	\draw[->] (U6) to node [right] {$m$} (L6);
	\draw[->] (L4) to node [above] {$\varepsilon$} (L5);
	\draw[->] (L5) to node [above] {$e$} (L6);
\end{tikzpicture}
\end{center}
\end{itemize}
\end{defn}

\begin{ex}\label{ex:hopf}
Equip the cotangent bundle $T^*G$ of a Lie group $G$ with the comonoid structure of Example \ref{ex:std-com} and the monoid structure coming from its symplectic groupoid structure over $\g^*$. Then these two structures together with the ``antipode'' $T^*i$, where $i: G \to G$ is inversion, make $T^*G$ into a Hopf algebra object in the symplectic category. 

Note that we have the same structure on the cotangent bundle of a more general Lie groupoid, but this will not form a Hopf algebra object; in particular, the antipode conditions fail owing to the fact that groupoid product and is not defined on all of $G \times G$.
\end{ex}

We will return to such structures, and generalizations, in the next chapter.

\section{Actions in the Symplectic Category}
We now turn to Hamiltonian actions of symplectic groupoids. As in any monoidal category, we can define the notion of an action of a monoid object:

\begin{defn}
Let $S$ be a symplectic monoid. An \emph{action} of $S$ on a symplectic manifold $Q$ in the symplectic category is a canonical relation $\tau: S \times Q \to Q$ so that the diagrams
\begin{center}
\begin{tikzpicture}[>=angle 90]
	\node (UL) at (0,1) {$S \times S \times Q$};
	\node (UR) at (3,1) {$S \times Q$};
	\node (LL) at (0,-1) {$S \times Q$};
	\node (LR) at (3,-1) {$Q$};
	
	\node at (5,0) {$\text{and}$};
	
	\node (UL2) at (7,1) {$P$};
	\node (UR2) at (10,1) {$S \times Q$};
	\node (LR2) at (8.5,-1) {$Q,$};

	\tikzset{font=\scriptsize};
	\draw[->] (UL) to node [above] {$id \times \tau$} (UR);
	\draw[->] (UL) to node [left] {$m \times id$} (LL);
	\draw[->] (UR) to node [right] {$\tau$} (LR);
	\draw[->] (LL) to node [above] {$\tau$} (LR);
	
	\draw[->] (UL2) to node [above] {$e \times id$} (UR2);
	\draw[->] (UR2) to node [right] {$\tau$} (LR2);
	\draw[->] (UL2) to node [left] {$id$} (LR2);
\end{tikzpicture}
\end{center}
which say that $\tau$ is compatible with the product and unit of $S$ respectively, commute. We again require that all compositions above be strongly transversal.
\end{defn}

\begin{ex}
The cotangent lift $T^*\tau: T^*G \times T^*M \to T^*M$ of Proposition \ref{prop:action-lift} defines an action in the symplectic category. The diagrams in the above definition commute simply because $T^*$ preserves commutative diagrams. As we will see, this action completely encodes the induced lifted Hamiltonian action of $G$ on $T^*M$.
\end{ex}

As a generalization of this example, suppose that we have a Hamiltonian action of $G$ on $P$ with equivariant momentum map $\mu: P \to \g^*$. Then the relation
\[
((g,\mu(p)),p) \mapsto gp
\]
defines an action $T^*G \times P \to P$ of $T^*G$ on $P$ in the symplectic category. This is essentially the content of Proposition \ref{prop:act}. To be precise, given smooth maps $G \times P \to P$ and $\mu: P \to \g^*$, the above relation defines an action in $\Symp$ if and only if $G \times P \to P$ is a Hamiltonian action with momentum map $\mu$; the diagrams which say that $T^*G \times P \to P$ is an action are equivalent to $G \times P \to P$ being an action and $\mu$ being equivariant, and the condition that $T^*G \times P \to P$ be a canonical relation is equivalent to $\mu$ being a momentum map.

This observation generalizes as follows:

\begin{thrm}
A Hamiltonian action of a symplectic groupoid $S \rightrightarrows P$ is the same as an action in the symplectic category of the corresponding symplectic monoid.
\end{thrm}

\begin{proof}
Suppose that $S \rightrightarrows P$ acts in a Hamiltonian way on a symplectic manifold $Q$. Then the graph of the action map $\tau: S \times_P Q \to Q$ is a lagrangian submanifold of $\overline{S} \times \overline{Q} \times Q$. Viewing this as a canonical relation $S \times Q \to Q$, it is then straightforward to check that this defines  an action of $S$ on $Q$ in the symplectic category.

Conversely, suppose that $\tau: S \times Q \to Q$ is an action of $S$ on $Q$ in the symplectic category, so that the following diagrams commute:
\begin{center}
\begin{tikzpicture}[>=angle 90]
	\node (UL) at (0,1) {$S \times S \times Q$};
	\node (UR) at (3,1) {$S \times Q$};
	\node (LL) at (0,-1) {$S \times Q$};
	\node (LR) at (3,-1) {$Q$};
	
	\node at (5,0) {$\text{and}$};
	
	\node (UL2) at (7,1) {$Q$};
	\node (UR2) at (10,1) {$S \times Q$};
	\node (LR2) at (8.5,-1) {$Q.$};

	\tikzset{font=\scriptsize};
	\draw[->] (UL) to node [above] {$id \times \tau$} (UR);
	\draw[->] (UL) to node [left] {$m \times id$} (LL);
	\draw[->] (UR) to node [right] {$\tau$} (LR);
	\draw[->] (LL) to node [above] {$\tau$} (LR);
	
	\draw[->] (UL2) to node [above] {$e \times id$} (UR2);
	\draw[->] (UR2) to node [right] {$\tau$} (LR2);
	\draw[->] (UL2) to node [left] {$id$} (LR2);
\end{tikzpicture}
\end{center}
We first extract a moment map $J: Q \to P$. Given $q \in Q$, the commutativity of the second diagram implies that there exists $p \in P$ such that
\[
\tau: (p,q) \mapsto q.
\]
Suppose that $p, p' \in P$ both have this property. Then $(p,p',q,q)$ lies in the composition $\tau \circ (id \times \tau)$ of the first diagram, whose commutativity then requires that this also lie in the composition $\tau \circ (m \times id)$. In particular, this requires that $(p,p')$ lie in the domain of $m$; since $m$ is the product of a symplectic groupoid structure on $S$, this forces $p = p'$ since the only way two units of a groupoid are composable is when they are the same. Thus, the element $p \in P$ such that $\tau: (p,q) \mapsto q$ is unique, and we define $J(q)$ to be this element.

We now claim that the action $\tau: S \times Q \to Q$ is single-valued. Indeed, suppose that
\[
\tau: (s,q) \mapsto q' \text{ and } \tau: (s,q) \mapsto q''.
\]
The associativity of $\tau$ then implies that
\[
\tau: (s^{-1},q') \mapsto q \text{ and } (s^{-1},q'') \mapsto q
\]
where $s^{-1}$ is the inverse of $s$ under the symplectic groupoid structure on $S$. On the one hand, the composition $\tau \circ (m \times id)$ in the first diagram above then gives
\[
\tau \circ (m \times id): (s,s^{-1},q') \mapsto (J(q'),q') \mapsto q',
\]
Note in particular that this is the only possible element an element of the form $(s,s^{-1},q')$ can map to under this composition. On the other hand, the composition $\tau \circ (id \times \tau)$ gives
\[
\tau \circ (id \times \tau): (s,s^{-1},q') \mapsto (s,q) \mapsto q''.
\]
Thus since these two compositions agree, we must have $q' = q''$, so $\tau$ is single-valued as claimed.

It is then straightforward to check that $\tau: S \times Q \to Q$ together with the moment map $J: Q \to P$ defines a Hamiltonian action of $S \rightrightarrows P$ on $Q$.
\end{proof}

In particular, the example of the cotangent lift of a group action has the following generalization:

\begin{prop}\label{prop:lifted-cot-action}
Suppose that a groupoid $G$ acts on a manifold $N$, and consider the map defining the action as a relation $\tau: G \times N \to N$. Then the cotangent lift $T^*\tau: T^*G \times T^*N \to T^*N$ defines an action in the symplectic category and hence an action of the symplectic groupoid $T^*G \rightrightarrows A^*$ on $T^*N$.
\end{prop}

\section{Fiber Products}
As we have seen, the symplectic category has many nice properties and applications. However, it has some problems (apart from the problem of having well-defined compositions) as well. In particular, consider the following setup.
Suppose that $f,g: M \to N$ are surjective submersions. Then we have the fiber product diagram
\begin{center}
\begin{tikzpicture}[thick]
	\node (UL) at (0,1) {$M_1 \times_N M_2$};
	\node (UR) at (3,1) {$M_2$};
	\node (LL) at (0,-1) {$M_1$};
	\node (LR) at (3,-1) {$N$};

	\tikzset{font=\scriptsize};
	\draw[->] (UL) to node [above] {$pr_2$} (UR);
	\draw[->] (UL) to node [left] {$pr_1$} (LL);
	\draw[->] (UR) to node [right] {$g$} (LR);
	\draw[->] (LL) to node [above] {$f$} (LR);
\end{tikzpicture}
\end{center}
in $\Man$. Applying the cotangent functor produces the diagram
\begin{center}
\begin{tikzpicture}[thick]
	\node (UL) at (0,1) {$T^*(M_1 \times_N M_2)$};
	\node (UR) at (4,1) {$T^*M_2$};
	\node (LL) at (0,-1) {$T^*M_1$};
	\node (LR) at (4,-1) {$T^*N$};

	\tikzset{font=\scriptsize};
	\draw[->] (UL) to node [above] {$T^*pr_2$} (UR);
	\draw[->] (UL) to node [left] {$T^*pr_1$} (LL);
	\draw[->] (UR) to node [right] {$T^*g$} (LR);
	\draw[->] (LL) to node [above] {$T^*f$} (LR);
\end{tikzpicture}
\end{center}
in the symplectic category. In this section, we are interested in the extent to which this can be viewed as a ``fiber-product'' diagram in $\Symp$. The immediate motivation for considering this is the case where $f$ and $g$ are the source and target of a Lie groupoid, in which case the fiber product $G \times_M G$ consists of the pairs of composable arrows.

To naively construct the fiber product of the canonical relations $T^*f$ and $T^*g$ above, we might try to proceed as in the case of actual maps and form
\[
T^*M_1 \times^{naive}_{T^*N} T^*M_2 := \{((p,\xi),(q,\eta)) \in T^*M_1 \times T^*M_2\ |\ T^*f(p,\xi) = T^*g(q,\eta)\}.
\]
(Recall that if $f,g$ are surjective submersions, $T^*f$ and $T^*g$ are reductions so that they are, in particular, single-valued.) The condition above can equivalently be expressed as
\[
T^*M_1 \times^{naive}_{T^*N} T^*M_2 = (T^*g)^t \circ (T^*f).
\]
Note that, since $T^*f$ is a reduction, this composition is strongly transversal and so the result is a smooth manifold.

The question is now whether or not this manifold is isomorphic to $T^*(M_1 \times_N M_2)$. One case in which this indeed works is when $T^*f$ and $T^*g$ are actual maps:

\begin{prop}
Consider the category whose objects are symplectic manifolds and morphisms are symplectic maps. Then the above diagram is a fiber product diagram in this category.
\end{prop}

\begin{proof}
When considering only symplectic maps, the above ``naive'' fiber product is an actual fiber product, so we need only show that it is isomorphic to $T^*(M_1 \times_N M_2)$.

Consider the inclusion of the naive fiber product into $\overline{T^*M_1} \times T^*M_2$:
\[
(T^*g)^t \circ (T^*f) \hookrightarrow \overline{T^*M_1} \times T^*M_2.
\]
Composing this with the canonical relations
\[
\overline{T^*M_1} \times T^*M_2 \to T^*(M_1 \times M_2) \red T^*(M_1 \times_N M_2),
\]
where the first is induced by the Schwartz transform of $T^*M_1$ and the second is  obtained by reducing the coisotropic $T^*(M_1 \times M_2)|_{M_1 \times_N M_2}$, gives a smooth relation
\[
(T^*g)^t \circ (T^*f) \to T^*(M_1 \times_N M_2).
\]
It is then simple to check that this relation is actually a map, and indeed bijective. Thus $T^*(M_1 \times_N M_2)$ is the fiber product of $T^*f$ and $T^*g$ in the category of symplectic maps.
\end{proof}

However, the situation is not so good in general. First, the composition
\[
(T^*g)^t \circ (T^*f) \to T^*(M_1 \times_N M_2)
\]
of the above proposition, while still a map, is only an inclusion in general. Second, the naive fiber product $(T^*g)^t \circ (T^*f)$ does not have an obvious symplectic structure and so is not necessarily an object in $\Symp$.

Instead, we can abandon the idea of trying to use the naive fiber product as above, and ask whether $T^*(M_1 \times_N M_2)$ is still the correct fiber product in the symplectic category. One quickly realizes again that this is not the case, and that fiber products do not exist in the symplectic category in general except under special circumstances---say the fiber product of symplectomorphisms. Indeed, one can check that even the simplest type of canonical relations---morphisms to a point--do not admit a fiber product. This is true even for the category of set-theoretic relations between sets; in other words, this is a drawback of working with relations in general and is not unique to smooth nor canonical relations.

This lack of fiber products in general, and in particular in the Lie groupoid case mentioned above, was one issue that led to some of the ideas considered in the next chapter. Another idea, to which we return in Chapter~\ref{chap:stacks}, is to dispense of fiber products altogether and instead work with simplicial objects.

\section{Simplicial Symplectic Manifolds}
\begin{defn}
A \emph{simplicial symplectic manifold} is a simplicial object in the symplectic category, i.e. a functor
\[
P: \Delta^{op} \to \textbf{Symp},
\]
where $\Delta$ is the category whose objects are sets $[n] := \{0,1,\ldots,n\}$ and morphisms $[n] \to [m]$ are order-preserving maps.
\end{defn}

In concrete terms, the above definition boils down to the following: for each $n \ge 0$ we have a symplectic manifold $P_{n} := P([n])$, canonical relations
\begin{center}
\begin{tikzpicture}[>=angle 90]
	\def\A{.3};
	\def\B{1.4};
	\def\C{2.5};

	\node at (0,0) {$\cdots$};
	\node at (1.1,0) {$P_2$};
	\node at (2.2,0) {$P_1$};
	\node at (3.3,0) {$P_0$};

	\tikzset{font=\scriptsize};
	\draw[->] (\A,+.18) -- (\A+.5,0+.18);
	\draw[->] (\A,+.06) -- (\A+.5,+.06);
	\draw[->] (\A,-.06) -- (\A+.5,-.06);
	\draw[->] (\A,-.18) -- (\A+.5,0-.18);
	
	\draw[->] (\B,+.12) -- (\B+.5,+.12);
	\draw[->] (\B,0) -- (\B+.5,0);
	\draw[->] (\B,-.12) -- (\B+.5,-.12);
	
	\draw[->] (\C,0+.06) -- (\C+.5,0+.06);
	\draw[->] (\C,0-.06) -- (\C+.5,0-.06);
\end{tikzpicture}
\end{center}
called the \emph{face} morphisms, and canonical relations
\begin{center}
\begin{tikzpicture}[>=angle 90]
	\def\A{.3};
	\def\B{1.4};
	\def\C{2.5};

	\node at (0,0) {$\cdots$};
	\node at (1.1,0) {$P_2$};
	\node at (2.2,0) {$P_1$};
	\node at (3.3,0) {$P_0$};

	\tikzset{font=\scriptsize};
	\draw[<-] (\A,+.12) -- (\A+.5,0+.12);
	\draw[<-] (\A,0) -- (\A+.5,0);
	\draw[<-] (\A,-.12) -- (\A+.5,0-.12);
	
	\draw[<-] (\B,+.06) -- (\B+.5,+.06);
	\draw[<-] (\B,-.06) -- (\B+.5,-.06);
	
	\draw[<-] (\C,0) -- (\C+.5,0);
\end{tikzpicture}
\end{center}
called the \emph{degeneracy} morphisms, so that the following holds (which we recall from \cite{CZ}): denoting the face morphisms coming out of $P_n$ by $d^n_i$ and the degeneracy morphisms coming out of $P_n$ by $s^n_i$, we require that
\[
d^{n-1}_i d^n_j = d^{n-1}_{j-1} d^n_i \text{ for } i < j, s^n_i s^{n-1}_j = s^n_{j+1} s^{n-1}_i \text{ if } i \le j,
\]
\[
d^n_i s^{n-1}_j = s^{n-2}_{j-1} d^{n-1}_i \text{ for } i < j, d^n_j s^{n-1}_j = id = d^n_{j+1} s^{n-1}_j, d^n_i s^{n-1}_j = s^{n-2}_j d^{n-1}_{i-1} \text{ for } i > j+1.
\]
This identities simply that the $d^n_i$ and $s^n_i$ behave as if they were the face and degeneracy maps of simplices. We also require that all compositions involved be strongly transversal, and will denote the above simplicial symplectic manifold simply by $P_\bullet$.

\begin{ex}\label{ex:triv-simp-symp}
For any symplectic manifold $P$, we have the trivial simplicial symplectic manifold
\begin{center}
\begin{tikzpicture}[>=angle 90]
	\def\A{.3};
	\def\B{1.4};
	\def\C{2.5};

	\node at (0,0) {$\cdots$};
	\node at (1.1,0) {$P$};
	\node at (2.2,0) {$P$};
	\node at (3.3,0) {$P$};

	\tikzset{font=\scriptsize};
	\draw[->] (\A,+.18) -- (\A+.5,0+.18);
	\draw[->] (\A,+.06) -- (\A+.5,+.06);
	\draw[->] (\A,-.06) -- (\A+.5,-.06);
	\draw[->] (\A,-.18) -- (\A+.5,0-.18);
	
	\draw[->] (\B,+.12) -- (\B+.5,+.12);
	\draw[->] (\B,0) -- (\B+.5,0);
	\draw[->] (\B,-.12) -- (\B+.5,-.12);
	
	\draw[->] (\C,0+.06) -- (\C+.5,0+.06);
	\draw[->] (\C,0-.06) -- (\C+.5,0-.06);
\end{tikzpicture}
\end{center}
where all relations are $id$.
\end{ex}

\begin{ex}\label{ex:cot-simp-symp}
Let $G \rightrightarrows M$ be a Lie groupoid and form the simplicial nerve\footnote{We refer to \cite{Meh} for the definition of the simplicial nerve of a groupoid, and related simplicial constructions.}
\begin{center}
\begin{tikzpicture}[>=angle 90]
	\def\A{.3};
	\def\B{2.6};
	\def\C{3.7};

	\node at (0,0) {$\cdots$};
	\node at (1.7,0) {$G \times_M G$};
	\node at (3.4,0) {$G$};
	\node at (4.6,0) {$M.$};

	\tikzset{font=\scriptsize};
	\draw[->] (\A,+.18) -- (\A+.5,0+.18);
	\draw[->] (\A,+.06) -- (\A+.5,+.06);
	\draw[->] (\A,-.06) -- (\A+.5,-.06);
	\draw[->] (\A,-.18) -- (\A+.5,0-.18);
	
	\draw[->] (\B,+.12) -- (\B+.5,+.12);
	\draw[->] (\B,0) -- (\B+.5,0);
	\draw[->] (\B,-.12) -- (\B+.5,-.12);
	
	\draw[->] (\C,0+.06) -- (\C+.5,0+.06);
	\draw[->] (\C,0-.06) -- (\C+.5,0-.06);
\end{tikzpicture}
\end{center}
The middle morphism $m: G \times_{M} G \to G$ is multiplication, the morphisms $G \rightrightarrows M$ are the source and target, and $e: M \to G$ is the identity embedding.
Applying $T^{*}$ gives the simplicial symplectic manifold
\begin{center}
\begin{tikzpicture}[>=angle 90]
	\def\A{.3};
	\def\B{3.3};
	\def\C{4.8};

	\node at (0,0) {$\cdots$};
	\node at (2.1,0) {$T^*(G \times_M G)$};
	\node at (4.3,0) {$T^*G$};
	\node at (5.9,0) {$T^*M.$};

	\tikzset{font=\scriptsize};
	\draw[->] (\A,+.18) -- (\A+.5,0+.18);
	\draw[->] (\A,+.06) -- (\A+.5,+.06);
	\draw[->] (\A,-.06) -- (\A+.5,-.06);
	\draw[->] (\A,-.18) -- (\A+.5,0-.18);
	
	\draw[->] (\B,+.12) -- (\B+.5,+.12);
	\draw[->] (\B,0) -- (\B+.5,0);
	\draw[->] (\B,-.12) -- (\B+.5,-.12);
	
	\draw[->] (\C,0+.06) -- (\C+.5,0+.06);
	\draw[->] (\C,0-.06) -- (\C+.5,0-.06);
\end{tikzpicture}
\end{center}

Alternatively, we can view the relation $T^*m$ as a canonical relation $T^*G \times T^*G \to T^*G$, in which case we can form a simplicial symplectic manifold of the form
\begin{center}
\begin{tikzpicture}[>=angle 90]
	\def\A{.3};
	\def\B{3.3};
	\def\C{4.8};

	\node at (0,0) {$\cdots$};
	\node at (2.1,0) {$T^*G \times T^*G$};
	\node at (4.3,0) {$T^*G$};
	\node at (5.9,0) {$T^*M,$};

	\tikzset{font=\scriptsize};
	\draw[->] (\A,+.18) -- (\A+.5,0+.18);
	\draw[->] (\A,+.06) -- (\A+.5,+.06);
	\draw[->] (\A,-.06) -- (\A+.5,-.06);
	\draw[->] (\A,-.18) -- (\A+.5,0-.18);
	
	\draw[->] (\B,+.12) -- (\B+.5,+.12);
	\draw[->] (\B,0) -- (\B+.5,0);
	\draw[->] (\B,-.12) -- (\B+.5,-.12);
	
	\draw[->] (\C,0+.06) -- (\C+.5,0+.06);
	\draw[->] (\C,0-.06) -- (\C+.5,0-.06);
\end{tikzpicture}
\end{center}
where the degree $n$ piece for $n > 1$ is $(T^*G)^n$. Now the first and third degree $2$ face morphisms $T^*G \times T^*G \to T^*G$ are the ``projections'' obtained by taking the cotangent lifts of the smooth relations
\[
(g,h) \mapsto g \text{ and } (g,h) \mapsto h, \text{ for $g,h \in G$ such that } r(g) = \ell(h).
\]
A simple calculation shows that these are respectively equal to the compositions
\begin{center}
\begin{tikzpicture}[>=angle 90]
	\node (1) at (0,1) {$T^*G \times T^*G$};
	\node (2) at (4,1) {$T^*G \times T^*M$};
	\node (3) at (8.5,1) {$T^*G \times T^*G$};
	\node (4) at (12,1) {$T^*G$};

	\tikzset{font=\scriptsize};
	\draw[->] (1) to node [above] {$id \times T^*\ell$} (2);
	\draw[->] (2) to node [above] {$id \times (T^*r)^t$} (3);
	\draw[->] (3) to node [above] {$(T^*\Delta)^t$} (4);
\end{tikzpicture}
\end{center}
and
\begin{center}
\begin{tikzpicture}[>=angle 90]
	\node (1) at (0,1) {$T^*G \times T^*G$};
	\node (2) at (4,1) {$T^*G \times T^*M$};
	\node (3) at (8.5,1) {$T^*G \times T^*G$};
	\node (4) at (12,1) {$T^*G.$};

	\tikzset{font=\scriptsize};
	\draw[->] (1) to node [above] {$T^*r \times id$} (2);
	\draw[->] (2) to node [above] {$(T^*\ell)^t \times id$} (3);
	\draw[->] (3) to node [above] {$(T^*\Delta)^t$} (4);
\end{tikzpicture}
\end{center}
The rest of the face and degeneracy maps can be similarly described.
\end{ex}

\begin{ex}\label{ex:bad-simp-symp}
Given a Lie groupoid $G$, we may also attempt to construct a simplicial symplectic manifold of the form
\begin{center}
\begin{tikzpicture}[>=angle 90]
	\def\A{.3};
	\def\B{3.3};
	\def\C{4.8};

	\node at (0,0) {$\cdots$};
	\node at (2.1,0) {$T^*G \times T^*G$};
	\node at (4.3,0) {$T^*G$};
	\node at (5.6,0) {$pt$};

	\tikzset{font=\scriptsize};
	\draw[->] (\A,+.18) -- (\A+.5,0+.18);
	\draw[->] (\A,+.06) -- (\A+.5,+.06);
	\draw[->] (\A,-.06) -- (\A+.5,-.06);
	\draw[->] (\A,-.18) -- (\A+.5,0-.18);
	
	\draw[->] (\B,+.12) -- (\B+.5,+.12);
	\draw[->] (\B,0) -- (\B+.5,0);
	\draw[->] (\B,-.12) -- (\B+.5,-.12);
	
	\draw[->] (\C,0+.06) -- (\C+.5,0+.06);
	\draw[->] (\C,0-.06) -- (\C+.5,0-.06);
\end{tikzpicture}
\end{center}
Indeed, we simply consider the previous example in the case where $G \rightrightarrows pt$ is a Lie group and note that all the resulting structure morphisms make sense for more general groupoids.

However, some of these compositions turn out to not be strongly transversal (such as the composition of the degeneracy $pt \to T^*G$ with either face map $T^*G \to pt$), and so this is not an allowed simplicial object in our category. We will see in Chapter $5$ the problems that non-strongly transversal compositions may cause.
\end{ex}

\begin{ex}\label{ex:action-grpd}
For a Hamiltonian action of $G$ on $P$, form the corresponding action $\tau: T^{*}G \times P \to P$ in $\Symp$. Then we have a simplicial symplectic manifold
\begin{center}
\begin{tikzpicture}[>=angle 90]
	\def\A{.3};
	\def\B{4};
	\def\C{6.3};

	\node at (0,0) {$\cdots$};
	\node at (2.4,0) {$T^*G \times T^*G \times P$};
	\node at (5.4,0) {$T^*G \times P$};
	\node at (7,0) {$P$};

	\tikzset{font=\scriptsize};
	\draw[->] (\A,+.18) -- (\A+.5,0+.18);
	\draw[->] (\A,+.06) -- (\A+.5,+.06);
	\draw[->] (\A,-.06) -- (\A+.5,-.06);
	\draw[->] (\A,-.18) -- (\A+.5,0-.18);
	
	\draw[->] (\B,+.12) -- (\B+.5,+.12);
	\draw[->] (\B,0) -- (\B+.5,0);
	\draw[->] (\B,-.12) -- (\B+.5,-.12);
	
	\draw[->] (\C,0+.06) -- (\C+.5,0+.06);
	\draw[->] (\C,0-.06) -- (\C+.5,0-.06);
\end{tikzpicture}
\end{center}
where the face maps in degree $2$ are
\[
id \times \tau,\ T^{*}m \times id,\ G \times id \times id,
\]
the face maps in degree $1$ are $\tau$ and $G \times id$, the degeneracy $P \to T^{*}G \times P$ is $\mathfrak g^{*} \times id$, and the other morphisms can be guessed from these.
\end{ex}

\begin{defn}
A \emph{$*$-structure} on a simplicial symplectic manifold $P_\bullet$ is a simplicial isomorphism
\[
I: P_\bullet \to (P_\bullet)^{op}
\]
such that $I^2 = id$, where $(P_\bullet)^{op}$ denotes the simplicial symplectic manifold obtained by reversing the order of the face and degeneracy maps in $P_\bullet$.
\end{defn}

The previous examples all admit such structures: for a trivial simplicial symplectic manifold it is the identity morphism, for the simplicial structures obtained from the cotangent bundle of a groupoid it is the simplicial morphism induced by the cotangent lifts of the maps
\[
(g_1,\ldots,g_n) \mapsto (g_n^{-1},\ldots,g_1^{-1}),
\]
and in the last example it is induced by a combination of the previous one and the morphism
\[
((g,-Ad_{g^{-1}}^*\theta + \mu(p)), p) \mapsto ((g^{-1},\xi), gp)
\]
in degree $2$ where $Ad^*$ denotes the coadjoint action of $G$ on $\g^*$.

\section{Towards Groupoids in $\Symp$}
Let us recall our original motivation: to understand the object $T^*G \rightrightarrows T^*M$ in $\Symp$ resulting from applying $T^*$ to a groupoid $G \rightrightarrows M$ in $\Man$. The notion of a \emph{groupoid object} in a category requires the existence of certain fiber products; indeed, the groupoid product is defined only as a map on a certain fiber product $G \times_M G$. The analogous construction in the symplectic category would require us to view $T^*(G \times_M G)$ as a fiber product of certain canonical relations.

As we have seen in the previous sections, making this precise is not in general possible. One way around this will be to work with the simplicial objects of the previous section---a point of view we consider briefly in the next chapter but more so in Chapter $5$.

However, there is another approach we can take: working with relations instead of maps, we can consider groupoid products $G \times_M G \to G$ as relations $G \times G \to G$, and hence we can consider the corresponding lift $T^*G \times T^*G \to T^*G$ as the ``groupoid product'' of $T^*G \rightrightarrows T^*M$. All commutative diagrams required in the definition of a groupoid object then hold, with one final caveat: these diagrams require the use of a ``diagonal'' morphism $T^*G \to T^*G \times T^*G$, which is extra data in the symplectic category (due to the lack of fiber products) as opposed to a canonical construction in $\Man$.

Taking this extra structure into account will lead us to the notion of a \emph{symplectic hopfoid} in the next chapter. As we will see, the compatibility between the ``product'' $T^*G \times T^*G \to T^*G$ and ``diagonal'' $T^*G \to T^*G \times T^*G$  required to have a ``groupoid-like'' structure is no accident: it reflects the fact that $T^*G$ is not simply a symplectic groupoid, but rather a \emph{symplectic double groupoid}.
\chapter{Double Groupoids}\label{chap:dbl-grpds}

Symplectic double groupoids arose \cite{LW} in attempts to integrate Poisson-Lie groups. Indeed, given a Poisson-Lie group $P$ with $P$ an integrable Poisson manifold, the integrating symplectic groupoid $S \rightrightarrows P$ often has an additional symplectic groupoid structure over $P^*$, the \emph{dual} Poisson-Lie group of $P$. This observation generalizes in many situations to general Poisson groupoids $P \rightrightarrows M$.

The aim of this chapter is to derive a characterization of symplectic double groupoids in terms of the symplectic category. The resulting object in the symplectic category, which we call a \emph{symplectic hopfoid}, should be thought of as a kind of ``groupoid'' or ``Hopf-algebroid'' object. The main example of this correspondence comes from the cotangent bundle of a Lie groupoid.

\section{Preliminaries on Double Groupoids}
Intuitively, a double Lie groupoid is a groupoid object in the category of Lie groupoids. To be precise:

\begin{defn}
A \emph{double Lie groupoid} is a diagram
\begin{center}
\begin{tikzpicture}[>=angle 90,scale=.85]
	\node at (0,0) {$H$};
	\node at (2,0) {$M$};
	\node at (0,2) {$D$};
	\node at (2,2) {$V$};

	\tikzset{font=\scriptsize};
	\draw[->] (.4,.08) -- (1.6,.08);
	\draw[->] (.4,-.08) -- (1.6,-.08);
	
	\draw[->] (-.08,1.6) -- (-.08,.4);
	\draw[->] (.08,1.6) -- (.08,.4);
	
	\draw[->] (.4,2+.08) -- (1.6,2+.08);
	\draw[->] (.4,2-.08) -- (1.6,2-.08);
	
	\draw[->] (2-.08,1.6) -- (2-.08,.4);
	\draw[->] (2+.08,1.6) -- (2+.08,.4);
\end{tikzpicture}
\end{center}
of Lie groupoids such that the structure maps of the top and bottom groupoids give homomorphisms from the left groupoid to the right groupoid, and vice-versa. We will refer to the  four groupoid structures involved as the top, bottom, left, and right groupoids. For technical reasons, we will assume that the \emph{double source map}
\[
D \to H \times V
\]
is a surjective submersion. Also, we often refer to $D$ itself as the double Lie groupoid and to $M$ as its double base. When $M$ is a point, we will call $D$ a \emph{double Lie group}.
\end{defn}

\begin{rmk}
Differentiating the above groupoid structures, the infinitesimal object corresponding to a double Lie groupoid is known as a \emph{double Lie algebroid}. We refer to \cite{M2} and \cite{M3} for further details. The ``intermediate'' object between double Lie groupoids and double Lie algebroids---known as an $\mathcal L\mathcal A$-groupoid---is also studied in \cite{S}.
\end{rmk}

\begin{rmk}
The double groupoid structure on $D$ should not depend on the manner in which we have chosen to draw the diagram above. In other words, we will think of
\begin{center}
\begin{tikzpicture}[>=angle 90,scale=.85]
	\node at (0,0) {$V$};
	\node at (2,0) {$M$};
	\node at (0,2) {$D$};
	\node at (2,2) {$H$};

	\tikzset{font=\scriptsize};
	\draw[->] (.4,.08) -- (1.6,.08);
	\draw[->] (.4,-.08) -- (1.6,-.08);
	
	\draw[->] (-.08,1.6) -- (-.08,.4);
	\draw[->] (.08,1.6) -- (.08,.4);
	
	\draw[->] (.4,2+.08) -- (1.6,2+.08);
	\draw[->] (.4,2-.08) -- (1.6,2-.08);
	
	\draw[->] (2-.08,1.6) -- (2-.08,.4);
	\draw[->] (2+.08,1.6) -- (2+.08,.4);
\end{tikzpicture}
\end{center}
as representing the same double groupoid as above, and will call this the \emph{transpose} of the previous structure.
\end{rmk}

\begin{ex}\label{ex:dmain}
For any Lie groupoid $G \rightrightarrows M$, there is a double Lie groupoid structure on
\begin{center}
\begin{tikzpicture}[>=angle 90]
	\node (LL) at (0,0) {$M$};
	\node (LR) at (2,0) {$M.$};
	\node (UL) at (0,2) {$G$};
	\node (UR) at (2,2) {$G$};

	\tikzset{font=\scriptsize};
	\draw[->] (.4,.07) -- (1.6,.07);
	\draw[->] (.4,-.07) -- (1.6,-.07);
	
	\draw[->] (-.07,1.6) -- (-.07,.4);
	\draw[->] (.07,1.6) -- (.07,.4);
	
	\draw[->] (.4,2.07) -- (1.7,2.07);
	\draw[->] (.4,2-.07) -- (1.7,2-.07);
	
	\draw[->] (2-.07,1.6) -- (2-.07,.4);
	\draw[->] (2+.07,1.6) -- (2+.07,.4);
\end{tikzpicture}
\end{center}
Here, the left and right sides are the given groupoid structures, while the top and bottom are trivial groupoids.
\end{ex}

\begin{ex}\label{ex:dinertia}
Again for any Lie groupoid $G \rightrightarrows M$, there is a double Lie groupoid structure on
\begin{center}
\begin{tikzpicture}[>=angle 90]
	\node (LL) at (0,0) {$M \times M$};
	\node (LR) at (3,0) {$M.$};
	\node (UL) at (0,2) {$G \times G$};
	\node (UR) at (3,2) {$G$};

	\tikzset{font=\scriptsize};
	\draw[->] (.9,.07) -- (2.6,.07);
	\draw[->] (.9,-.07) -- (2.6,-.07);
	
	\draw[->] (-.07,1.6) -- (-.07,.4);
	\draw[->] (.07,1.6) -- (.07,.4);
	
	\draw[->] (.9,2.07) -- (2.5,2.07);
	\draw[->] (.9,2-.07) -- (2.5,2-.07);
	
	\draw[->] (3-.07,1.6) -- (3-.07,.4);
	\draw[->] (3.07,1.6) -- (3.07,.4);
\end{tikzpicture}
\end{center}
Here, the right side is the given groupoid structure, the top and bottom are pair groupoids, and the left is a product groupoid.
\end{ex}

For much of what follows, we will need a consistent labeling of the structure maps involved in the various groupoids above. First, the source, target, unit, inverse, and product of the right groupoid $V \rightrightarrows M$ are respectively
\[
r_V(\cdot), \ell_V(\cdot), 1^V_\cdot, i_V(\cdot), m_V(\cdot,\cdot) \text{ or } \cdot \circ_V \cdot
\]
The structure maps of $H \rightrightarrows M$ will use the same symbols with $V$ replaced by $H$. Now, the structure maps of the top and left groupoids will use the same symbol as those of the \emph{opposite} structure with a tilde on top; so for example, the structure maps of $D \rightrightarrows H$ are
\[
\widetilde{r}_V(\cdot), \widetilde{\ell}_V(\cdot), \widetilde{1}^V_\cdot, \widetilde{i}_V(\cdot), \wt m_V(\cdot,\cdot) \text{ or } \cdot \widetilde{\circ}_V \cdot
\]

To emphasize: the structure maps of the groupoid structure on $D$ \emph{over} $V$ using an $H$, and those of the groupoid structure on $D$ \emph{over} $H$ use a $V$. This has a nice practical benefit in that it is simpler to keep track of the various relations these maps satisfy; for example, the maps $\widetilde{r}_{H}$ and $r_{H}$ give the groupoid homomorphism
\begin{center}
\begin{tikzpicture}[>=angle 90,scale=.75]
	\node (LL) at (0,0) {$H$};
	\node (LR) at (2,0) {$M,$};
	\node (UL) at (0,2) {$D$};
	\node (UR) at (2,2) {$V$};

	\tikzset{font=\scriptsize};
	\draw[->] (UL) to node [above] {$\widetilde{r}_{H}$} (UR);
	
	\draw[->] (-.07,1.6) -- (-.07,.4);
	\draw[->] (.07,1.6) -- (.07,.4);
	
	\draw[->] (2-.07,1.6) -- (2-.07,.4);
	\draw[->] (2+.07,1.6) -- (2+.07,.4);
	
	\draw[->] (LL) -- node [above] {$r_{H}$} (LR);
\end{tikzpicture}
\end{center}
so for instance we have: $\ell_V(\wt r_H(s)) = r_H(\wt \ell_V(s))$, $i_V(\wt r_H(s)) = \wt r_H(\wt i_V(s))$, $\wt r_H(\wt 1^V_v) = 1^V_{r_H(v)}$, etc. We will make extensive use of such identities, and in this sense we can think of the tilde'd maps as \emph{covers} of the untilde'd maps.

We will denote elements of $D$ as squares with sides labeled by the possible sources and targets:
\begin{center}
\begin{tikzpicture}[scale=.75]
	\draw (0,0) rectangle (2,2);
	\node at (1,1) {$s$};

	\tikzset{font=\scriptsize};
	\node at (-.6,1) {$\wt\ell_V(s)$};
	\node at (1,2.4) {$\wt r_{H}(s)$};
	\node at (2.6,1) {$\wt r_V(s)$};
	\node at (1,-.4) {$\wt\ell_{H}(s)$};
\end{tikzpicture}
\end{center}
so that the left and right sides are the target and source of the left groupoid structure while the top and bottom sides are the source and target of the top groupoid structure. This lends itself well to compositions: two squares $s$ and $s'$ are composable in the left groupoid if the right side of the first is the left of the second, i.e. if $\widetilde{r}_{V}(s) = \widetilde{\ell}_{V}(s')$, and the composition $s \,\wt\circ_{V} s'$ in the left groupoid is given by ``horizontal concatenation'':
\begin{center}
\begin{tikzpicture}[scale=.75]
	\node at (1,1) {$s$};
	\node at (3,1) {$s'$};
	\node at (6,1) {$=$};
	\node at (9,1) {$s\, \wt\circ_V s'$};

	\tikzset{font=\scriptsize};
	\draw (0,0) rectangle (2,2);
	\node at (-.6,1) {$\wt\ell_V(s)$};
	\node at (1,2.4) {$\wt r_{H}(s)$};
	\node at (1,-.4) {$\wt\ell_{H}(s)$};
	
	\draw (2,0) rectangle +(2,2);
	\node at (3,2.4) {$\wt r_{H}(s')$};
	\node at (4.7,1) {$\wt r_V(s')$};
	\node at (3,-.4) {$\wt\ell_{H}(s')$};
	
	\draw (8,0) rectangle +(2,2);
	\node at (7.4,1) {$\wt\ell_V(s)$};
	\node at (9,2.4) {$\wt r_{H}(s) \circ_V \wt r_{H}(s')$};
	\node at (10.7,1) {$\wt r_V(s')$};
	\node at (9,-.4) {$\wt\ell_{H}(s) \circ_V \wt\ell_{H}(s')$};
\end{tikzpicture}
\end{center}
Similarly, the groupoid composition in the top groupoid is given by ``vertical concatenation''.

With these notations, the compatibility between the two groupoid products on $D$ can then be expressed as saying that composing vertically and then horizontally in
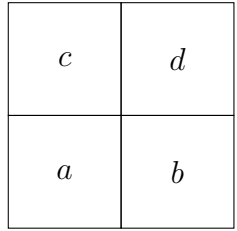
\begin{center}
\begin{tikzpicture}
	\draw (0,0) rectangle (1.5,1.5);
	\node at (.75,.75) {$a$};
	\draw (1.5,0) rectangle (3,1.5);
	\node at (2.25,.75) {$b$};
	\draw (0,1.5) rectangle (1.5,3);
	\node at (.75,2.25) {$c$};
	\draw (1.5,1.5) rectangle (3,3);
	\node at (2.25,2.25) {$d$};
\end{tikzpicture}
\end{center}
produces the same result as composing horizontally and then vertically, whenever all compositions involved are defined.

The core of a double Lie groupoid, introduced in \cite{BM}, is a certain submanifold of $D$ which encodes part of the double groupoid structure:

\begin{defn}
The \emph{core} of a double Lie groupoid $D$ is the submanifold $C$ (the condition of the double source map ensures that this is a submanifold) of elements of $D$ whose sources are both units; that is, the set of elements of the form
\begin{center}
\begin{tikzpicture}[scale=.75]
	\node at (1,1) {$s$};
	
	\tikzset{font=\scriptsize};
	\draw (0,0) rectangle (2,2);
	\node at (-.7,1) {$\wt\ell_{V}(s)$};
	\node at (1,2.4) {$1^V_m$};
	\node at (2.6,1) {$1^{H}_m$};
	\node at (1,-.4) {$\wt\ell_H(s)$};
\end{tikzpicture}
\end{center}
\end{defn}

\begin{thrm}[Brown-Mackenzie, \cite{BM}]
The core of a double Lie groupoid $D$ has a natural Lie groupoid structure over the double base $M$.
\end{thrm}

In the same paper, Brown and Mackenzie also show that certain double Lie groupoids can in fact be recovered from a certain \emph{core diagram}. We will not need this notion here, and refer to the above mentioned paper for further details.

The groupoid structure on the core comes from a combination of the two groupoid structures on $D$. Explicitly, the groupoid product of two elements $s$ and $s'$ in the core can be expressed as composing vertically and then horizontally (or equivalently  horizontally and then vertically) in the following diagram:
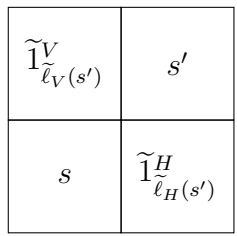
\begin{center}
\begin{tikzpicture}
	\draw (0,0) rectangle (1.5,1.5);
	\node at (.75,.75) {$s$};
	\draw (1.5,0) rectangle (3,1.5);
	\node at (2.25,.75) {$\wt 1^{H}_{\wt\ell_{H}(s')}$};
	\draw (0,1.5) rectangle (1.5,3);
	\node at (.75,2.25) {$\wt 1^V_{\wt\ell_V(s')}$};
	\draw (1.5,1.5) rectangle (3,3);
	\node at (2.25,2.25) {$s'$};
\end{tikzpicture}
\end{center}

The core groupoid of Example \ref{ex:dmain} is the trivial groupoid $M \rightrightarrows M$ while that of Example \ref{ex:dinertia} is $G \rightrightarrows M$ itself.

\section{Symplectic Double Groupoids}
Now we put symplectic structures on double Lie groupoids:

\begin{defn}
A \emph{symplectic double groupoid} is a double Lie groupoid
\begin{center}
\begin{tikzpicture}[>=angle 90,scale=.85]
	\node at (0,0) {$P^*$};
	\node at (2,0) {$M$};
	\node at (0,2) {$S$};
	\node at (2,2) {$P$};

	\tikzset{font=\scriptsize};
	\draw[->] (.4,.08) -- (1.6,.08);
	\draw[->] (.4,-.08) -- (1.6,-.08);
	
	\draw[->] (-.08,1.6) -- (-.08,.4);
	\draw[->] (.08,1.6) -- (.08,.4);
	
	\draw[->] (.4,2+.08) -- (1.6,2+.08);
	\draw[->] (.4,2-.08) -- (1.6,2-.08);
	
	\draw[->] (2-.08,1.6) -- (2-.08,.4);
	\draw[->] (2+.08,1.6) -- (2+.08,.4);
\end{tikzpicture}
\end{center}
where $S$ is a symplectic manifold and the top and left groupoid structures are symplectic groupoids. As for general double groupoids, when $M$ is a point we will call $S$ a \emph{symplectic double group}.
\end{defn}

The notations for the total spaces of the right and bottom groupoids come from the following fact:

\begin{thrm}[Weinstein \cite{W2}, Mackenzie \cite{M}]
If $S$ is a symplectic double groupoid, then $P$ and $P^*$ are Poisson groupoids (for the induced Poisson structures) in duality, meaning that their Lie algebroids are isomorphic to the duals of one another.
\end{thrm}

Indeed, symplectic double groupoids first arose as a tool for integrating Poisson groupoids. The problem of determining when a Poisson groupoid $P \rightrightarrows M$, with $P$ integrable as a Poisson manifold, arises from a symplectic double groupoid is a deep one, with still no satisfactory answer in general; see for example \cite{S}.

The symplectic structure on $S$ endows the core with additional structure as follows:

\begin{thrm}[Mackenzie \cite{M}]
The core $C$ of a symplectic double groupoid $S$ is a symplectic submanifold of $S$ and the induced groupoid structure on $C \rightrightarrows M$ is a symplectic groupoid.
\end{thrm}

Later we present a new proof of this result via a reduction procedure.

\begin{ex}\label{ex:main}
For any groupoid $G \rightrightarrows M$, there is a symplectic double groupoid structure on $T^*G$ of the form
\begin{center}
\begin{tikzpicture}[>=angle 90]
	\node (LL) at (0,0) {$A^*$};
	\node (LR) at (2,0) {$M$};
	\node (UL) at (0,2) {$T^*G$};
	\node (UR) at (2,2) {$G$};

	\tikzset{font=\scriptsize};
	\draw[->] (.4,.07) -- (1.6,.07);
	\draw[->] (.4,-.07) -- (1.6,-.07);
	
	\draw[->] (-.07,1.6) -- (-.07,.4);
	\draw[->] (.07,1.6) -- (.07,.4);
	
	\draw[->] (.6,2.07) -- (1.7,2.07);
	\draw[->] (.6,2-.07) -- (1.7,2-.07);
	
	\draw[->] (2-.07,1.6) -- (2-.07,.4);
	\draw[->] (2+.07,1.6) -- (2+.07,.4);
\end{tikzpicture}
\end{center}
Here, $A$ is the Lie algebroid of $G \rightrightarrows M$, the right groupoid structure is the given one on $G$, the top and bottom are the natural groupoid structures on vector bundles given by fiber-wise addition, and the left structure is the induced symplectic groupoid structure on the cotangent bundle of a Lie groupoid. The core of this symplectic double groupoid is symplectomorphic to $T^*M$, and the core groupoid is simply $T^*M \rightrightarrows M$. Note, in particular, that when $G$ is a Lie group $T^*G$ is a symplectic double group.
\end{ex}

\begin{ex}\label{ex:inertia}
Again for any groupoid $G \rightrightarrows M$, there is a symplectic double groupoid structure on $\overline{T^*G} \times T^*G$ of the form
\begin{center}
\begin{tikzpicture}[>=angle 90]
	\node (LL) at (0,0) {$\overline{A^*} \times A^*$};
	\node (LR) at (3,0) {$A^*$};
	\node (UL) at (0,2) {$\overline{T^*G} \times T^*G$};
	\node (UR) at (3,2) {$T^*G$};

	\tikzset{font=\scriptsize};
	\draw[->] (.9,.07) -- (2.6,.07);
	\draw[->] (.9,-.07) -- (2.6,-.07);
	
	\draw[->] (-.07,1.6) -- (-.07,.4);
	\draw[->] (.07,1.6) -- (.07,.4);
	
	\draw[->] (1.2,2.07) -- (2.5,2.07);
	\draw[->] (1.2,2-.07) -- (2.5,2-.07);
	
	\draw[->] (3-.07,1.6) -- (3-.07,.4);
	\draw[->] (3.07,1.6) -- (3.07,.4);
\end{tikzpicture}
\end{center}
Here, the right side is the induced symplectic groupoid structure on $T^*G$, the top and bottom are pair groupoids, and the left is a product groupoid. The core is symplectomorphic to $T^*G$, and the core groupoid is $T^*G \rightrightarrows A^*$.
\end{ex}

Both of the above examples are special cases of the following result due to Mackenzie:

\begin{thrm}[Mackenzie \cite{M}]\label{thrm:ctdbl}
Let $D$ be a double Lie groupoid. Then the cotangent bundle $T^*D$ has a natural symplectic double groupoid structure
\begin{center}
\begin{tikzpicture}[>=angle 90]
	\node (LL) at (0,0) {$A^*H$};
	\node (LR) at (3,0) {$A^*C$};
	\node (UL) at (0,2) {$T^*D$};
	\node (UR) at (3,2) {$A^*V$};

	\tikzset{font=\scriptsize};
	\draw[->] (.5,.07) -- (2.5,.07);
	\draw[->] (.5,-.07) -- (2.5,-.07);
	
	\draw[->] (-.07,1.6) -- (-.07,.4);
	\draw[->] (.07,1.6) -- (.07,.4);
	
	\draw[->] (.6,2.07) -- (2.5,2.07);
	\draw[->] (.6,2-.07) -- (2.5,2-.07);
	
	\draw[->] (3-.07,1.6) -- (3-.07,.4);
	\draw[->] (3+.07,1.6) -- (3+.07,.4);
\end{tikzpicture}
\end{center}
where $A^*H$ and $A^*V$ are the duals of the Lie algebroids of $H \rightrightarrows M$ and $V \rightrightarrows M$ respectively, and $A^*C$ is the dual of the Lie algebroid of the core groupoid $C \rightrightarrows M$. The core of this symplectic double groupoid is symplectomorphic to $T^*C$ where $C$ is the core of $D$.
\end{thrm}

An explicit description of the groupoid structures on the right and bottom groupoids can also be found in \cite{GM}. In this thesis, in addition to the descriptions of the top and left cotangent groupoids, we will only need to use the units of the right and bottom groupoids: the unit of $A^*V \rightrightarrows A^*C$ is  induced by that of $T^*D \rightrightarrows A^*H$, and so comes from the identification of $A^*H$ with the conormal bundle $N^*H \subset T^*D$, and similarly the unit of $A^*H \rightrightarrows A^*C$ comes from the identification of $A^*V$ with $N^*V$. Example \ref{ex:main} arises from applying this theorem to the double groupoid
\begin{center}
\begin{tikzpicture}[>=angle 90]
	\node (LL) at (0,0) {$M$};
	\node (LR) at (2,0) {$M,$};
	\node (UL) at (0,2) {$G$};
	\node (UR) at (2,2) {$G$};

	\tikzset{font=\scriptsize};
	\draw[->] (.4,.07) -- (1.6,.07);
	\draw[->] (.4,-.07) -- (1.6,-.07);
	
	\draw[->] (-.07,1.6) -- (-.07,.4);
	\draw[->] (.07,1.6) -- (.07,.4);
	
	\draw[->] (.4,2.07) -- (1.7,2.07);
	\draw[->] (.4,2-.07) -- (1.7,2-.07);
	
	\draw[->] (2-.07,1.6) -- (2-.07,.4);
	\draw[->] (2+.07,1.6) -- (2+.07,.4);
\end{tikzpicture}
\end{center}
and Example \ref{ex:inertia} arises from the double groupoid
\begin{center}
\begin{tikzpicture}[>=angle 90]
	\node (LL) at (0,0) {$M \times M$};
	\node (LR) at (3,0) {$M.$};
	\node (UL) at (0,2) {$G \times G$};
	\node (UR) at (3,2) {$G$};

	\tikzset{font=\scriptsize};
	\draw[->] (.9,.07) -- (2.6,.07);
	\draw[->] (.9,-.07) -- (2.6,-.07);
	
	\draw[->] (-.07,1.6) -- (-.07,.4);
	\draw[->] (.07,1.6) -- (.07,.4);
	
	\draw[->] (.9,2.07) -- (2.5,2.07);
	\draw[->] (.9,2-.07) -- (2.5,2-.07);
	
	\draw[->] (3-.07,1.6) -- (3-.07,.4);
	\draw[->] (3.07,1.6) -- (3.07,.4);
\end{tikzpicture}
\end{center}

\section{Realizing the Core via Reduction}
Here we describe a procedure for producing the core of a symplectic double groupoid, and indeed the symplectic groupoid structure on the core, via symplectic reduction. This will eventually provide an interpretation of symplectic double groupoids in terms of the symplectic category, and it will lead to us thinking of the core as the ``base'' of a certain ``Hopf algebroid''-like structure.

Let $S$ be a symplectic double groupoid. The unit submanifold $1^P M$ of the groupoid structure on $P$ is coisotropic in $P$ for the induced Poisson structure. Hence the preimage under the source of the top groupoid (which is a Poisson map)
\[
X := \widetilde{r}_{P^*}^{-1}(1^{P} M) \subseteq S
\]
is a coisotropic submanifold of $S$. In square notation, $X$ consists of those elements of the form:
\begin{center}
\begin{tikzpicture}[scale=.75]
	\node at (1,1) {$s$};

	\tikzset{font=\scriptsize};
	\draw (0,0) rectangle (2,2);
	\node at (-.6,1) {$\wt\ell_P(s)$};
	\node at (1,2.4) {$1^P_m$};
	\node at (2.6,1) {$\wt r_P(s)$};
	\node at (1,-.4) {$\wt\ell_{P^*}(s)$};
\end{tikzpicture}
\end{center}
Note, in particular, that the core $C$ of $S$ sits inside of $X$. Now, there is nothing special about using the top and right groupoids to do this, and by doing the same to the transpose of $S$ we produce the coisotropic submanifold
\[
Y := \wt r_P^{-1}(1^{P^*}M) \subseteq S
\]
of $S$, which also contains the core. We can then ask about the reductions of these two coisotropics by their characteristic foliations.

The following observation was a key motivation:

\begin{ex}
Consider the symplectic double groupoid of Example \autoref{ex:main}. The submanifold $X$ in this case is the restricted cotangent bundle $T^*G|_M \subseteq T^*G$. As described in Chapter 2, the reduction of this coisotropic by its characteristic foliation is just the core $T^*M$, and the resulting reduction relation
\[
T^*G \red T^*M
\] 
is the transpose of the cotangent lift $T^*e$ of the unit embedding $e: M \to G$ of the groupoid $G \rightrightarrows M$.

Now, performing this procedure using the left groupoid is more interesting. Recall that the source map $\wt r$ of the groupoid $T^*G \rightrightarrows A^*$ is determined by the requirement that
\[
\wt r(g,\xi)\big|_{\ker d\ell_{e(r(g))}} = (dL_g)_{e(r(g))}^*\left(\xi\big|_{\ker d\ell_g}\right)
\]
where $L_g: \ell^{-1}(r(g)) \to \ell^{-1}(\ell(g))$ is left groupoid multiplication by $g$ and we identify $A^* \cong N^*M$. Thus, we see that $(g,\xi)$ maps to a unit $(r(g),0) \in A^*$ of the bottom groupoid (note that the unit of the bottom groupoid is given by the zero section $0_M$ of $A^*$) if and only if $\xi|_{\ker d\ell_g} = 0$, and hence equivalently if and only if $\xi$ is in the image of $d\ell_g^*$. 

Therefore, the coisotropic $Y := \wt r^{-1}(0_M)$ in $T^*G$ is equal to $N^*\F_\ell$, where $\F_\ell$ is the foliation of $G$ given by the $\ell$-fibers of the groupoid $G \rightrightarrows M$. According to Chapter 3, the reduction of this is then also the core $T^*M$, where the reduction relation $T^*G \red T^*M$ turns out to be the cotangent lift $T^*\ell$. 
\end{ex}

These results generalize in the following way to an arbitrary symplectic double groupoid. We first have the following explicit description of the characteristic foliation $X^\perp$ of $X := \widetilde{r}_{P^*}^{-1}(1^P M)$:

\begin{lem}\label{lem:fol}
The leaf of the foliation $X^\perp$ containing $s \in X$ is given by
\[
X^\perp_s := \{ s\ \wt\circ_{P^*}\,\wt 1^P_\lambda\ |\ \lambda \in \ell_{P^*}^{-1}(m)\},
\]
where $m = r_P(\wt r_{P^*}(s))$.
\end{lem}

\begin{proof}
To show that the characteristic foliation is as claimed, we must show that
\[
T_s (X^\perp_s) = (T_s X)^\perp
\]
where $(T_s X)^\perp$ is the symplectic orthogonal of $T_s X$. A dimension count shows that these two spaces have the same dimension, so we need only show that the former is contained in the latter. Thus we must show that if $Y \in T_s(X^\perp_s)$, then
\[
\omega_s(Y,V) = 0
\]
for any $V \in T_s X$ where $\omega$ is the symplectic form on $S$. Since
\[
T_s X = (d\wt r_{P^*})_s^{-1}\left(T_{1^P_m}(1^P M)\right),
\]
this means that $V$ is such that $(d\wt r_{P^*})_s V \in T_{1^P_m}(1^P M)$.

To do so, we use the following explicit description of $T_s(X^\perp_s)$. The elements of this are parametrized by $\ell_{P^*}^{-1}(m)$, so we have that $X^\perp_s$ is the image of the map
\[
\ell_{P^*}^{-1}(m) \to S
\]
given by the composition
\[
\lambda \mapsto \wt 1^P_\lambda \mapsto \wt L_s^{P^*} (\wt 1^P_\lambda)
\]
where $\wt L_s^{P^*}$ is left-multiplication by $s$ in the top groupoid. Taking differentials at $\lambda = 1^{P^*}_m$ then gives the explicit description of $T_s(X^\perp_s)$ we want; in particular, we can write $Y \in T_s(X^\perp_s)$ as
\begin{equation}\label{tang-dist}
Y = (d \wt L_s^{P^*})_{\wt 1^P(1^{P^*}_m)}(d \wt 1^P)_{1^{P^*}_m} W
\end{equation}
for some $W \in \ker (d\ell_{P^*})_{1^{P^*}_m}$.

First we consider the case where $s$ is a unit of the top groupoid structure, so suppose that $s = \wt 1^{P^*}(1^P_m)$. Using the splitting
\[
TG|_P = TP \oplus \ker d\wt r_{P^*}|_P,
\]
we can write $V$ above as
\[
V = (d\wt 1^{P^*})_{1^P_m}(d\wt r_{P^*})_s V + [V- (d\wt 1^{P^*})_{1^P_m}(d\wt r_{P^*})_s V],
\]
where the first term is tangent to the units of the top groupoid and the second term is in $\ker (d\wt r_{P^*})_s$, and so tangent to the $r$-fiber through $s$. Then we have
\[
\omega_s (Y,V) = \omega_s ((d \wt 1^P)_{1^{P^*}_m} W,(d\wt 1^{P^*})_{1^P_m}(d\wt r_{P^*})_s V) +  \omega_s ((d \wt 1^P)_{1^{P^*}_m} W,[V- (d\wt 1^{P^*})_{1^P_m}(d\wt r_{P^*})_s V]).
\]
Since $W \in \ker(d\ell_{P^*})_{1^{P^*}_m}$, one can check that $Y = (d \wt 1^P)_{1^{P^*}_m} W$ is tangent to the $\ell$-fiber through $s$, so the second term above vanishes since $\ell$ and $r$-fibers of a symplectic groupoid are symplectically dual to each other. By the defining property of $V$, we have $(d\wt r_{P^*})_s V = (d 1^P)_m v$ for some $v \in T_m M$. Using this and the fact that $\wt 1^{P^*} \circ 1^P = \wt 1^P \circ 1^{P^*}$ (which follows from the double groupoid compatibilities), we can write the first term above as
\[
\omega_s ((d \wt 1^P)_{1^{P^*}_m} W,(d \wt 1^P)_{1^{P^*}_m}(d 1^{P^*})_m v),
\]
which vanishes since the embedding $\wt 1^P: P^* \to S$ is lagrangian. Thus $\omega_s(Y,V) = 0$ as was to be shown. 

Now using (local) lagrangian bisections we can transport this argument to an arbitrary $s \in \wt r_{P^*}^{-1}(1^P M)$.
\end{proof}

\begin{rmk}
To emphasize, it is the entire double groupoid structure that makes this explicit description possible; such a description is not necessarily available given solely a symplectic realization $S \to P$.
\end{rmk}

In square notation, the leaf of the characteristic foliation through $s$ consists of elements of the form
\begin{center}
\begin{tikzpicture}[scale=.85]
	\node at (1,1) {$s\,\,\wt\circ_{P^*}\wt 1^P_\lambda$};

	\tikzset{font=\scriptsize};
	\draw (0,0) rectangle (2,2);
	\node at (-1.2,1) {$\wt\ell_P(s) \circ_{P^*} \lambda$};
	\node at (1,2.4) {$1^P_{r_{P^*}(\lambda)}$};
	\node at (3.2,1) {$\wt r_P(s) \circ_{P^*} \lambda$};
	\node at (1,-.4) {$\wt\ell_{P^*}(s)$};
\end{tikzpicture}
\end{center}
where $\lambda \in \ell_{P^*}^{-1}(m)$. Switching the roles of the top and left groupoids, we get that the leaves of the characteristic foliation of $Y := \wt r_P^{-1}(1^{P^*}M)$ are given by
\[
Y^\perp_s := \{s\ \wt\circ_P\,\wt 1^{P^*}_{\lambda}\ |\ \lambda \in \ell_P^{-1}(m)\},
\]
where $m = r_{P^*}(\wt r_P(s))$

Now, returning to the characteristic leaves of $X$, note that such a leaf intersects the core in exactly one point, since there is only one choice of $\lambda$ which will make $\wt r_P(s) \circ_{P^*} \lambda$ a unit, namely $\lambda = i_{P^*}(\wt r_P(s))$. Thus, the core forms a cross section to the characteristic foliation of $X$ and we conclude that the leaf space $X/X^\perp$ of this characteristic foliation can be identified with the core. This leaf space is then naturally symplectic, and we have:

\begin{prop}
The symplectic structure on the core obtained via the above reduction agrees with Mackenzie's symplectic structure.
\end{prop}

\begin{proof}
First, let $i: C \to S$ be the inclusion of the core into $S$. Then if $\omega$ is the symplectic form on $S$, $i^*\omega$ is Mackenzie's symplectic structure on $C$.

Now, let $\pi: X \to C$ be the surjective submersion sending $s \in X$ to the characteristic leaf containing it, where we have identified the leaf space $X/X^\perp$ with $C$ in the above manner. Let $s: C \to X$ be the inclusion of the core into $X$; this is a section of $\pi$. Finally, let $j: X \to S$ be the inclusion of $X$ into $S$.

The symplectic form $\omega_C$ on $C$ obtained via reduction is characterized by the property that $j^*\omega = \pi^*\omega_C$. Since $i = j \circ s$, we have
\begin{align*}
i^*\omega &= (j \circ s)^*\omega \\
&= s^*(j^*\omega) \\
&= s^*(\pi^*\omega_C) \\
&= (\pi \circ s)^*\omega_C,
\end{align*}
which equals $\omega_C$ since $s$ is a section of $\pi$. This proves the claim.
\end{proof}

Explicitly, the reduction relation $S \red C$, which we will call $E^t$ for reasons to be made clear later, obtained by reducing $X$ is given by
\begin{displaymath}
E^t: s \mapsto s\, \wt\circ_{P^*} \wt 1^P_{i_{P^*}(\wt r_P(s))} \text{ for } s \in \wt r_{P^*}^{-1}(1^P M).
\end{displaymath}
For future reference, the transposed relation $E: C \cored S$ (which is a coreduction) is given by
\begin{equation}\label{E}
E: s \mapsto s \wt\circ_{P^*} \wt 1^P_\lambda, \text{ for $\lambda \in P^*$ such that } \ell_{P^*}(\lambda) = r_P(\wt r_{P^*}(s)).
\end{equation}

The same then holds if we consider the tranposed double groupoid structure, so that the reduction of the coisotropic $Y = \wt r_P^{-1}(1^{P^*}M) \subseteq S$ can also be identified with the core via the same maps, simply exchanging the roles of $P$ and $P^*$. The reduction relation $S \red C$ obtained by reducing $Y$ will be called $L$ and is explicitly given by
\begin{equation}\label{L}
L: s \mapsto s\,\wt\circ_P\,\wt 1^{P^*}_{i_{P}(\wt r_{P^*}(s))} \text{ for } s \in \wt r_P^{-1}(1^{P^*} M).
\end{equation}

Similar results hold for the preimages of units under the target maps. To be clear, let $Z$ now be $\wt\ell_{P^*}^{-1}(1^P M)$, the preimage of the units of $P$ under the target map of the top groupoid. This is again coisotropic in $S$, and the leaf of the characteristic foliation through a point $s \in Z$ is now given by
\[
\{\wt 1^P_{\lambda}\,\wt\circ_{P^*}\,s\ |\ \lambda \in r_{P^*}^{-1}(m)\}
\]
where $m = \ell_P(\wt\ell_{P^*}(s))$. Similar to the above, we can now easily identify the reduction of $Z$ with the set of elements of $S$ of the form
\begin{center}
\begin{tikzpicture}[scale=.75]
	\node at (1,1) {$s$};

	\tikzset{font=\scriptsize};
	\draw (0,0) rectangle (2,2);
	\node at (-.6,1) {$1^{P^*}_m$};
	\node at (1,2.4) {$\wt r_{P^*}(s)$};
	\node at (2.7,1) {$\wt r_P(s)$};
	\node at (1,-.4) {$1^P_m$};
\end{tikzpicture}
\end{center}
which we might call the ``left-core'' $C_L$ of $S$ to distinguish it from the ``right-core'' $C_R\ (:= C)$ previously defined. However, the left-core $C_L$ can be identified with $C_R$ using the composition of the two groupoid inverses on $S$:
\[
\wt i_P \circ \wt i_{P^*}: C_L \to C_R,
\]
so we again get that the reduction of $Z$ is symplectomorphic to the core $C$. Considering the transpose of $S$, we find that the same is true for the preimage of the units of $P^*$ under the target of the left groupoid.

\begin{ex}
Let us return to Example \ref{ex:main}. The target of the top groupoid is the same as the source, so the above reduction procedure produces the same reduction relation
\[
(T^*e)^t: T^*G \red T^*M
\]
as before. A similar computation to that carried out for the source of the left groupoid $T^*G \rightrightarrows A^*$ shows that the coisotropic submanifold $\wt \ell^{-1}(0_M)$ of $T^*G$ is $N^*\F_r$, where $\F_r$ is the foliation of $G$ given by the $r$-fibers of $G \rightrightarrows M$, and that the reduction relation
\[
T^*G \red T^*M
\]
obtained by reducing $N^*\F_r$ is then the cotangent lift $T^*r$.
\end{ex}

We thus have multiple ways of recovering the core $S$ by reducing certain coisotropic submanifolds, and each such way produces a canonical relation $S \red C$. Moreoever, we can recover the core groupoid product on $C$ via this reduction procedure:

\begin{prop}
Let $M: S \times S \to S$ be the canonical relation corresponding to the product $\wt m_P$ of the left groupoid structure on $S$. Then the composition
\begin{center}
\begin{tikzpicture}[>=angle 90]
	\node (UL) at (0,0) {$C \times C$};
	\node (UR) at (3,0) {$S \times S$};
	\node (LL) at (6,0) {$S$};
	\node (LR) at (8,0) {$C$};

	\tikzset{font=\scriptsize};
	\draw[->] (UL) to node [above] {$E \times E$} (UR);
	\draw[->] (UR) to node [above] {$M$} (LL);
	\draw[->] (LL) to node [above] {$L$} (LR);
\end{tikzpicture}
\end{center}
is Mackenzie's groupoid product on $C$.
\end{prop}

\begin{proof}
Let $s, s' \in C$. Then applying the relation $E \times E$ gives
\[
(s,s') \mapsto (s\,\wt\circ_{P^*}\,\wt 1^P_\lambda, s'\,\wt\circ_{P^*}\,\wt 1^P_{\lambda'}).
\]
Now, these are composable under $\wt m_P$ when
\[
\wt r_P(s\,\wt\circ_{P^*}\,\wt 1^P_\lambda) = \wt r_P(s) \circ_{P^*} \lambda = \lambda \text{ equals } \wt \ell_P(s'\,\wt\circ_{P^*}\,\wt 1^P_{\lambda'}) = \wt\ell_P(s') \circ_{P^*} \lambda',
\]
where we have used the fact that $\wt r_P(s)$ is a unit. From this we get the condition that
\[
\lambda = \wt\ell_P(s') \circ_{P^*}\lambda'.
\]
Then the relation $M$ produces
\[
(s\,\wt\circ_{P^*}\,\wt 1^P_\lambda) \wt\circ_P (s'\,\wt\circ_{P^*}\,\wt 1^P_{\lambda'}) = \left(s\,\wt\circ_{P^*}\,\wt 1^P_{\wt\ell_P(s') \circ_{P^*}\lambda'}\right) \wt\circ_P (s'\,\wt\circ_{P^*}\,\wt 1^P_{\lambda'}).
\]
Now, to be in the domain of the final relation $L$, we need the source $\wt r_P$ of the above to be a unit, but $\wt r_P$ of the above is
\[
\wt r_P(s') \circ_{P^*} \lambda' = \lambda'
\]
since $\wt r_P(s')$ is a unit. Thus our expression is in the domain of the final relation exactly when $\lambda'$ is a unit, and the final relation gives
\[
\left(s \wt\circ_{P^*} \wt 1^P_{\wt\ell_P(s')}\right) \wt\circ_P s',
\]
which is the product of $s$ and $s'$ under Mackenzie's groupoid structure on $C$.
\end{proof}

Indeed, using the characterization of symplectic groupoids in terms of symplectic monoids in the symplectic category, we can now give an alternate proof that $C$ with the above product is a symplectic groupoid as follows.

The unit $e: pt \to C$ is obtained as the composition
\begin{center}
\begin{tikzpicture}[>=angle 90]
	\node (UL) at (0,1) {$pt$};
	\node (UM) at (2,1) {$S$};
	\node (UR) at (4,1) {$C,$};

	\tikzset{font=\scriptsize};
	\draw[->] (UL) to node [above] {$P^*$} (UM);
	\draw[->] (UM) to node [above] {$L$} (UR);
\end{tikzpicture}
\end{center}
where $P^*$ denotes the image of the lagrangian embedding $\wt 1^P: P^* \to S$. Now, the left unit property of $e$ follows from the commutativity of the diagram
\begin{center}
\begin{tikzpicture}[>=angle 90]
	\node (L1) at (0,0) {$C$};
	\node (L2) at (3,0) {$S$};
	\node (U1) at (0,2) {$C \times C$};
	\node (U2) at (3,2) {$S \times S$};
	\node (U3) at (6,1) {$S$};
	\node (U4) at (9,1) {$C$};

	\tikzset{font=\scriptsize};
	\draw[->] (U1) to node [above] {$E \times E$} (U2);
	\draw[->] (U2) to node [above] {$M$} (U3);
	\draw[->] (U3) to node [above] {$L$} (U4);
	\draw[->] (L1) to node [above] {$E$} (L2);
	\draw[->] (L1) to node [left] {$e \times id$} (U1);
	\draw[->] (L2) to node [left] {$P^* \times id$} (U2);
	\draw[->] (L2) to node [below] {$id$} (U3);
	\draw[->] (U1) to [bend left=30] node [above] {$m$} (U4);
	\draw[->] (L1) to [bend right=20] node [below] {$id$} (U4);
\end{tikzpicture}
\end{center}
where $m$ is Mackenzie's groupoid product on $C$ and the composition $L \circ id \circ E = L \circ E$ is the identity on $C$ by Proposition \ref{prop:sections}, which we prove in the coming sections. A similar diagram gives the right unit property of $e$. Finally, the $*$-structure $s: \overline{C} \to C$ on $(C,m,e)$ is given by the composition
\begin{center}
\begin{tikzpicture}[>=angle 90]
	\node (1) at (0,1) {$\overline{C}$};
	\node (2) at (2,1) {$\overline{S}$};
	\node (3) at (4,1) {$S$};
	\node (4) at (6,1) {$C,$};

	\tikzset{font=\scriptsize};
	\draw[->] (1) to node [above] {$E$} (2);
	\draw[->] (2) to node [above] {$\wt i_P$} (3);
	\draw[->] (3) to node [above] {$E^t$} (4);
\end{tikzpicture}
\end{center}
and the strong positivity of this $*$-structure follows from the corresponding diagram for $S$ with the left symplectic groupoid structure.

In summary, we have shown:
\begin{thrm}\label{thrm:main}
Let $S$ be a symplectic double groupoid with core $C$. Then the reductions of the coisotropic submanifolds
\[
\wt r_{P^*}^{-1}(1^P M),\ \wt r_P^{-1}(1^{P^*}M),\ \wt \ell_{P^*}^{-1}(1^P M), \text{ and } \wt \ell_P^{-1}(1^{P^*}M)
\]
of $S$ are symplectomorphic to the core $C$ of $S$. Furthermore, under these reductions the groupoid structures on $S$ descend to Mackenzie's symplectic groupoid structure on $C$.
\end{thrm}

\section{Symplectic Double Groups}
It will be instructive to see the results of the above constructions in the special case of a symplectic double group. We will see that we recover a result of Zakrzewski \cite{SZ2} concerning Hopf algebra objects in the symplectic category.

Suppose that $S$ is a symplectic double group, so that the double base is a point. In particular, then $P$ and $P^*$ are Poisson-Lie groups. The coisotropic $X \subseteq S$ in this case is simply the fiber of $\wt r_{P^*}$ over the identity element of $P$. This fiber can be identified with $P^*$ and is lagrangian. Hence its reduction is a point, as the core of a symplectic double group should be. The reduction relation
\[
S \to pt
\]
resulting from this is the lagrangian submanifold of $S$ given by the unit embedding of $P^*$ into $S$.

The transpose of the above relation, together with the graph of the left groupoid product $S \times S \to S$ give $S$ the structure of a symplectic monoid. The compatability between this structure and the comonoid structure resulting from the top groupoid structure on $S$ can then be expressed by simply saying that $S$ together with these structures is a Hopf algebra object in $\Symp$, where the antipode $S \to S$ is the composition of the two groupoid inverses of $S$; this is then a generalization of Example \ref{ex:hopf}: In fact, we have the following:

\begin{thrm}[Zakrzewski, \cite{SZ2}]
$S$ endowed with two symplectic groupoid structures is a symplectic double group if and only if $S$ endowed with the above mentioned monoid and comonoid structures forms a Hopf algebra object in $\Symp$.
\end{thrm}

To be clear, what we call here a ``Hopf algebra object'' in the symplectic category is what Zakzrewski calls an ``$S^*$-group'', and his result is phrased in a slightly different manner.

Most of the rest of this chapter is concerned with finding an analog of this result for general symplectic double groupoids. Indeed, if $S$ is a symplectic double groupoid, we naturally have induced symplectic monoid and comonoid structures on $S$, and we can ask whether these satisfy the same Hopf-condition as above. As already mentioned in the previous chapter, a simple computation shows that this is the case only when $S$ is a symplectic double group---the existence of a double base which is not a point prevents this in general. What we consider instead then is a ``Hopf algebroid''-like structure, which should be thought of as similar to a Hopf algebra object except with a more general base.

\section{Double Structures $\to$ Symplectic Category}
Let us return to the structures of Example \autoref{ex:main}. To recall, performing the above procedure in this case produced the relations $(T^*e)^t: T^*G \to T^*M$, $T^*r: T^*G \to T^*M$, and $T^*\ell: T^*G \to T^*r$ obtained by taking the cotangent lifts of the groupoid structure maps of $G \rightrightarrows M$. Also note that the cotangent lift $T^*i: T^*G \to T^*G$ of the groupoid inverse of $G$ is just the composition of the two symplectic groupoid inverses on $T^*G$. In addition, the two canonical relations $T^*m: T^*G \times T^*G \to T^*G$ and $T^*\Delta: T^*G \to T^*G \times T^*G$ (where $\Delta: G \to G \times G$ is the usual diagonal map) come from the two symplectic groupoid products on $T^*G$.

The above canonical relations are precisely the ones which appear in the diagrams obtained by applying $T^*$ to the commutative diagrams appearing in the definition of a Lie groupoid. This suggests that the structure $T^*G \rightrightarrows T^*M$ obtained should be viewed as the analog of a ``groupoid'' in the symplectic category, as mentioned before. Since the ``diagonal'' $T^*\Delta$ is extra required data, we prefer to think of it as a coproduct and view $T^*G \rightrightarrows T^*M$ as the analog of a ``Hopf algebroid'' in the symplectic category.

This all generalizes in the following way for an arbitrary symplectic double groupoid; we will denote the resulting structure by $S \rightrightarrows C$. First, the coreduction $E: C \to S$ is given in \ref{E}; this is the analog of $T^*e: T^*M \cored T^*G$ in Example \ref{ex:main}. Second, the reduction $L: S \to C$ is given in \ref{L}; this is the analog of $T^*\ell: T^*G \red T^*M$ in Example \ref{ex:main}. We will think of $E$ as the ``unit'' and $L$ the ``target'' of $S \rightrightarrows C$.

Now, as we noted before, the reduction of $\wt \ell_P^{-1}(1^{P^*}M)$ naturally gives the left-core $C_L$, so to get a morphism to $C$ we must post-compose the resulting reduction relation $S \red C_L$ with the composition $\wt i_P \circ \wt i_{P^*}$ of the two inverses on $S$. Note that since the inverse of a symplectic groupoid is an anti-symplectomorphism, this composition of two inverses is a symplectomorphism, so the relation $R: S \red C$ so obtained is indeed a canonical relation. The relation $R: S \red C_L \cong C$ so obtained is explicitly given by
\begin{equation}\label{R}
R: s \mapsto \wt i_P \wt i_{P^*}(s)\ \wt\circ_P\ \wt 1^{P^*}_{\wt \ell_{P^*}(s)}.
\end{equation}
In Example \ref{ex:main} this relation became $T^*r: T^*G \red T^*M$, and we will think of $R$ as the ``source'' of $S \rightrightarrows C$.

Finally, we will let $I: S \to S$ be the composition $\wt i_P \circ \wt i_{P^*}$, which as noted above is a symplectomorphism; in Example \ref{ex:main} this is $T^*i$. We then have the following observations:

\begin{prop}\label{prop:sections}
The canonical relations above satisfy the following identities: $L \circ E = id_C$, $R \circ E = id_C$, $R = L \circ I$, $L = R \circ I$. All of these compositions are strongly transversal.
\end{prop}

\begin{proof}
Suppose that $s \in C$, so that $\wt r_P(s)$ and $\wt r_{P^*}(s)$ are both units. Then the relation $E$ sends this to 
\[
E: s \mapsto s\,\wt\circ_{P^*}\,\wt 1^P_\lambda
\]
for $\lambda$ in the same leaf as $s$ of the characteristic foliation of $\wt r_{P^*}^{-1}(1^P M)$. Now, this is in the domain of $L$ when
\[
\wt r_P(s\,\wt\circ_{P^*}\,\wt 1^P_\lambda) = \wt r_P(s)\,\circ_{P^*}\,\lambda = \lambda
\]
is a unit. Thus there is only such $\lambda$ so that $s\,\wt\circ_{P^*}\,\wt 1^P_\lambda \in \dom L$, from which strong transversality will follow, and then it is straightforward to check that applying the relation $L$ will give back $s$. Hence $L \circ E = id_C$.

More interestingly, $s\,\wt\circ_{P^*}\,\wt 1^P_\lambda$ is in the domain of $R$ when
\[
\wt \ell_P(s\,\wt\circ_{P^*}\,\wt 1^P_\lambda) = \wt \ell_P(s)\circ_{P^*}\lambda
\]
is a unit, which requires that $\lambda = i_{P^*}(\wt\ell_P(s))$. Again, from this strong transversality follows and we see that the composition $E \circ R$ is
\[
s \mapsto \wt i_P \wt i_{P^*}\left(s\,\wt\circ_{P^*}\,\wt 1^P_{i_{P^*}(\wt\ell_P(s))}\right) \wt\circ_P\,\wt 1^{P^*}_{\wt \ell_{P^*}(s)} = \left(\wt 1^P_{\wt\ell_P(s)}\,\wt\circ_{P^*}\,\wt i_P\wt i_{P^*}(s)\right)\wt\circ_P\,\wt 1^{P^*}_{\wt\ell_{P^*}(s)}.
\]

We claim that this result is simply $s$. To see this, we express the result as

\begin{center}
\begin{tikzpicture}[>=angle 90]
	\draw (0,0) rectangle (2,2);
	\node at (1,1) {$\wt 1^P_{\wt\ell_P(s)}$};
	\draw (0,2) rectangle (2,4);
	\node at (1,3) {$\wt i_P\wt i_{P^*}(s)$};
	\draw (2,0) rectangle (4,4);
	\node at (3,2) {$\wt 1^{P^*}_{\wt\ell_{P^*}(s)}$};
	
	\node at (5,2) {$=$};
	
	\draw (6,0) rectangle (8,2);
	\node at (7,1) {$\wt 1^P_{\wt\ell_P(s)}$};
	\draw (8,0) rectangle (10,2);
	\node at (9,1) {$s$};
	\draw (6,2) rectangle (8,4);
	\node at (7,3) {$\wt i_P\wt i_{P^*}(s)$};
	\draw (8,2) rectangle (10,4);
	\node at (9,3) {$\wt i_{P^*}(s)$};
	
	\node at (11,2) {$=$};
	
	\draw (12,0) rectangle (14,2);
	\node at (13,1) {$s$};
	\draw (12,2) rectangle (14,4);
	\node at (13,3) {$\wt 1^P_{i_{P^*}(\wt r_P(s))}$};
\end{tikzpicture}
\end{center}
where in the first step we decompose vertically and in the second we compose horizontally. Now, since $\wt r_P(s)$ is a unit, the element on top in the last term is of the form
\[
\wt 1^P_{i_{P^*}(1^{P^*}_m)} = \wt 1^P 1^{P^*}_m = \wt 1^{P^*}1^P_m.
\]
The claim then follows by composing the last term vertically. The computations which show that $R = L \circ I$ and $L = R \circ I$ are similar.
\end{proof}

Now, there are two groupoid multiplications on $S$. As before we will denote by $M: S \times S \to S$ the canonical relation given by $\wt m_P$, and will think of $M$ as a ``product''. Similarly, we will denote by $\Delta: S \to S \times S$ the transpose of canonical relation obtained from $\wt m_{P^*}$, and will think of $\Delta$ as a ``coproduct''.

The following propositions then express some compatibilities between these relations and those previously defined. In particular, in the example of $T^*G \rightrightarrows T^*M$, the diagrams considered in the next two are precisely the ones obtained by applying the cotangent functor to those in condition (iii, left and right units) defining a Lie groupoid.

\begin{prop}
The following diagram commutes and all compositions are strongly transversal:
\begin{center}
\begin{tikzpicture}[>=angle 90]
	\node (UL) at (0,1) {$S \times S$};
	\node (UM) at (3,1) {$C \times S$};
	\node (UR) at (6,1) {$S \times S$};
	\node (LL) at (0,-1) {$S$};
	\node (LR) at (6,-1) {$S$};

	\tikzset{font=\scriptsize};
	\draw[->] (UL) to node [above] {$L \times id$} (UM);
	\draw[->] (UM) to node [above] {$E \times id$} (UR);
	\draw[->] (LL) to node [left] {$\Delta$} (UL);
	\draw[->] (UR) to node [right] {$M$} (LR);
	\draw[->] (LL) to node [above] {$id$} (LR);
\end{tikzpicture}
\end{center}
\end{prop}

\begin{proof}
For $s \in S$, the composition $M \circ (E \times id) \circ (L \times id) \circ \Delta$ is
\[
s \mapsto [(s_1\ \wt\circ_P\ \wt 1^{P^*}_{i_P(\wt r_{P^*}(s_1))})\ \wt\circ_{P^*}\ \wt 1^P_\lambda]\ \wt\circ_P\ s_2
\]
where $s = s_1\ \wt\circ_{P^*}\ s_2$ for composable---with respect to the groupoid structure whose product is $\Delta^t$---$s_1$ and $s_2 \in S$; this is simply saying that $s$ is in the domain of $\Delta$. We must show that the resulting expression is just $s$ itself. This follows from the compositions:

\begin{center}
\begin{tikzpicture}
	\draw (0,0) rectangle (2,2);
	\node at (1,1) {$s_1$};
	\draw (2,0) rectangle (4,2);
	\node at (3,1) {$\,\wt 1^{P^*}_{i_P(\wt r_{P^*}(s_1))}$};
	\draw (0,2) rectangle (4,4);
	\node at (2,3) {$\wt 1^P_\lambda$};
	\draw (4,0) rectangle (6,4);
	\node at (5,2) {$s_2$};
	
	\node at (7,2) {$=$};
	
	\draw (8,0) rectangle (10,2);
	\node at (9,1) {$s_1$};
	\draw (10,0) rectangle (12,2);
	\node at (11,1) {$\,\wt 1^{P^*}_{i_P(\wt r_{P^*}(s_1))}$};
	\draw (12,0) rectangle (14,2);
	\node at (13,1) {$\wt 1^{P^*}_{\wt \ell_{P^*}(s_2)}$};
	\draw (8,2) rectangle (10,4);
	\node at (9,3) {$s_2$};
	\draw (10,2) rectangle (12,4);
	\node at (11,3) {$\wt i_P(s_2)$};
	\draw (12,2) rectangle (14,4);
	\node at (13,3) {$s_2$};
\end{tikzpicture}
\end{center}
where we decompose the top left box horizontally and the right box vertically, and

\begin{center}
\begin{tikzpicture}[scale=.91]
	\draw (8,0) rectangle (10,2);
	\node at (9,1) {$s_1$};
	\draw (10,0) rectangle (12,2);
	\node at (11,1) {$\,\wt 1^{P^*}_{i_P(\wt r_{P^*}(s_1))}$};
	\draw (12,0) rectangle (14,2);
	\node at (13,1) {$\wt 1^{P^*}_{\wt\ell_{P^*}(s_2)}$};
	\draw (8,2) rectangle (10,4);
	\node at (9,3) {$s_2$};
	\draw (10,2) rectangle (12,4);
	\node at (11,3) {$\wt i_P(s_2)$};
	\draw (12,2) rectangle (14,4);
	\node at (13,3) {$s_2$};
	
	\node at (15,2) {$=$};
	
	\draw (16,0) rectangle (18,2);
	\node at (17,1) {$s_1$};
	\draw (18,0) rectangle (20,2);
	\node at (19,1) {$\,\wt 1^{P^*}_{r_P(\wt\ell_{P^*}(s_2))}$};
	\draw (16,2) rectangle (18,4);
	\node at (17,3) {$s_2$};
	\draw (18,2) rectangle (20,4);
	\node at (19,3) {$\wt 1^P_{\wt r_P(s_2)}$};
	
	\node at (21,2) {$=$};
	
	\draw (22,0) rectangle (24,2);
	\node at (23,1) {$s_1$};
	\draw (22,2) rectangle (24,4);
	\node at (23,3) {$s_2$};
	\draw (24,0) rectangle (26,4);
	\node at (25,2) {$\wt 1^P_{\wt r_P(s_2)}$};
\end{tikzpicture}
\end{center}
where we first compose the right four boxes horizontally (using the fact that $\wt r_{P^*}(s_1) = \wt\ell_{P^*}(s_2)$) and then compose the right two boxes vertically. The resulting composition is then $s_1\,\wt\circ_{P^*}\,s_2 = s$ as was to be shown.
\end{proof}

The proof of the following result is very similar to that above:

\begin{prop}
The following diagram commutes and all compositions are strongly transversal:
\begin{center}
\begin{tikzpicture}[>=angle 90]
	\node (UL) at (0,1) {$S \times S$};
	\node (UM) at (3,1) {$S \times C$};
	\node (UR) at (6,1) {$S \times S$};
	\node (LL) at (0,-1) {$S$};
	\node (LR) at (6,-1) {$S$};

	\tikzset{font=\scriptsize};
	\draw[->] (UL) to node [above] {$id \times R$} (UM);
	\draw[->] (UM) to node [above] {$id \times (I \circ E)$} (UR);
	\draw[->] (LL) to node [left] {$\Delta$} (UL);
	\draw[->] (UR) to node [right] {$M$} (LR);
	\draw[->] (LL) to node [above] {$id$} (LR);
\end{tikzpicture}
\end{center}
\end{prop}

In particular, these results already imply that there is a simplicial symplectic manifold structure on $S \rightrightarrows C$:

\begin{cor}
There is a simplicial symplectic manifold
\begin{center}
\begin{tikzpicture}[>=angle 90]
	\def\A{.3};
	\def\B{2.2};
	\def\C{3.3};

	\node at (0,0) {$\cdots$};
	\node at (1.5,0) {$S \times S$};
	\node at (3.0,0) {$S$};
	\node at (4.1,0) {$C,$};

	\tikzset{font=\scriptsize};
	\draw[->] (\A,+.18) -- (\A+.5,0+.18);
	\draw[->] (\A,+.06) -- (\A+.5,+.06);
	\draw[->] (\A,-.06) -- (\A+.5,-.06);
	\draw[->] (\A,-.18) -- (\A+.5,0-.18);
	
	\draw[->] (\B,+.12) -- (\B+.5,+.12);
	\draw[->] (\B,0) -- (\B+.5,0);
	\draw[->] (\B,-.12) -- (\B+.5,-.12);
	
	\draw[->] (\C,0+.06) -- (\C+.5,0+.06);
	\draw[->] (\C,0-.06) -- (\C+.5,0-.06);
\end{tikzpicture}
\end{center}
generalizing that of Example~\ref{ex:cot-simp-symp}, where the degree $n$ piece is the product of $n$ copies of $S$. Here the degree $1$ face morphisms are $L$ and $R$, the degree $2$ face morphisms are ``projection onto the first factor''\footnote{Here and below, the term ``projection'' is motivated by Example~\ref{ex:cot-simp-symp}, where these relations are indeed cotangent lifts of projections. The expression for these morphisms is also motivated by that example.}:
\begin{center}
\begin{tikzpicture}[>=angle 90]
	\node (1) at (0,1) {$S \times S$};
	\node (2) at (3,1) {$S \times C$};
	\node (3) at (6,1) {$S \times S$};
	\node (4) at (8,1) {$S,$};

	\tikzset{font=\scriptsize};
	\draw[->] (1) to node [above] {$id \times L$} (2);
	\draw[->] (2) to node [above] {$id \times R^t$} (3);
	\draw[->] (3) to node [above] {$\Delta^t$} (4);
\end{tikzpicture}
\end{center}
the product $M$, and ``projection onto the second factor'':
\begin{center}
\begin{tikzpicture}[>=angle 90]
	\node (1) at (0,1) {$S \times S$};
	\node (2) at (3,1) {$S \times C$};
	\node (3) at (6,1) {$S \times S$};
	\node (4) at (8,1) {$S.$};

	\tikzset{font=\scriptsize};
	\draw[->] (1) to node [above] {$R \times id$} (2);
	\draw[->] (2) to node [above] {$L^t \times id$} (3);
	\draw[->] (3) to node [above] {$\Delta^t$} (4);
\end{tikzpicture}
\end{center}
The degree $0$ degeneracy morphism is $E$, and the degree $1$ degeneracy morphisms are obtained by composing along the left and top sides of the diagrams in the previous propositions.
\end{cor}

\begin{proof}
The remaining face and degenaracy morphisms can be guessed from these. The previous propositions give part of the simplicial identities, and the rest can be checked via a straightforward, if tedious, computation in a similar manner.
\end{proof}

The following brings in an antipode-like condition on $I$. In particular, in the case of $T^*G \rightrightarrows T^*M$, these diagrams are those obtained by applying $T^*$ to the diagrams in condition (iv, left and right inverse) in Definition~\ref{grpd} of a Lie groupoid:

\begin{prop}
The following diagram commutes and all compositions are strongly transversal:
\begin{center}
\begin{tikzpicture}[>=angle 90]
	\node (UL) at (0,1) {$S \times S$};
	\node (UR) at (6,1) {$S \times S$};
	\node (LL) at (0,-1) {$S$};
	\node (LM) at (3,-1) {$C$};
	\node (LR) at (6,-1) {$S$};

	\tikzset{font=\scriptsize};
	\draw[->] (UL) to node [above] {$id \times I$} (UR);
	\draw[->] (LL) to node [left] {$\Delta$} (UL);
	\draw[->] (UR) to node [right] {$M$} (LR);
	\draw[->] (LL) to node [above] {$L$} (LM);
	\draw[->] (LM) to node [above] {$E$} (LR);
\end{tikzpicture}
\end{center}
\end{prop}

\begin{proof}
Let $s \in S$. The composition $M \circ (id \times I) \circ \Delta$ then looks like
\[
s \mapsto s_1\ \wt\circ_P\ \wt i_P \wt i_{P^*}(s_2)
\]
where $s = s_1 \wt\circ_{P^*} s_2$. The composition $E \circ T$ is
\[
s \mapsto (s\ \wt\circ_P\ \wt 1^{P^*}_{i_P(\wt r_{P^*}(s))})\ \wt\circ_{P^*}\ \wt 1^P_\lambda.
\]
We must show that anything of the form resulting from the first relation is equivalently of the form resulting from the second. In particular, writing $s$ as $s = s_1 \wt\circ_{P^*} s_2$, we must show that
\[
s_1\ \wt\circ_P\ \wt i_P \wt i_{P^*}(s_2) =  [(s_1 \wt\circ_{P^*} s_2)\ \wt\circ_P\ \wt 1^{P^*}_{i_P(\wt r_{P^*}(s_2))}]\ \wt\circ_{P^*}\ \wt 1^P_\lambda.
\]

For this we proceed as follows. The expression on the right is

\begin{center}
\begin{tikzpicture}
	\draw (0,0) rectangle (2,2);
	\node at (1,1) {$s_1$};
	\draw (0,2) rectangle (2,4);
	\node at (1,3) {$s_2$};
	\draw (2,0) rectangle (4,4);
	\node at (3,2) {$\,\wt 1^{P^*}_{i_P(\wt r_{P^*}(s_2))}$};
	\draw (0,4) rectangle (4,6);
	\node at (2,5) {$\wt 1^P_\lambda$};
	
	\node at (5,3) {$=$};
	
	\draw (6,0) rectangle (8,2);
	\node at (7,1) {$s_1$};
	\draw (6,2) rectangle (8,4);
	\node at (7,3) {$s_2$};
	\draw (8,0) rectangle (10,2);
	\node at (9,1) {$\,\wt 1^{P^*}_{i_P(\wt r_{P^*}(s_2))}$};
	\draw (8,2) rectangle (10,4);
	\node at (9,3) {$\,\wt 1^{P^*}_{i_P(\wt r_{P^*}(s_2))}$};
	\draw (6,4) rectangle (8,6);
	\node at (7,5) {$\wt i_{P^*}(s_2)$};
	\draw (8,4) rectangle (10,6);
	\node at (9,5) {$\wt i_P\wt i_{P^*}(s_2)$};
	
	\node at (11,3) {$=$};
	
	\draw (12,1) rectangle (14,3);
	\node at (13,2) {$s_1$};
	\draw (14,1) rectangle (16,3);
	\node at (15,2) {$\,\wt 1^{P^*}_{i_P(\wt r_{P^*}(s_2))}$};
	\draw (12,3) rectangle (14,5);
	\node at (13,4) {$\wt 1^{P^*}_{\wt\ell_{P^*}(s_2)}$};
	\draw (14,3) rectangle (16,5);
	\node at (15,4) {$\wt i_P\wt i_{P^*}(s_2)$};
\end{tikzpicture}
\end{center}
where in the first step we have decomposed the right box vertically and the top box horizontally, and in the second we have composed the top four boxes vertically. The desired equality now follows by composing the remaining boxes vertically.
\end{proof}

Again, the proof of the following is very similar to that of the above proposition.

\begin{prop}
The following diagram commutes and all compositions are strongly transversal:
\begin{center}
\begin{tikzpicture}[>=angle 90]
	\node (UL) at (0,1) {$S \times S$};
	\node (UR) at (6,1) {$S \times S$};
	\node (LL) at (0,-1) {$S$};
	\node (LM) at (3,-1) {$C$};
	\node (LR) at (6,-1) {$S$};

	\tikzset{font=\scriptsize};
	\draw[->] (UL) to node [above] {$I \times id$} (UR);
	\draw[->] (LL) to node [left] {$\Delta$} (UL);
	\draw[->] (UR) to node [right] {$M$} (LR);
	\draw[->] (LL) to node [above] {$R$} (LM);
	\draw[->] (LM) to node [above] {$I \circ E$} (LR);
\end{tikzpicture}
\end{center}
\end{prop}

\begin{rmk}
Returning to the simplicial point of view, one can check that $I: S \to S$ extends to a $*$-structure on the simplicial symplectic manifold resulting from $S \rightrightarrows C$. In particular, the $*$-structure in degree $2$ is the canonical relation $S \times S \to S \times S$ given by $(I \times I) \times \sigma$, where $\sigma: S \times S \to S \times S$ is the symplectomorphism switching components.
\end{rmk}

\section{Symplectic Hopfoids}
As stated before, the results of the previous propositions suggest that the structure $S \rightrightarrows C$ resulting from a symplectic double groupoid is similar to that of a ``groupoid''. This motivates the following definition, where the diagrams which appear---specifically the ones in conditions (vii) and (viii)---should be viewed as the analogs of the corresponding diagrams in Definition \ref{grpd} of a Lie groupoid.

\begin{defn}\label{defn:symp-hopf}
A \emph{symplectic hopfoid} $S \rightrightarrows C$ consists of the following data:
\begin{itemize}
\item strongly positive symplectic $*$-comonoids $S$ and $C$ called the total and base space respectively,
\item reductions $L, R: S \to C$ called the target and source respectively,
\item a coreduction $E: C \to S$ called the unit,
\item a canonical relation $M: S \times S \to S$ called the product, and
\item a symplectomorphism $I: S \to S$ called the antipode
\end{itemize}
satisfying the following requirements:
\begin{enumerate}[(i)]
\item $L \circ E = id_C = R \circ E$,
\item $L \circ I = R$ and $R \circ I = L$,
\item $L$ and $R$ preserve the counits of $S$ and $C$,
\item $M$ is associative,
\item $I^2 = id_S$, $I$ commutes with the $*$-structure of $S$, and the diagram
\begin{center}
\begin{tikzpicture}[>=angle 90]
	\node (U1) at (0,1) {$S \times S$};
	\node (U2) at (3,1) {$S \times S$};
	\node (U3) at (6,1) {$S \times S$};
	\node (L1) at (0,-1) {$S$};
	\node (L3) at (6,-1) {$S$};

	\tikzset{font=\scriptsize};
	\draw[->] (U1) to node [above] {$\sigma$} (U2);
	\draw[->] (U2) to node [above] {$I \times I$} (U3);
	\draw[->] (U1) to node [left] {$m$} (L1);
	\draw[->] (U3) to node [right] {$m$} (L3);
	\draw[->] (L1) to node [above] {$I$} (L3);
\end{tikzpicture}
\end{center}
where $\sigma: S \times S \to S \times S$ is the symplectomorphism exchanging components, commutes,
\item the diagram
\begin{center}
\begin{tikzpicture}[>=angle 90]
	\node (U1) at (0,1) {$S \times S$};
	\node (U2) at (3,1) {$S$};
	\node (U3) at (6,1) {$S \times S$};
	\node (L1) at (0,-1) {$S \times S \times S \times S$};
	\node (L3) at (6,-1) {$S \times S \times S \times S$};

	\tikzset{font=\scriptsize};
	\draw[->] (U1) to node [above] {$m$} (U2);
	\draw[->] (U2) to node [above] {$\Delta$} (U3);
	\draw[->] (U1) to node [left] {$\Delta \times \Delta$} (L1);
	\draw[->] (L3) to node [right] {$m \times m$} (U3);
	\draw[->] (L1) to node [above] {$id \times \sigma \times id$} (L3);
\end{tikzpicture}
\end{center}
where $\sigma: S \times S \to S \times S$ is the symplectomorphism exchanging components, commutes,
\item the diagrams
\begin{center}
\begin{tikzpicture}[>=angle 90]
	\node (UL) at (0,1) {$S \times S$};
	\node (UM) at (2.5,1) {$C \times S$};
	\node (UR) at (5,1) {$S \times S$};
	\node (LL) at (0,-1) {$S$};
	\node (LR) at (5,-1) {$S$};
	
	\node at (6.5,0) {$\text{and}$};
	
	\node (UL2) at (8,1) {$S \times S$};
	\node (UM2) at (10.5,1) {$S \times C$};
	\node (UR2) at (13.5,1) {$S \times S$};
	\node (LL2) at (8,-1) {$S$};
	\node (LR2) at (13.5,-1) {$S$};

	\tikzset{font=\scriptsize};
	\draw[->] (UL) to node [above] {$L \times id$} (UM);
	\draw[->] (UM) to node [above] {$E \times id$} (UR);
	\draw[->] (LL) to node [left] {$\Delta$} (UL);
	\draw[->] (UR) to node [right] {$M$} (LR);
	\draw[->] (LL) to node [above] {$id$} (LR);
	
	\draw[->] (UL2) to node [above] {$id \times R$} (UM2);
	\draw[->] (UM2) to node [above] {$id \times (I \circ E)$} (UR2);
	\draw[->] (LL2) to node [left] {$\Delta$} (UL2);
	\draw[->] (UR2) to node [right] {$M$} (LR2);
	\draw[->] (LL2) to node [above] {$id$} (LR2);
\end{tikzpicture}
\end{center}
commute, and
\item the diagrams
\begin{center}
\begin{tikzpicture}[>=angle 90]
	\node (UL) at (0,1) {$S \times S$};
	\node (UR) at (5,1) {$S \times S$};
	\node (LL) at (0,-1) {$S$};
	\node (LM) at (2.5,-1) {$C$};
	\node (LR) at (5,-1) {$S$};
	
	\node at (6.5,0) {$\text{and}$};
	
	\node (UL2) at (8,1) {$S \times S$};
	\node (UR2) at (13,1) {$S \times S$};
	\node (LL2) at (8,-1) {$S$};
	\node (LM2) at (10.5,-1) {$C$};
	\node (LR2) at (13,-1) {$S$};

	\tikzset{font=\scriptsize};
	\draw[->] (UL) to node [above] {$id \times I$} (UR);
	\draw[->] (LL) to node [left] {$\Delta$} (UL);
	\draw[->] (UR) to node [right] {$M$} (LR);
	\draw[->] (LL) to node [above] {$L$} (LM);
	\draw[->] (LM) to node [above] {$E$} (LR);
	
	\draw[->] (UL2) to node [above] {$I \times id$} (UR2);
	\draw[->] (LL2) to node [left] {$\Delta$} (UL2);
	\draw[->] (UR2) to node [right] {$M$} (LR2);
	\draw[->] (LL2) to node [above] {$R$} (LM2);
	\draw[->] (LM2) to node [above] {$I \circ E$} (LR2);
\end{tikzpicture}
\end{center}
commute,
\end{enumerate}
where $\Delta := \Delta_S$ is the coproduct of the comonoid structure on $S$. We require that all compositions above be strongly transversal.
\end{defn}

\begin{rmk}
The term ``hopfoid'' is meant to suggest that this type of structure is a mix of a ``grouopid'' structure and a ``Hopf algebroid'' structure. The terms ``target'', ``source'', ``unit'', ``product'', ``coproduct'', and ``antipode'' used in the above definition are motivated by this perspective. 
\end{rmk}

In summary then, we have the following:

\begin{thrm}
Let $S$ be a symplectic double groupoid. Then the object $S \rightrightarrows C$ resulting from the constructions of the previous sections is a symplectic hopfoid.
\end{thrm}

\begin{proof}
The remaining conditions are straightforward to check; in particular, (vi) is precisely the compatibility between the two groupoid structures on $S$.
\end{proof}

We note that the symplectic hopfoid obtained by applying the procedure of the previous section to the transposed double groupoid is simply the one obtained by taking the transpose of all structure morphisms of $S \rightrightarrows C$; we will say that the symplectic hopfoids obtained in this manner are \emph{dual} to one another. To be clear, the transpose $E^t: S \to C$ of the unit of $S \rightrightarrows C$ is the target of the dual symplectic hopfoid, the transpose of the target $L: S \to C$ is the unit of the dual hopfoid, and the transposes of the product and coproduct of $S \rightrightarrows C$ are respectively the coproduct and product of the dual hopfoid.

\section{Back to General Double Groupoids}
Here we note that the above constructions can also be carried out even when there is no symplectic structure present---i.e. for a general double Lie groupoid. 
From such a double groupoid $D$ with core $C$, we can produce what we will call a \emph{Lie hopfoid} structure on $D \rightrightarrows C$ in the category of smooth manifolds and smooth relations, defined via the same data and diagrams as that for a symplectic hopfoid.\footnote{In this case we require that the target and source be surmersions, that the unit be a cosurmersion, and that the antipode be a diffeomorphism.}

Indeed, suppose that $D$ is a double Lie groupoid. The description of the leaves of Lemma \ref{lem:fol} still defines a foliation on $\wt r_{H}^{-1}(1^V M)$ whose leaf space may be identified with the core $C$ in the same manner. Similarly, the core can also be realized as the leaf space of certain foliations on
\[
\wt r_V^{-1}(1^H M),\ \wt \ell_H^{-1}(1^V M), \text{ and } \wt \ell_V^{-1}(1^H M),
\]
given by the same expressions as those in the symplectic case. Then all the constructions of the previous sections produce the Lie hopfoid structure $D \rightrightarrows C$ we are interested in. For example, the ``target'' $L_D: D \to C$ obtained by realizing $C$ as the leaf space of the correct foliation on $\wt r_V^{-1}(1^H M)$ is the smooth relation
\[
L_D: s \mapsto s\,\wt\circ_V\,\wt 1^H_{i_V(\wt r_H(s))} \text{ for } s \in \wt r_V^{-1}(1^H M),
\]
which is the analog of equation \ref{L} in the symplectic case.

\begin{rmk}
We should note that although the constructions of this chapter now work more generally in the double Lie groupoid case, the expressions for the various foliations used only came about by first considering the symplectic case where these foliations are the characteristic foliations of certain coisotropic submanifolds. Indeed, it is unclear that one would have guessed to consider these foliations in the first place without the symplectic case as a guide.
\end{rmk}

These two ``hopfoid'' structures are well-behaved with respect to the cotangent functor in the following manner:

\begin{thrm}\label{thrm:lie-hopf}
Suppose that $D$ is a double Lie groupoid with core $C$ and endow $T^*D$ with the induced symplectic double groupoid structure of Theorem \ref{thrm:ctdbl} with core $T^*C$. Then the induced symplectic hopfoid structure on $T^*D \rightrightarrows T^*C$ is simply the one obtained by applying the cotangent functor to the Lie hopfoid $D \rightrightarrows C$.
\end{thrm}

In particular, consider the double groupoid of Example \ref{ex:dmain}. The Lie hopfoid corresponding to this turns out to be the groupoid $G \rightrightarrows M$ itself, where the structure maps are now viewed as smooth relations. The above result  generalizes the fact that the structure $T^*G \rightrightarrows T^*M$ obtained by applying the cotangent functor to this hopfoid is the symplectic hopfoid obtained from the induced symplectic double groupoid structure on $T^*G$.

\begin{proof}
First we show that the cotangent lift $T^*L_D$ of the target relation $L_D$ of $D \rightrightarrows C$ is the target relation of the symplectic hopfoid $T^*D \rightrightarrows T^*C$ constructed from the symplectic double groupoid $T^*D$. Set $X := \wt r_V^{-1}(1^H M)$, and let $\pi: X \to C$ be the surjective submersion
\[
s \mapsto s\,\wt\circ_V\,\wt 1^H_{i_V(\wt r_H(s))}
\]
and $i: X \to D$ the inclusion. Viewing both of these as relations, the target relation $L_D: D \to C$ is then the composition
\begin{center}
\begin{tikzpicture}[>=angle 90]
	\node (1) at (0,1) {$D$};
	\node (2) at (2,1) {$X$};
	\node (3) at (4,1) {$C.$};

	\tikzset{font=\scriptsize};
	\draw[->] (1) to node [above] {$i^t$} (2);
	\draw[->] (2) to node [above] {$\pi$} (3);
\end{tikzpicture}
\end{center}
Hence the cotangent lift $T^*L_D: T^*D \to T^*C$ is the composition (of reductions)
\begin{center}
\begin{tikzpicture}[>=angle 90]
	\node (1) at (0,1) {$T^*D$};
	\node (2) at (2,1) {$T^*X$};
	\node (3) at (4,1) {$T^*C,$};

	\tikzset{font=\scriptsize};
	\draw[->>] (1) to node [above] {$red$} (2);
	\draw[->>] (2) to node [above] {$red_\pi$} (3);
\end{tikzpicture}
\end{center}
where the first reduction is the one obtained by reducing the coisotropic submanifold $T^*D|_X$ of $T^*D$ and the second is the one obtained by reducing the coisotropic $N^*\F_\pi$ of $T^*X$, where $\F_\pi$ is the foliation on $X$ given by the fibers of $\pi$. Note that the domain of $T^*L_D$ is simply $N^*\F_\pi$.

To show that this composition is the target $L: T^*D \red T^*C$ of the symplectic hopfoid $T^*D \rightrightarrows T^*C$, we compute the domain of the latter. Recall that the domain of this canonical relation consists of those $(s,\xi) \in T^*D$ which map under the source $\alpha: T^*D \to A^*H$ of the left groupoid structure on $T^*D$ to a unit of the bottom groupoid structure. Identifying $A^*H$ with $N^*H$ in $T^*D$, a calculation using the explicit descriptions of $\alpha$ and of the unit $A^*C \to A^*H$ of the bottom groupoid shows that $\dom L$ then consists of those $(s,\xi) \in T^*D$ such that $s \in X$ and
\[
(dL_s)_{\wt 1^V 1^H_m}^*\left(\xi|_{\ker(d\wt\ell_V)_s}\right) \text{ restricted to $V$ is zero},
\]
where $\wt r_V(s) = 1^H_m$, $L_s$ is left multiplication by $s$ in the left groupoid structure on $T^*D$, and by ``restricted to $V$'' we mean restricted to the elements of $\ker(d\wt \ell_V)_{\wt 1^V1^H_m}$ coming from $TV \subset TD$ under the decomposition
\[
TD|_H = TH \oplus \ker(d\wt\ell_V)|_H.
\]
Explicitly, such elements are of the form
\[
(d\wt 1^H)_{1^V_m}u - (d\wt 1^V)_{1^H_m}(d\wt\ell_V)_{\wt 1^H 1^V_m}(d\wt 1^H)_{1^V_m}u \text{ for } u \in TV,
\]
which, using $\wt\ell_V \circ \wt 1^H = 1^H \circ \ell _V$, we can write as
\[
(d\wt 1^H)_{1^V_m}\left[u - (d 1^V)_m(d \ell_V)_{1^V_m}u\right].
\]
The expression in square brackets above gives all of $\ker(d\ell_V)_{1^Vm}$ when varying $u$, and so finally we find that the condition on $\xi$ is that
\[
\xi \text{ vanishes on } (dL_s)_{\wt 1^V 1^H_m}(d\wt 1^H)_{1^V_m}(\ker(d\ell_V)_{1^V_m}).
\]
Hence $(s,\xi) \in T^*D$ is in $\dom L$ if and only if $s \in X$ and $\xi$ satisfies the above condition.

According to the analog of equation \ref{tang-dist} in the general double groupoid case giving a description of the tangent distribution to $\F_\pi$, we see that $(s,\xi) \in T^*D$ is in the domain of $L$ if and only if $(s,\xi)$ is in $N^*\F_\pi$. Thus the canonical relations $T^*L_D$ and $L$ have the same domain, and since both are reductions we conclude that $T^*L_D = L$.

A similar argument shows that the cotangent lifts of the source $R_D: D \to C$ and unit $E_D: C \to D$ of $D \rightrightarrows C$ are the source and unit of $T^*D \rightrightarrows T^*C$ respectively. From the definition of the symplectic groupoid structure on the cotangent bundle of a groupoid it follows that the product, coproduct, and antipode of the symplectic hopfoid $T^*D \rightrightarrows T^*C$ are indeed the cotangent lifts of the product, coproduct, and antipode of the Lie hopfoid $D \rightrightarrows C$, and the theorem is proved.
\end{proof}

\begin{ex}\label{ex:dinertia-hopf}
Consider the double groupoid of Example~\ref{ex:dinertia}. The resulting Lie hopfoid structure $G \times G \rightrightarrows G$ is the following. First, the target and source relations $G \times G \to G$ are respectively given by
\[
(g,h) \mapsto gh^{-1} \text{ for $g,h \in G$ such that $r(g)=r(h)$}
\]
and
\[
(g,h) \mapsto h^{-1}g \text{ for $g,h \in G$ such that $\ell(g)=\ell(h)$}.
\]
The unit relation $G \to G \times G$ and antipode $G \to G$ are respectively $g \mapsto \{g\} \times M$ and $(g,h) \mapsto (h^{-1},g^{-1})$. Finally, the product $G \times G \times G \times G \to G \times G$ is
\[
(g,h,a,b) \mapsto (gh,ab) \text{ for composable $g,h$ and composable $a,b$}
\]
and the coproduct $G \times G \to G \times G \times G \times G$ is
\[
(g,h) \mapsto (g,k,k,h) \text{ where $k$ is any element of $G$}.
\]

According to the above theorem, the symplectic hopfoid corresponding to the symplectic double groupoid of Example~\ref{ex:inertia} is the cotangent lift of this Lie hopfoid.\footnote{At first glance, these relations seem to come out of nowhere, but we will see in the next chapter that they are actually the structure maps of the so-called \emph{inertia groupoid} of $G \rightrightarrows M$ in disguise.}
\end{ex}

\section{Symplectic Category $\to$ Double Structures}
Suppose that $S \rightrightarrows C$ is a symplectic hopfoid. In this section we show how to recover from this data a symplectic double groupoid structure on $S$.

As a guide to how this procedure works, let us first consider the case where $C = pt$. In this case, since $L$ and $R$ preserve counits it follows that $L$ and $R$ are both the counit $\varepsilon: S \to pt$ of the given comonoid structure on $S$. Call the lagrangian submanifold of $S$ given by this morphism $P$, and denote by $P^*$ the lagrangian submanifold given by the unit $E: pt \to S$. It follows from the calculations below that $(S,M,E)$ is then a strongly positive symplectic $*$-monoid, so that we now have two symplectic groupoid structures on $S$; one with base $P$ and one with base $P^*$. This in turn induces Poisson-Lie group structures on $P$ and $P^*$, so that $S$ becomes a symplectic double group.

In general we proceed as follows. First, the given strongly positive $*$-comonoid structure on $S$ gives rise to a symplectic groupoid $S \rightrightarrows P$, where $P$ is the lagrangian submanifold of $S$ representing the counit $\varepsilon_S: S \to pt$. Keeping with our notation for symplectic double groupoids, we will denote the structure maps of this symplectic groupoid as
\[
\wt\ell_{P^*}, \wt r_{P^*}, \wt 1^{P^*}, \text{etc.}
\]

Now, denote the image of the composition
\begin{center}
\begin{tikzpicture}[>=angle 90]
	\node (UL) at (0,1) {$pt$};
	\node (UM) at (2,1) {$C$};
	\node (UR) at (4,1) {$S$};

	\tikzset{font=\scriptsize};
	\draw[->] (UL) to node [above] {$\varepsilon_C^t$} (UM);
	\draw[->] (UM) to node [above] {$E$} (UR);
\end{tikzpicture}
\end{center}
by $P^*$; this is a lagrangian submanifold of $S$ since the above composition is strongly transversal given that $E$ is a coreduction. As the notation suggests, this will form the base of the second symplectic groupoid structure on $P$.

\begin{prop}
$(S,m,e)$, where $m$ is the product of the symplectic hopfoid structure and $e: pt \to S$ is the morphism given by $P^*$, is a symplectic monoid.
\end{prop}

\begin{proof}
By assumption, $m$ is a associative. The unit properties of $e$ follow from the diagrams in condition (vii) in the definition of a symplectic hopfoid as follows.

First, we note the following: to say that $L$ preserves counits means that $\varepsilon_C \circ L = \varepsilon_S$, which after composing both sides on the right by $L^t$ and taking tranposes of the results implies that
\[
L \circ \varepsilon_S^t = \varepsilon_C^t.
\]
Here we have used the fact that $L \circ L^t = id_C$ since $L$ is a reduction.

Now, consider the first diagram in condition (vii):
\begin{center}
\begin{tikzpicture}[>=angle 90]
	\node (UL) at (0,1) {$S \times S$};
	\node (UM) at (2.5,1) {$C \times S$};
	\node (UR) at (5,1) {$S \times S$};
	\node (LL) at (0,-1) {$S$};
	\node (LR) at (5,-1) {$S.$};

	\tikzset{font=\scriptsize};
	\draw[->] (UL) to node [above] {$L \times id$} (UM);
	\draw[->] (UM) to node [above] {$E \times id$} (UR);
	\draw[->] (LL) to node [left] {$\Delta$} (UL);
	\draw[->] (UR) to node [right] {$M$} (LR);
	\draw[->] (LL) to node [above] {$id$} (LR);
\end{tikzpicture}
\end{center}
For any $s \in S$ the canonical relation $\Delta$ in particular contains an element of the form $(s,\wt 1^{P^*}_\lambda,s)$, which follows from the fact that $\Delta^t$ is in fact a groupoid product. This element is in the canonical relation
\[
\varepsilon_S^t \times id: S \to S \times S,
\]
and thus the composition $(L \times id) \circ \Delta$ contains $(L \times id) \circ (\varepsilon_S^t \times id)$, which from the above noted fact is $\varepsilon_C^t \times id$. Thus the commutativity of the above diagram implies the commutativity of
\begin{center}
\begin{tikzpicture}[>=angle 90]
	\node (UL) at (0,1) {$C \times S$};
	\node (UR) at (5,1) {$S \times S$};
	\node (LL) at (0,-1) {$S$};
	\node (LR) at (5,-1) {$S.$};

	\tikzset{font=\scriptsize};
	\draw[->] (UL) to node [above] {$E \times id$} (UR);
	\draw[->] (LL) to node [left] {$\varepsilon_C^t \times id$} (UL);
	\draw[->] (UR) to node [right] {$M$} (LR);
	\draw[->] (LL) to node [above] {$id$} (LR);
\end{tikzpicture}
\end{center}
Recalling the definition of $e = E \circ \varepsilon_C^t$, this then says that $e$ is a left unit for $m$. The right unit property of $e$ follows similarly from the second diagram in condition (vii), and we thus conclude that $(S,m,e)$ is a symplectic monoid.
\end{proof}

Now, it is simple to check that $I \circ \wt i_{P^*}$ is a $*$-structure on the above monoid using condition (v) in the definition of a symplectic hopfoid and the corresponding $*$-properties of the $*$-structure on $(S,\Delta_S,\varepsilon_S)$.

\begin{prop}
The above $*$-structure is strongly positive.
\end{prop}

\begin{proof}
Consider the first diagram in condition (viii) in the definition of symplectic hopfoid. It is a simple check to see that the following then commutes:
\begin{center}
\begin{tikzpicture}[>=angle 90]
	\node (P) at (-2,-2) {$pt$};
	\node (UL) at (0,1) {$S \times S$};
	\node (UR) at (5,1) {$S \times S$};
	\node (LL) at (0,-1) {$S$};
	\node (LM) at (2.5,-1) {$C$};
	\node (LR) at (5,-1) {$S$};

	\tikzset{font=\scriptsize};
	\draw[->] (P) to node [above] {$\varepsilon_S^t$} (LL);
	\draw[->] (P) to [bend left=15] node [left] {$\wt i_{P^*}$} (UL);
	\draw[->] (P) to [bend right=15] node [above] {$\varepsilon_C^t$} (LM);
	\draw[->] (UL) to node [above] {$id \times I$} (UR);
	\draw[->] (LL) to node [left] {$\Delta$} (UL);
	\draw[->] (UR) to node [right] {$M$} (LR);
	\draw[->] (LL) to node [above] {$L$} (LM);
	\draw[->] (LM) to node [above] {$E$} (LR);
	\draw[->] (P) to [bend right=15] node [below] {$E \circ \varepsilon_C^t$} (LR);
\end{tikzpicture}
\end{center}
where $\wt i_{P^*}: \overline{S} \to S$ is now thought of as a morphism $pt \to S \times S$ via its graph. Now, we can factor the top row as
\begin{center}
\begin{tikzpicture}[>=angle 90]
	\node (UL) at (0,1) {$S \times S$};
	\node (UM) at (3,1) {$S \times \overline{S}$};
	\node (UR) at (7,1) {$S \times S$};

	\tikzset{font=\scriptsize};
	\draw[->] (UL) to node [above] {$id \times \wt i_{P^*}$} (UM);
	\draw[->] (UM) to node [above] {$id \times (I \circ \wt i_{P^*})$} (UR);
\end{tikzpicture}
\end{center}
using the fact that $\wt i_{P^*} \circ \wt i_{P^*} = id$. Then the composition
\begin{center}
\begin{tikzpicture}[>=angle 90]
	\node (UL) at (0,1) {$pt$};
	\node (UM) at (3,1) {$S \times S$};
	\node (UR) at (6,1) {$S \times \overline{S}$};

	\tikzset{font=\scriptsize};
	\draw[->] (UL) to node [above] {$\wt i_{P^*}$} (UM);
	\draw[->] (UM) to node [above] {$id \times \wt i_{P^*}$} (UR);
\end{tikzpicture}
\end{center}
is precisely the morphism $pt \to S \times \overline{S}$ given by the diagonal of $S \times \overline{S}$. Thus the above diagram becomes
\begin{center}
\begin{tikzpicture}[>=angle 90]
	\node (UL) at (0,1) {$S \times \overline{S}$};
	\node (UR) at (4,1) {$S \times S$};
	\node (LL) at (0,-1) {$pt$};
	\node (LR) at (4,-1) {$S$};

	\tikzset{font=\scriptsize};
	\draw[->] (UL) to node [above] {$id \times (I \circ \wt i_{P^*})$} (UR);
	\draw[->] (LL) to node [left] {$$} (UL);
	\draw[->] (UR) to node [right] {$m$} (LR);
	\draw[->] (LL) to node [above] {$E \circ \varepsilon_C^t$} (LR);
\end{tikzpicture}
\end{center}
which is precisely the strong positivity requirement on $(S,m,e)$.
\end{proof}

We thus conclude that the monoid $(S,m,e)$ produces an additional symplectic groupoid structure on $S$. Restricting the groupoid structure on $S \rightrightarrows P$ to $P^*$ produces a groupoid structure on $P^* \rightrightarrows M$, and restricting the groupoid structure on $S \rightrightarrows P^*$ to $P$ gives a groupoid structure $P \rightrightarrows M$.

These give the four groupoid structures in the definition of a symplectic double groupoid, and it is straightforward to check that the various compatibilities hold; in particular, condition (vi) defining a symplectic hopfoid precisely says that
\[
m: S \times S \to S
\]
is a groupoid morphism for the groupoid product $\Delta^t$ on $S$ and vice-versa. We thus have:

\begin{thrm}\label{thrm:main-cor}
There is a $1$-$1$ correspondence (up to symplectomorphism) between symplectic double groupoids and pairs of dual symplectic hopfoids.
\end{thrm}

\begin{proof}
It remains only to check that the operations of producing a symplectic hopfoid from a symplectic double groupoid and conversely producing a symplectic double groupoid from a symplectic hopfoid are inverse to each other. This is a straightforward verification. The need to use pairs of dual hopfoids in this correspondence comes from the fact that we make no distinction between a double groupoid and its tranpose.
\end{proof}

The above theorem then has an obvious generalization to arbitrary double Lie groupoids and pairs of dual Lie hopfoids.
\chapter{Symplectic Orbifolds}\label{chap:stacks}

In this chapter we show how the structures obtained in the previous section may be used to produce symplectic-like structures on stacks; under certain assumptions, what we produce are precisely symplectic orbifolds. This procedure also works, in nice situations, more generally for simplicial symplectic manifolds, suggesting that such structures may be useful in the study of general symplectic orbifolds.

\section{Induced Groupoids}
The structures produced in the previous section are not honest Lie groupoids, in particular because the source and target relations are only partially defined. However, we may restrict these relations to their common domains in an attempt to produce an actual groupoid structure.

As a guide, let us first see what happens with the structure $T^*G \rightrightarrows T^*M$ arising from the standard symplectic double groupoid structure on $T^*G$ for a Lie groupoid $G \rightrightarrows M$, which as we have seen is also the structure obtained after application of the cotangent functor. We will assume that $G \rightrightarrows M$ is a regular groupoid, so that the orbits are of constant dimension. Let $\F$ denote the (regular) foliation of $M$ by the orbits of $G$, and let $\F_\ell$ and $\F_r$ respectively denote foliations of $G$ given by the $\ell$ and $r$-fibers. The domains of the relations $T^*\ell$ and $T^*r$ are respectively the conormal bundles $N^*\F_\ell$ and $N^*\F_r$.

Let $g \in G$, $\xi \in T^*_{\ell(g)}M$, and $\eta \in T^*_{r(g)}M$. Then we have $d\ell_g^*\xi \in N^*\F_\ell$ and $dr_g^*\eta \in N^*\F_r$ by definition. The following observation can be proved via a simple argument using bisections, but we will soon derive it from a more general statement:

\begin{lem}\label{lem:cot-ind}
In the above setup, $d\ell_g^*\xi \in N^*\F_r$ if and only if $\xi \in N^*\F$, and similarly $dr_g^*\eta \in N^*\F_\ell$ if and only if $\eta \in N^*\F$.
\end{lem}

The relations $T^*\ell$ and $T^*r$ are single-valued since $\ell$ and $r$ are submersions, and their common domain is $N^*\F_\ell \cap N^*\F_r$. The above lemma implies that an element of $T^*G$ is in this common domain if and only if the element of $T^*M$ it maps to under $T^*\ell$ and $T^*r$ is an element of the conormal bundle $N^*\F$, which is smooth since $\F$ is a regular foliation. Hence, we obtain from the symplectic hopfoid structure on $T^*G \rightrightarrows T^*M$ an honest groupoid structure on
\[
N^*\F_\ell \cap N^*\F_r \rightrightarrows N^*\F.
\]
Explicitly, the target and source of an element $(g,d\ell_g^*\xi=dr_g^*\eta)$ are respectively $(\ell(g),\xi)$ and $(r(g),\eta)$, and the product is given by
\[
(g,d\ell_g^*\xi=dr_g^*\eta)\cdot(h,d\ell_h^*\eta=dr_h^*\gamma) = (gh,d\ell_{gh}^*\xi = dr_{gh}^*\gamma).
\]

We see that an element $(g,d\ell_g^*\xi=dr_g^*\eta)$ of $N^*\F_\ell \cap N^*\F_r$ is completely determined by $g$ and $\eta \in T_{r(g)}^*M$, and thus $N^*\F_\ell \cap N^*\F_r$ is isomorphic to the fiber product
\[
G \times_M N^*\F
\]
of $r: G \to M$ and the canonical projection $N^*\F \to M$. Under this identification, the groupoid structure on $N^*\F_\ell \cap N^*\F_r \rightrightarrows N^*\F$ becomes the action groupoid
\[
G \ltimes_M N^*\F \rightrightarrows  N^*\F
\]
for the action of $G$ on $N^*\F$ given by $g \cdot (r(g), \eta) = (\ell(g),\xi)$ where $\xi \in T_{\ell(g)}^*M$ is the unique element such that $d\ell_g^*\xi = dr_g^*\eta$, which exists by the previous lemma.

Thus, restricting the structure relations of $T^*G \rightrightarrows T^*M$ to those elements of $T^*G$ on which the source and target relations are both defined produces (in the regular case) an honest Lie groupoid. This generalizes to the symplectic hopfoid structure $S \rightrightarrows C$ resulting from any symplectic double groupoid $S$ as follows:\footnote{We use the same notation for the structure maps of a symplectic double groupoid as in the previous chapter.}

\begin{thrm}\label{thrm:ind}
Let $S$ be a symplectic double groupoid with core $C$, and suppose that the restriction $\wt \ell_P|_C: C \to P^*$ is transverse to $1^{P^*} M \subseteq P^*$. Then the symplectic hopfoid structure on $S \rightrightarrows C$ gives rise to a Lie groupoid structure on
\[
D := \wt r_P^{-1}(1^{P^*}M) \cap \wt \ell_P^{-1}(1^{P^*}M) \rightrightarrows Y
\]
where $Y$ is the preimage of $1^{P^*}M$ under the restriction $\wt\ell_P|_C$.
\end{thrm}

The transversality condition is needed to guarantee that $Y$ is indeed a submanifold of $C$, and we will assume that this is the case throughout the rest of this chapter. In particular, in the case of the symplectic hopfoid $T^*G \rightrightarrows T^*M$, one can check that this condition is equivalent to $G \rightrightarrows M$ being regular.

We prove the above theorem in stages, and will call $D \rightrightarrows Y$ the \emph{induced groupoid} of $S \rightrightarrows C$. Note that elements $s \in D$ look like
\begin{center}
\begin{tikzpicture}[scale=.75]
	\node at (1,1) {$s$};

	\tikzset{font=\scriptsize};
	\draw (0,0) rectangle +(2,2);
	\node at (-.5,1) {$1^{P^*}_m$};
	\node at (2.5,1) {$1^{P^*}_q$};
\end{tikzpicture}
\end{center}
and elements of $Y$ consist of such squares with the top side also being a unit. From this it is easy to show the following, which explains why the base of the induced groupoid is the given submanifold $Y$ of the core:

\begin{lem}
The image of an element $s \in D$ under (the restrictions of) the target and source relations of $S \rightrightarrows C$ lies in $Y$, and such elements of $S$ are the only ones with this property.
\end{lem}

Lemma~\ref{lem:cot-ind} is then a consequence of this, where $Y$ in that case is precisely $N^*\F$. As opposed to using bisections, we now obtain this result simply by considering the form (in square notation) which elements of $D$ and $Y$ take.

Next we have the following technical lemma, which will be important in later computations:

\begin{lem}
For any $s \in D$,
\[
\wt 1^{P^*}_{\wt r_{P^*}(s)}\,\wt\circ_P\,\,\wt i_P\wt i_{P^*}(s)\,\,\wt\circ_P\,\wt 1^{P^*}_{\wt\ell_{P^*}(s)} = s
\]
\end{lem}

\begin{proof}
The composition of the first two terms can be written as
\begin{center}
\begin{tikzpicture}[scale=.9]
	\draw (0,1) rectangle (2,3);
	\node at (1,2) {$\wt 1^{P^*}_{\wt r_{P^*}(s)}$};
	\draw (2,1) rectangle (4,3);
	\node at (3,2) {$\wt i_P \wt i_{P^*}(s)$};
	
	\node at (4.5,2) {$=$};

	\draw (5,0) rectangle (7,2);
	\node at (6,1) {$\wt i_{P^*}(s)$};
	\draw (7,0) rectangle (9,2);
	\node at (8,1) {$\wt i_P \wt i_{P^*}(s)$};
	\draw (5,2) rectangle (7,4);
	\node at (6,3) {$s$};
	\draw (7,2) rectangle (9,4);
	\node at (8,3) {$\wt 1^{P^*}_{i_P(\wt\ell_{P^*}(s))}$};
	
	\node at (9.5,2) {$=$};
	
	\draw (10,0) rectangle (14,2);
	\node at (12,1) {$\wt 1^P 1^{P^*}_m$};
	\draw (10,2) rectangle (12,4);
	\node at (11,3) {$s$};
	\draw (12,2) rectangle (14,4);
	\node at (13,3) {$\wt 1^{P^*}_{i_P(\wt\ell_{P^*}(s))}$};
\end{tikzpicture}
\end{center}
where in the first step we decompose vertically and in the second step we compose horizontally. Similarly, we can decompose
\begin{center}
\begin{tikzpicture}[scale=.75]
	\draw (0,1) rectangle (2,3);
	\node at (1,2) {$\wt 1^{P^*}_{\wt\ell_{P^*}(s)}$};
	
	\node at (2.5,2) {$=$};

	\draw (3,0) rectangle (5,2);
	\node at (4,1) {$\wt 1^{P^*}_{\wt\ell_{P^*}(s)}$};
	\draw (3,2) rectangle (5,4);
	\node at (4,3) {$\wt 1^{P^*}_{\wt\ell_{P^*}(s)}$};
\end{tikzpicture}
\end{center}
Thus we have that $\wt 1^{P^*}_{\wt r_{P^*}(s)}\,\wt\circ_P\,\,\wt i_P\wt i_{P^*}(s)\,\,\wt\circ_P\,\wt 1^{P^*}_{\wt\ell_{P^*}(s)}$ is given by
\begin{center}
\begin{tikzpicture}[scale=.9]
	\draw (0,0) rectangle (4,2);
	\node at (2,1) {$\wt 1^P 1^{P^*}_m$};
	\draw (0,2) rectangle (2,4);
	\node at (1,3) {$s$};
	\draw (2,2) rectangle (4,4);
	\node at (3,3) {$\wt 1^{P^*}_{i_P(\wt\ell_{P^*}(s))}$};
	\draw (4,0) rectangle (6,2);
	\node at (5,1) {$\wt 1^{P^*}_{\wt\ell_{P^*}(s)}$};
	\draw (4,2) rectangle (6,4);
	\node at (5,3) {$\wt 1^{P^*}_{\wt\ell_{P^*}(s)}$};
	
	\node at (6.5,2) {$=$};
	
	\draw (7,0) rectangle (11,2);
	\node at (9,1) {$\wt 1^{P^*}_{\wt\ell_{P^*}(s)}$};
	\draw (7,2) rectangle (9,4);
	\node at (8,3) {$s$};
	\draw (9,2) rectangle (11,4);
	\node at (10,3) {$\wt 1^{P^*} 1^P_q$};
	
	\node at (11.5,2) {$=$};
	
	\draw (12,1) rectangle (14,3);
	\node at (13,2) {$s$};
\end{tikzpicture}
\end{center}
as was to be shown, where in the first step we compose horizontally, and in the second step we compose horizontally (using the fact that $\wt 1^{P^*} 1^P_q = \wt 1^P 1^{P^*}_q$) and then vertically.
\end{proof}

We now describe the groupoid structure on $D \rightrightarrows Y$. First, the target $L(s)$ and source $R(s)$ of an element $s \in D$ are respectively
\[
s\, \wt \circ_P\, \wt 1^{P^*}_{\wt r_{P^*}(\wt i_P(s))}\ \text{ and }\ \wt i_P \wt i_{P^*} (s)\, \wt \circ_P\, \wt 1^{P^*}_{\wt \ell_{P^*}(s)}.
\]
These are just the images of $s$ under the (single-valued) reductions $L, R: S \to C$. Thus two elements $s, s' \in D$ will be composable when
\[
\wt i_P \wt i_{P^*} (s)\, \wt \circ_P\, \wt 1^{P^*}_{\wt \ell_{P^*}(s)} = s'\, \wt \circ_P\, \wt 1^{P^*}_{\wt r_{P^*}(\wt i_P(s'))}.
\]

\begin{lem}
For composable $s$ and $s'$, we have
\[
\wt 1^{P^*}_{\wt r_{P^*}(s)}\,\wt\circ_P\,s' = s\,\wt\circ_P\,\wt 1^{P^*}_{\wt r_{P^*}(s')}.
\]
\end{lem}

\begin{proof}
By the previous lemma we know that
\[
s = \wt 1^{P^*}_{\wt r_{P^*}(s)}\,\wt\circ_P\,\,\wt i_P\wt i_{P^*}(s)\,\,\wt\circ_P\,\wt 1^{P^*}_{\wt\ell_{P^*}(s)},
\]
which, using our assumption that $s$ and $s'$ are composable, we can write as
\[
s = \wt 1^{P^*}_{\wt r_{P^*}(s)} \wt \circ_P s'\, \wt \circ_P\, \wt 1^{P^*}_{\wt r_{P^*}(\wt i_P(s'))}.
\]
Composing both sides on the right by $\wt i_P\left(\wt 1^{P^*}_{\wt r_{P^*}(\wt i_P(s'))}\right) = \wt 1^{P^*}_{\wt r_{P^*}(s')}$ under $\wt\circ_P$ then gives the desired equality.
\end{proof}

We define the product of composable $s$ and $s'$ to be the expression given by the lemma:
\[
s \circ s' := \wt 1^{P^*}_{\wt r_{P^*}(s)}\,\wt\circ_P\,s' = s\,\wt\circ_P\,\wt 1^{P^*}_{\wt r_{P^*}(s')},
\]
and claim that this indeed defines a groupoid structure on $D \rightrightarrows Y$. As a first step, we have the following verification:

\begin{prop}
The above product is associative, and for composable $s$ and $s'$ we have $L(s \circ s') = L(s)$ and $R(s \circ s') = R(s')$.
\end{prop}

\begin{proof}
Associativity follows from the fact that $\wt 1^{P^*}_{\wt r_{P^*}(s') \circ_P \wt r_{P^*}(s'')} = \wt 1^{P^*}_{\wt r_{P^*}(s')} \wt\circ_P \wt 1^{P^*}_{\wt r_{P^*}(s'')}$ and can be easily checked. We have
\begin{align*}
L(s \circ s') &= (s \circ s') \wt\circ_P \wt 1^{P^*}_{i_P(\wt r_{P^*}(s \circ s'))} \\
&= (s \circ s') \wt\circ_P \wt 1^{P^*}_{i_P(\wt r_{P^*}(s) \circ \wt r_{P^*}(s'))} \\
&= \left(s \wt\circ_P \wt 1^{P^*}_{\wt r_{P^*}(s')}\right) \wt\circ_P \wt 1^{P^*}_{i_P(\wt r_{P^*}(s')) \circ_P i_P(\wt r_{P^*}(s))} \\
&= s \wt\circ_P \left( \wt 1^{P^*}_{\wt r_{P^*}(s')} \wt\circ_P \wt 1^{P^*}_{i_P(\wt r_{P^*}(s'))}\right) \wt\circ_P \wt 1^{P^*}_{i_P(\wt r_{P^*}(s))} \\
&= s \wt\circ_P \left(\wt 1^{P^*}1^P_m\right) \wt\circ_P \wt 1^{P^*}_{i_P(\wt r_{P^*}(s))} \\
&= s \wt\circ_P \wt 1^{P^*}_{i_P(\wt r_{P^*}(s))} \\
&= L(s)
\end{align*}
where we have used the fact that $\wt 1^{P^*} \circ 1^P = \wt 1^P \circ 1^{P^*}$. Similarly, we have
\begin{align*}
R(s \circ s') &= \wt i_P \wt i_{P^*}(s \circ s') \wt\circ_P \wt 1^{P^*}_{\wt\ell_{P^*}(s \circ s')} \\
&= \wt i_P \wt i_{P^*}(s \circ s') \wt\circ_P \wt 1^{P^*}_{\wt r_{P^*}(s) \circ_P \wt \ell_{P^*}(s')} \\
&= \wt i_P \wt i_{P^*}\left(\wt 1^{P^*}_{\wt r_{P^*}(s)} \wt\circ_P s'\right) \wt\circ_P \wt 1^{P^*}_{\wt r_{P^*}(s) \circ_P \wt \ell_{P^*}(s')} \\
&= \wt i_P \wt i_{P^*}(s') \wt\circ_P \left(\wt 1^{P^*}_{i_P(\wt r_{P^*}(s))} \wt\circ_P \wt 1^{P^*}_{\wt r_{P^*}(s)}\right) \wt\circ_P \wt 1^{P^*}_{\wt\ell_{P^*}(s')} \\
&= \wt i_P \wt i_{P^*}(s') \wt\circ_P \left(\wt 1^{P^*}1^P_q\right) \wt\circ_P 1^{P^*}_{\wt\ell_{P^*}(s')} \\
&= \wt i_P \wt i_{P^*}(s') \wt\circ_P 1^{P^*}_{\wt\ell_{P^*}(s')} \\
&= R(s')
\end{align*}
as claimed.
\end{proof}

\begin{rmk}
In the structure $T^*G \rightrightarrows T^*M$ arising from applying the cotangent functor to a groupoid, the above product becomes exactly the aforementioned product on the groupoid
\[
G \ltimes_M N^*\F \rightrightarrows N^*\F.
\]
If instead we had attempted to define the product of $s$ and $s'$ simply via $s\,\wt\circ_P s'$, sources and targets would not be preserved as in the above proposition; this is already apparent in this example, and explains the need for defining the groupoid product as we have.
\end{rmk}

The inverse of $D \rightrightarrows Y$ is induced by the antipode of the symplectic hopfoid $S \rightrightarrows C$ and is given by the following:
\[
s^{-1} := \wt 1^{P^*}_{i_P(\wt r_{P^*}(s))\, \circ_P\, \wt\ell_{P^*}(s)}\, \wt\circ_P\, \wt i_P \wt i_{{P^*}(s)} = \wt 1^{P^*}_{i_P(\wt r_{P^*}(s))}\,\wt\circ_P\,s\,\wt\circ_P\,\wt 1^{P^*}_{i_P(\wt r_{P^*}(s))}.
\]
Note that if $\wt r_{P^*}(s) = \wt \ell_{P^*}(s)$, the above inverse is simply $\wt i_P \wt i_{P^*}(s)$; this is the case in Example \ref{ex:main} for instance. The unit is induced by the canonical relation $E: C \to S$, and is given by the inclusion of the core $C$ in $S$.

\begin{proof}[Proof of Theorem \ref{thrm:ind}]
The check that units and inverses behave as expected is a straightforward verification. Moreover, all the defined structure maps of $D \rightrightarrows Y$ are smooth since they can be expressed as compositions of the smooth structure maps of $S$; in particular, the fact that $L, R: S \to C$ are reductions implies that the induced maps they induce on their domains are surjective submersions.
\end{proof}

\begin{rmk}
As mentioned in the previous chapter, the same construction which produces a symplectic hopfoid from a symplectic double groupoid generalizes to produce a Lie hopfoid from a double Lie groupoid. In this latter case, the same construction as above also produces an honest Lie groupoid.
\end{rmk}

\section{Examples of Induced Groupoids}
To get a feel for the induced groupoid construction, we consider the following examples.

\begin{ex}
As mentioned, going through all this with the structure from Example \autoref{ex:main} produces
\[
G \ltimes_M N^*\F \rightrightarrows N^*\F.
\]
In fact, the action of $G$ on $N^*\F$ used here is well known and can be described as follows. A groupoid $G \rightrightarrows M$ in general acts neither on $T^*M$ nor $A^*$, but only acts ``up to homotopy''\footnote{See \cite{AC} and \cite{ELW}.} on the complex
\begin{center}
\begin{tikzpicture}[>=angle 90]
	\node (0) at (-2,1) {$0$};
	\node (1) at (0,1) {$T^*M$};
	\node (2) at (2,1) {$A^*$};
	\node (3) at (3.5,1) {$0$};

	\tikzset{font=\scriptsize};
	\draw[->] (0) to node [above] {$$} (1);
	\draw[->] (1) to node [above] {$\rho^*$} (2);
	\draw[->] (2) to node [above] {$$} (3);
\end{tikzpicture}
\end{center}
obtained by dualizing the anchor $\rho: A \to TM$ of the Lie algebroid of $G$. This action ``up to homotopy'' induces an honest action of $G$ on the cokernel $N^*\F$, which is the action appearing above.

Performing this construction on the transposed double groupoid structure on $T^*G$ produces the induced Lie groupoid
\[
A^* \oplus T^*M \rightrightarrows T^*M,
\]
which can be viewed as the action groupoid for the trivial action of $A^* \rightrightarrows M$ on $T^*M$.
\end{ex}

\begin{ex}\label{ex:inertia-induced}(Inertia Groupoid)
Consider the symplectic double groupoid $T^*G \times T^*G$ of Example \autoref{ex:inertia}, and the corresponding symplectic hopfoid $T^*G \times T^*G \rightrightarrows T^*G$. The procedure above produces the following.

We first focus on what is happening along zero sections---i.e. we compute the induced groupoid of the Lie hopfoid $G \times G \rightrightarrows G$ of Example~\ref{ex:dinertia-hopf}. We get a groupoid structure on
\[
G \times_{M \times M} G \rightrightarrows LG
\]
where $G \times_{M \times M} G = \{(g,h)\ |\ \ell(g)=\ell(h) \text{ and } r(g) = r(h)\}$ and $LG$ is the ``isotropy subgroupoid''  $\{g\ |\ \ell(g) = r(g)\}$. The source and target are
\[
L(g,h) = gh^{-1} \text{ and } R(g,h) = h^{-1}g,
\]
the unit is $E(k) = (k, 1)$, the inverse is $I(g,h) = (h^{-1}gh^{-1},h^{-1})$, and the product is defined for $(g,h),(a,b)$ such that $gb = ha$ and is
\[
(g,h) \cdot (a,b) = (ha,hb) = (gb,hb).
\]
In fact, this groupoid structure is well known: we have a diffeomorphism
\[
G \times_{M \times M} G \to G \times_M LG
\]
given by
\[
(g,h) \mapsto (h,h^{-1}g),
\]
and under this map, the groupoid structure above becomes the so-called \emph{inertia groupoid}\footnote{We refer to \cite{IM} for the definition and properties of the inertia groupoid of a groupoid.} structure on $G \ltimes_M LG \rightrightarrows LG$. 

According to the relation between Lie and symplectic hopfoids given in Theorem~\ref{thrm:lie-hopf} of the last section then, the induced groupoid of the symplectic hopfoid $T^*G \times T^*G \rightrightarrows T^*G$ is the induced groupoid of the cotangent lift of this inertia groupoid. The induced groupoid resulting from the transpose of the given double groupoid structure on $G \times G$ is the pair groupoid $M \times M \rightrightarrows M$.
\end{ex}

\begin{ex}
Let $S$ be a symplectic double group. The induced groupoid of the corresponding symplectic hopfoid $S \rightrightarrows pt$ is simply the Poisson-Lie group $P \rightrightarrows pt$, and the induced groupoid of the dual symplectic hopfoid is the Poisson-Lie group $P^* \rightrightarrows pt$.
\end{ex}

\begin{ex}
This example is inspired by the construction in \cite{ZZ}. Consider a crossed module $(G,H,t,\phi)$ of Lie groups, so $t: H \to G$ is a homomorphism and $\phi: G \to Aut(H)$ is a left action of $G$ on $H$ such that some compatibilities hold. Then we have the double groupoid
\begin{center}
\begin{tikzpicture}[>=angle 90]
	\node (LL) at (0,0) {$pt$};
	\node (LR) at (2.5,0) {$pt$};
	\node (UL) at (0,2) {$H \times G$};
	\node (UR) at (2.5,2) {$G$};

	\tikzset{font=\scriptsize};
	\draw[->] (.4,.07) -- (2.1,.07);
	\draw[->] (.4,-.07) -- (2.1,-.07);
	
	\draw[->] (-.07,1.6) -- (-.07,.4);
	\draw[->] (.07,1.6) -- (.07,.4);
	
	\draw[->] (.8,2.07) -- (2.1,2.07);
	\draw[->] (.8,2-.07) -- (2.1,2-.07);
	
	\draw[->] (2.5-.07,1.6) -- (2.5-.07,.4);
	\draw[->] (2.5+.07,1.6) -- (2.5+.07,.4);
\end{tikzpicture}
\end{center}
where the left side is the semi-direct product group structure on $H \times G$ coming from the action $\phi$ and the top is the action groupoid for the action of $H$ on $G$ induced by $t$. We then get two Lie hopfoid structures on $H \times G \rightrightarrows H$ by applying the construction of the previous chapter to this and its tranpose. The induced groupoids of these two hopfoids are respectively the trivial groupoid
\[
\ker t \rightrightarrows \ker t
\]
and the action groupoid
\[
H \times G \rightrightarrows H
\]
for the right action of $G$ on $H$ obtained by precomposing $\phi$ with the inversion map of $G$.
\end{ex}

It is unclear what relationship there is between the induced groupoids of a (symplectic) Lie hopfoid and its dual, other than the fact that they both arise from the same double groupoid. In the symplectic double group case, the resulting induced groupoids were dual Poisson-Lie groups---it would be interesting to know if there is a more general form of duality between Lie groupoids which accounts for these observations.

\section{Symplectic Structures}
The induced groupoids of symplectic hopfoids come equipped with additional structure, reflecting the fact that they arose from canonical relations. We first consider the induced groupoid $G \ltimes_M N^*\F \rightrightarrows N^*\F$. Here, $N^*\F$ is a coisotropic submanifold of $T^*M$, and so comes equipped with a closed $2$-form. We then have:

\begin{prop}\label{prop:cot}
The induced $2$-form on $N^*\F$ has the property that $L^*\omega = R^*\omega$, where $L$ and $R$ are the target and source of $G \ltimes_M N^*\F \rightrightarrows N^*\F$. Moreover, the orbits of this groupoid are the leaves of the characteristic foliation of $N^*\F$.
\end{prop}

The first condition on $\omega$ above should be viewed as a kind of ``invariance'' and the second as a kind of ``non-degeneracy''. Indeed, when the groupoid $G \rightrightarrows M$ presents an orbifold\footnote{We recall this notion in Definition~\ref{orbifold}.} $\X$, so is proper and has finite isotropy groups, $G \ltimes_M N^*\F \rightrightarrows N^*\F$ presents the cotangent orbifold $T^*\X$ of Proposition \ref{prop:cot-orb} and the above $2$-form is precisely the standard symplectic form on $T^*\X$:

\begin{defn}
Let $\X$ be an orbifold presented by a proper groupoid $G \rightrightarrows M$ with finite isotropy groups. A differential $2$-form $\omega \in \Omega^2(M)$ on $\X$ is \emph{non-degenerate} if $\ker\omega = T\F$, where $\F$ is the foliation of $M$ by orbits of $G$; in other words, $\omega$ is non-degenerate if the leaves of the characteristic foliation on $M$ induced by $\ker\omega$ are the orbits of the groupoid. A \emph{symplectic form} on $\X$ is a closed, non-degenerate $2$-form on $\X$, and an orbifold equipped with such a form is a \emph{symplectic orbifold.}
\end{defn}

The following result justifies this definition by implying that it is Morita-invariant:

\begin{prop}[\cite{LM}]
If $\omega$ is a non-degenerate form, in the sense above, on a groupoid presenting an orbifold, then the form corresponding to it by Proposition~\ref{prop:diff-form} on a Morita equivalent groupoid is also non-degenerate.
\end{prop}

We now show that Proposition \ref{prop:cot} generalizes to the induced groupoids of symplectic hopfoids. As a consequence, when such an induced groupoid presents an orbifold, it will in fact naturally be a symplectic orbifold.

We first note with the following general facts:

\begin{lem}
Suppose that $\Lambda: X \to Y$ is a canonical relation between symplectic manifolds $X$ and $Y$ whose domain is a submanifold of $X$. Then for any $p \in \dom \Lambda$, the image $\Lambda(p)$ is a leaf of the characteristic foliation of the coisotropic submanifold $\im \Lambda \subseteq Y$.
\end{lem}

\begin{proof}
This follows from the factorization of $\Lambda$ given in the commutative diagram \eqref{factor} by following $p \in \dom\Lambda$ around the diagram.
\end{proof}

\begin{prop}\label{prop:orbits}
Suppose that $L, R: S \to C$ are two reductions from a symplectic manifold $S$ to a symplectic manifold $C$, and that $I: S \to S$ is a symplectomorphism such that $I^2 = id$ and $R = L \circ I$. Then the domain and image of the composition
\[
L \circ R^t: C \to C
\]
are both the same coisotropic submanifold $Y$ of $C$, and the leaf $Y^\perp_p$ of the characteristic foliation of $Y$ containing $p$ is the image of $p$ under the above composition:
\[
Y^\perp_p = (L \circ R^t)(p).
\]
\end{prop}

\begin{proof}
We claim that $L \circ R^t = R \circ L^t$, which together with the previous lemma then gives the result. This is an easy check:
\[
L \circ R^t = L \circ (L \circ I)^t = L \circ (I^t \circ L^t) = (L \circ I) \circ L^t = R \circ L^t,
\]
where we have used that $I^t = I$ since $I$ is a symplectomorphism which equals its own inverse.
\end{proof}

\begin{thrm}\label{thrm:symp-orb}
Let $S \rightrightarrows C$ be the symplectic hopfoid corresponding to a symplectic double groupoid $S$, and let $D \rightrightarrows Y$ be its induced groupoid. Then $Y$ is a coisotropic submanifold of $C$ whose characteristic leaves are the orbits of the induced groupoid. Moreover, the restriction of the symplectic form on $C$ to $Y$ has the property that its pullback under the target and source of the induced groupoid agree.
\end{thrm}

\begin{proof}
The source, target, and antipode of $S \rightrightarrows C$ satisfy the conditions of the previous proposition, and the base $Y$ of the induced groupoid is the common domain and image of $L \circ R^t$. Thus $Y$ is coisotropic in $C$ and the characteristic leaf through $p \in Y$ is the image
\[
(L \circ R^t)(p) \subset Y,
\]
which is precisely the orbit $L(R^{-1}(p))$ of $D \rightrightarrows Y$ through $p$.

Now, the fact that
\[
L: S \to C
\]
is a canonical relation implies (via the same argument showing that a map is symplectic if and only if its graph is lagrangian) that $L^*\omega = \omega_1$, and for the same reason $R^*\omega = \omega_1$, where here $L,R: D \to Y$ are the target and source of the induced grouopid. Thus we have $L^*\omega = R^*\omega$ as claimed.
\end{proof}

In particular, when the induced groupoid presents an orbifold, the above conditions on $\omega$ are precisely the defining properties of an orbifold symplectic structure:

\begin{cor}
In the setting of the Theorem~\ref{thrm:symp-orb}, when the induced groupoid is proper with finite isotropy groups, the orbifold $\X$ it presents is symplectic.
\end{cor}

In the general case, the structure resulting from Theorem~\ref{thrm:symp-orb} should be thought of as a ``symplectic-like'' structure on the stack presented by the induced groupoid. We return to this point of view in the final chapter.

\section{Generalizing to Simplicial Symplectic Manifolds}
The procedures of the previous sections may, in certain situations, be generalized to arbitrary simplicial symplectic manifolds. Indeed, as a first example, suppose that the Lie group $G$ has a locally free, proper Hamiltonian action on the symplectic manifold $P$ with equivariant momentum map $\mu: P \to \g^*$. Then we have the simplicial symplectic manifold
\begin{center}
\begin{tikzpicture}[>=angle 90]
	\def\A{.3};
	\def\B{4};
	\def\C{6.3};

	\node at (0,0) {$\cdots$};
	\node at (2.4,0) {$T^*G \times T^*G \times P$};
	\node at (5.4,0) {$T^*G \times P$};
	\node at (7,0) {$P$};

	\tikzset{font=\scriptsize};
	\draw[->] (\A,+.18) -- (\A+.5,0+.18);
	\draw[->] (\A,+.06) -- (\A+.5,+.06);
	\draw[->] (\A,-.06) -- (\A+.5,-.06);
	\draw[->] (\A,-.18) -- (\A+.5,0-.18);
	
	\draw[->] (\B,+.12) -- (\B+.5,+.12);
	\draw[->] (\B,0) -- (\B+.5,0);
	\draw[->] (\B,-.12) -- (\B+.5,-.12);
	
	\draw[->] (\C,0+.06) -- (\C+.5,0+.06);
	\draw[->] (\C,0-.06) -- (\C+.5,0-.06);
\end{tikzpicture}
\end{center}
of Example \ref{ex:action-grpd}. Recall that the degree $1$ face morphisms are
\[
\tau: ((g,\mu(p)),p) \mapsto gp \text{ and } G \times id: ((g,0),p) \mapsto p,
\]
and so we see that their common domain is determined by the requirement that $\mu(p) = 0$ and hence is $G \times \mu^{-1}(0)$, where we note that $\mu^{-1}(0)$ is a smooth manifold since the action is locally free.

Thus restricting to the domains on which the degree $1$ face morphisms are both defined, as in the construction of induced groupoids, produces in degrees $1$ and $0$:
\[
G \times \mu^{-1}(0) \rightrightarrows \mu^{-1}(0),
\]
where we identify the zero section of $T^*G$ with $G$ itself and the maps are simply the induced action of $G$ on $\mu^{-1}(0)$ and projection onto the second factor. Similarly, the degree $0$ degeneracy morphism
\[
P \to T^*G \times P,\ p \mapsto ((1,\eta),p) \text{ where } \eta \in \g^*
\]
of the simplicial symplectic manifold structure restricts to the map
\[
\mu^{-1}(0) \to G \times \mu^{-1}(0),\ p \mapsto (g,p).
\]
These are precisely the target, source, and unit of the action groupoid $G \times \mu^{-1}(0) \rightrightarrows \mu^{-1}(0)$, and indeed carrying out the same procedure (i.e. restricting to common domains on which canonical relations are defined) with the rest of the simplicial symplectic manifold above produces the simplicial nerve
\begin{center}
\begin{tikzpicture}[>=angle 90]
	\def\A{.3};
	\def\B{4};
	\def\C{6.8};

	\node at (0,0) {$\cdots$};
	\node at (2.4,0) {$G \times G \times \mu^{-1}(0)$};
	\node at (5.7,0) {$G \times \mu^{-1}(0)$};
	\node at (8,0) {$\mu^{-1}(0)$};

	\tikzset{font=\scriptsize};
	\draw[->] (\A,+.18) -- (\A+.5,0+.18);
	\draw[->] (\A,+.06) -- (\A+.5,+.06);
	\draw[->] (\A,-.06) -- (\A+.5,-.06);
	\draw[->] (\A,-.18) -- (\A+.5,0-.18);
	
	\draw[->] (\B,+.12) -- (\B+.5,+.12);
	\draw[->] (\B,0) -- (\B+.5,0);
	\draw[->] (\B,-.12) -- (\B+.5,-.12);
	
	\draw[->] (\C,0+.06) -- (\C+.5,0+.06);
	\draw[->] (\C,0-.06) -- (\C+.5,0-.06);
\end{tikzpicture}
\end{center}
of this action groupoid. From the general theory of Hamiltonian reduction, the base $\mu^{-1}(0)$ is a coisotropic submanifold of $P$ whose characteristic leaves are precisely the oribits of the action groupoid. This groupoid is in fact proper with finite isotropy groups, so we find that it presents a symplectic orbifold.

\begin{ex}\label{ex:triv-simp-symp}
For a very special case of the above, suppose that $G$ is the trivial group. Then the simplicial symplectic manifold $T^*G \times P \rightrightarrows P$ is nothing but the trivial one
\begin{center}
\begin{tikzpicture}[>=angle 90]
	\def\A{.3};
	\def\B{1.4};
	\def\C{2.5};

	\node at (0,0) {$\cdots$};
	\node at (1.1,0) {$P$};
	\node at (2.2,0) {$P$};
	\node at (3.3,0) {$P$};

	\tikzset{font=\scriptsize};
	\draw[->] (\A,+.18) -- (\A+.5,0+.18);
	\draw[->] (\A,+.06) -- (\A+.5,+.06);
	\draw[->] (\A,-.06) -- (\A+.5,-.06);
	\draw[->] (\A,-.18) -- (\A+.5,0-.18);
	
	\draw[->] (\B,+.12) -- (\B+.5,+.12);
	\draw[->] (\B,0) -- (\B+.5,0);
	\draw[->] (\B,-.12) -- (\B+.5,-.12);
	
	\draw[->] (\C,0+.06) -- (\C+.5,0+.06);
	\draw[->] (\C,0-.06) -- (\C+.5,0-.06);
\end{tikzpicture}
\end{center}
of Example~\ref{ex:triv-simp-symp}. The procedure above then produces the trivial simplicial manifold with $P$ in every degree, and the corresponding symplectic orbifold is nothing but the symplectic manifold $P$ itself.
\end{ex}

In general, suppose that $P_\bullet$ is a simplicial symplectic manifold equipped with a $*$-structure. Restricting to the domains on which the various canonical relations in this structure are defined also produces (under some smoothness assumptions) a simplicial manifold $N_\bullet$. To be precise, let us describe this construction in more detail.

First, Proposition~\ref{prop:orbits} implies that the domain and image of $L \circ R^t$, where $L, R: P_1 \to P_0$ are the degree $1$ face morphisms, are the same manifold $N_0$. Next, let $N_1$ denote the common domain $\dom L \cap \dom R$ of $L$ and $R$; then $L$ and $R$ restrict to maps
\[
N_1 \rightrightarrows N_0,
\]
and the degree $0$ degeneracy morphism $E: P_0 \to P_1$ restricts to a map $N_0 \to N_1$. Note that we are using the strong transversality of the compositions $L \circ E, R \circ E: P_0 \to P_0$ here to ensure that this map is well-defined. The following example makes this clear:

\begin{ex}
Let $G \rightrightarrows M$ be a Lie groupoid which is not a Lie group and consider the structure of Example~\ref{ex:bad-simp-symp}. Carrying out the above restrictions in this case produces
\[
G \rightrightarrows pt.
\]
However, the ``degeneracy map'' in this case is only the relation $pt \to G$ whose image is $M$, and so is not a true map. Hence we do not get an induced simplicial manifold in this case, as we expect since such a structure would imply that $G$ were actually a Lie group. The problem is that the composition $pt \to G \to pt$ is not strongly transversal in this more general setting, and so there is not a unique choice for the ``unit'' in $N_1$ corresponding to a point of $N_0$.
\end{ex}

Continuing in this manner (i.e. restricting to common domains on which canonical relations are defined) then produces a simplicial manifold
\begin{center}
\begin{tikzpicture}[>=angle 90]
	\def\A{.3};
	\def\B{1.4};
	\def\C{2.5};

	\node at (0,0) {$\cdots$};
	\node at (1.1,0) {$N_2$};
	\node at (2.2,0) {$N_1$};
	\node at (3.3,0) {$N_0$};

	\tikzset{font=\scriptsize};
	\draw[->] (\A,+.18) -- (\A+.5,0+.18);
	\draw[->] (\A,+.06) -- (\A+.5,+.06);
	\draw[->] (\A,-.06) -- (\A+.5,-.06);
	\draw[->] (\A,-.18) -- (\A+.5,0-.18);
	
	\draw[->] (\B,+.12) -- (\B+.5,+.12);
	\draw[->] (\B,0) -- (\B+.5,0);
	\draw[->] (\B,-.12) -- (\B+.5,-.12);
	
	\draw[->] (\C,0+.06) -- (\C+.5,0+.06);
	\draw[->] (\C,0-.06) -- (\C+.5,0-.06);
\end{tikzpicture}
\end{center}
which we denote by $N_\bullet$; the simplicial identities satisfied by the face and degeneracy maps follow from those satisfied by the canonical relations from which they arise. As with the base of the induced groupoid of Theorem~\ref{thrm:symp-orb}, the degree $0$ part $N_0$ is coisotropic in $P_0$ and hence comes equipped with a closed $2$-form $\omega$ such that $(d_0^1)^*\omega = (d_1^1)^*\omega$ where $d_0^1, d_1^1: N_1 \to N_0$ are the degree $1$ face maps of $N_\bullet$.

When $N_\bullet$ is actually the nerve of a Lie groupoid $N_1 \rightrightarrows N_0$, we will call $N_\bullet$ the \emph{induced groupoid} of $P_\bullet$.
In this case, the same results as in the above cited theorem also hold: the characteristic leaves of $N_0$ are the orbits of $N_1 \rightrightarrows N_0$ and the pullbacks of $\omega$ under the source and target of $N_1 \rightrightarrows N_0$ agree. Thus, when the induced groupoid $N_1 \rightrightarrows N_0$ is proper with finite isotropy, the orbifold it presents is symplectic.

\section{Lifting Symplectic Orbifolds}
The observations of the previous section lead naturally to the question as to when a symplectic orbifold arises from the induced groupoid of a simplicial symplectic manifold. We finish this chapter with an approach to this question and a conjecture.

Suppose that $(X,\omega)$ is a \emph{presymplectic} manifold, meaning that $\omega$ is a closed $2$-form on $X$ of constant rank. We know by a theorem of Gotay (see for example the book \cite{GS}) that $X$ then embeds coisotropically into some symplectic manifold $U$---we briefly recall the construction. The kernel of the induced bundle map
\[
\omega: TX \to T^*X
\]
is a vector bundle and we consider its dual $(\ker\omega)^*$. The result is that there is a neighborhood $U$ of the zero section in $(\ker\omega)^*$ which carries a (non-canonical) symplectic structure, and the embedding $X \hookrightarrow U$ as the zero section is the coisotropic embedding we are looking for.

Now, suppose we have a surjective submersion $f: Y \to X$ and consider the presymplectic map
\[
f: (Y,f^*\omega) \to (X,\omega).
\]
As above, $Y$ embeds coisotropically into some symplectic manifold $V$. The bundle $\ker f^*\omega$ is given by
\[
\ker f^*\omega = \{v \in TY \ |\ f_* v \in \ker \omega\} = TY \times_{TX} \ker\omega,
\]
and so we have a natural map $f_*: \ker f^*\omega \to \ker\omega$. Dualizing gives a (possibly multi-valued) relation $(\ker f^*\omega)^* \leftarrow (\ker \omega)^*$, which we can also view as a relation in the opposite direction:
\[
(f^*)^t: (\ker f^*\omega)^* \to (\ker\omega)^*.
\]
Note that since $f$ is a surjective submersion, this relation is single-valued and surjective. Note also that this relation restricts to the map $f: Y \to X$ when looking at what happens along zero sections.

\begin{ques}
Does the relation $(f^*)^t$ restrict to a canonical relation (reduction) $V \to U$?
\end{ques}

It is unclear the extent to which this result would depend on the (non-canonical) symplectic structures on $V$ and $U$. However, if true for some correctly chosen symplectic structures, then a consequence would be:

\begin{conj}\label{conj:main}
Every symplectic orbifold arises from the induced groupoid of some simplicial symplectic manifold.
\end{conj}

\begin{rmk}
One may ask whether the use of general simplicial symplectic manifolds in this conjecture is truly necessary, or is the use of symplectic hopfoids specifically enough. It is easy to see that hopfoids are not enough already in Example~\ref{ex:triv-simp-symp}, where for a general symplectic manifold $P$ there will not be a symplectic hopfoid structure on $P \rightrightarrows P$ since there is no natural choice of ``product'' $P \times P \to P$.

Also, not every simplicial symplectic manifold can be expected to produce (ignoring smoothness issues) a Lie groupoid as its induced structure. Indeed, there should be some sort of ``Kan'' conditions one can impose on a simplicial symplectic manifold which would imply this was the case.
\end{rmk}

\begin{ex}
One case in which the above conjecture is true is the following. Let $(X,\omega)$ be a presymplectic manifold equipped with a proper, locally free action of a Lie group $G$ such that $\omega$ is $G$-invariant. Suppose in addition that the orbits of the $G$-action are the leaves of the characteristic foliation of $X$ induced by $\ker\omega$. Then the action groupoid
\[
G \times X \rightrightarrows X
\]
presents a symplectic orbifold. A result of Montgomery \cite{RM} implies the following version of Gotay's result: there exists a symplectic manifold $Q$ equipped with a Hamiltonian $G$-action so that $X$ is the zero level set of the equivariant moment map $\mu: Q \to \g^*$.

Now, from this Hamiltonian action we can form a simplicial symplectic manifold $T^*G \times Q \rightrightarrows Q$ as in Example \ref{ex:action-grpd}. The discussion of the previous section then implies that the induced groupoid of this structure is the action groupoid $G \times X \rightrightarrows X$, and we recover the symplectic structure on the corresponding orbifold.
\end{ex}

It would also be interesting to see if this example generalizes to more general $G$-actions on presymplectic manifolds $(X,\omega)$. In this case, the action does not necessarily present an orbifold, but we may still attempt to realize it as the induced groupoid of some simplicial symplectic manifold.
\chapter{Further Directions}\label{chap:further}

In this final chapter we describe possible further directions one could take with the material developed in this thesis. In particular, the original motivation of developing a theory of ``symplectic stacks'' remains a key goal.

\section{Cotangent and Symplectic Stacks}
Let us recall the original motivation for this thesis described in the introduction: to give candidates for the notions of cotangent and symplectic stacks suitable for quantization problems. To this end, we have given descriptions of cotangent orbifolds and a conjectural description of symplectic orbifolds in general using the symplectic category.

Given a regular Lie groupoid $G \rightrightarrows M$, we can still ask whether the induced groupoid $G \ltimes_M N^*\F \rightrightarrows N^*\F$ of the symplectic hopfoid $T^*G \rightrightarrows T^*M$ can be considered to be a presentation of the cotangent stack of the stack presented by $G \rightrightarrows M$. There is some motivation for this idea as follows.

First, view $G \rightrightarrows M$ as an action groupoid for $G$ acting on $M$ via the target map $\ell$. As a consequence of Proposition~\ref{prop:lifted-cot-action}, this action lifts to a Hamiltonian action of the symplectic groupoid $T^*G \rightrightarrows A^*$ on $T^*M$, and one can check that the moment map for this action is the dual $\rho^*: T^*M \to A^*$ of the anchor $\rho: A \to TM$ of the Lie algebroid of $G$. Then $N^*\F$ is precisely the kernel of this moment map, or equivalently, the dual to the cokernel of $\rho$. Comparing with the usual theory of Hamiltonian reduction, this suggests that $N^*\F$ is the ``zero level set'' of the moment map $\rho^*: T^*M \to A^*$ and hence that the ``symplectic quotient'' modeled by the groupoid
\[
G \ltimes_M N^*\F \rightrightarrows N^*\F
\]
is the cotangent bundle of the ``space'' (i.e. stack) modeled by $G \rightrightarrows M$. In particular, when $G = H \ltimes M \rightrightarrows M$ is the action groupoid of a proper, free action of $H$ on $M$, the above procedure indeed reproduces the well-known fact that $N^*\F = \mu^{-1}(0)$ for the moment map $\mu: T^*M \to \h^*$ of the lifted Hamiltonian action of $H$ on $T^*M$ and
\[
T^*(M/H) \text{ is symplectomorphic to } N^*\F/H.
\]
Thus, we are led to the idea of viewing the stack $N^*\F//G$ as the cotangent stack of $M//G$.

However, this point of view is incomplete since have only used information about the orbits of $G \rightrightarrows M$ while the stack $M//G$ this presents carries much more information; in particular, $M//G$ also encodes the isotropy groups of $G \rightrightarrows M$. Because of this, it is reasonable to expect that the full ``cotangent stack'' of $X//G$ should also involve an action groupoid of the form
\[
G \ltimes_M (\ker\rho)^* \rightrightarrows (\ker\rho)^*,
\]
where $\ker\rho$ (the kernel of the anchor $\rho: A \to TM$) is the bundle consisting of the Lie algebras of the isotropy groups of $G$. Such a structure is trivial in the orbifold case, which explains why our construction of cotangent stacks makes sense in that setting.

There is another problem with this point of view, in that it only works for regular groupoids since otherwise $N^*\F$ is not smooth. Perhaps the correct way to proceed, and indeed Weinstein's original idea mentioned in the introduction, is to view the symplectic hopfoid $T^*G \rightrightarrows T^*M$ itself as a presentation of the ``cotangent stack'' of $M//G$. This object makes sense for any Lie groupoid.

\begin{ques}
Is there a sense in which the symplectic hopfoid $T^*G \rightrightarrows T^*M$ can be viewed as a model for the cotangent stack of the stack $M//G$?
\end{ques}

As a first case, an answer to this question should describe the sense in which the hopfoid $T^*G \rightrightarrows pt$, for $G$ a Lie group, is a model for the cotangent stack of the stack $pt//G$.

Turning to symplectic structures, we can also ask what sort of ``symplectic'' structure $G \ltimes_M N^*\F \rightrightarrows N^*\F$ carries. As we saw in Chapter~\ref{chap:stacks}, the induced groupoid $D \rightrightarrows Y$ of a symplectic hopfoid $S \rightrightarrows C$ comes equipped with a closed $2$-form $\omega$ on its base whose pullbacks under the source and target maps agree and whose characteristic leaves are the groupoid orbits. This latter condition is precisely the sort of condition which would imply the quotient $Y/D$ had a symplectic structure if it were actually smooth. In general, however, the condition that $\omega$ has the same pullback by the source or target map is not even equivalent to it being $D$-invariant, and so may not naturally correspond to an object on the stack $Y//D$. Indeed, the obstruction to having a good notion of ``symplectic stack'' is finding a good definition which is Morita invariant.

In addition, we again have the issue that the above condition uses only the orbits of $D \rightrightarrows Y$, while the stack $Y//D$ also reflects the isotropy groups. In general then, it seems natural to require that the notion of a symplectic structure on a stack should somehow also depend on the isotropy. In the case of $G \ltimes_M N^*\F \rightrightarrows N^*\F$ for instance, this should possibly come from a certain structure on
\[
G \ltimes_M (\ker\rho)^* \rightrightarrows (\ker\rho)^*.
\]
It is not clear at all, however, what this structure should actually be.

As in the previous discussion, perhaps the correct way to proceed is not to pass to induced groupoids, but to view the symplectic hopfoid $S \rightrightarrows C$ itself as a model of a ``symplectic stack''. In particular, from this point of view the fact that the cotangent stack of a stack has a symplectic structure is already encoded in the fact that the hopfoid $T^*G \rightrightarrows T^*M$ lives in the symplectic category. More generally, one could attempt to view simplicial symplectic manifolds in general as models for ``symplectic-like spaces''.

\begin{ques}
Is there a sense in which symplectic hopfoids, or more generally simplicial symplectic manifolds, can be viewed as models for symplectic stacks? In the case of simplicial symplectic manifolds, what is the correct ``Kan'' condition to impose to capture the sense in which these generalize symplectic hopfoids?
\end{ques}

\section{Stacky Canonical Relations}
To make sense of the idea that a symplectic hopfoid should be viewed as a model for a symplectic stack, a first step would be to construct a notion of ``Morita equivalence'' between such structures.

The following should be an example of such an equivalence. Suppose that $B$ is a Morita equivalence between two Lie groupoids $G \rightrightarrows M$ and $H \rightrightarrows N$:
\begin{center}
\begin{tikzpicture}[thick]
	\node (UL) at (0,1) {$G$};
	\node (UM) at (2,1) {$B$};
	\node (UR) at (4,1) {$H$};
	\node (ML) at (0,.5) {$\downdownarrows$};
	\node (LL) at (0,0) {$M$};
	\node (MR) at (4,.5) {$\downdownarrows$};
	\node (LR) at (4,0) {$N.$};

	\tikzset{font=\scriptsize};
	\draw[->>] (UM) to node [above] {$J_G$} (LL);
	\draw[->>] (UM) to node [above] {$J_H$} (LR);
\end{tikzpicture}
\end{center}
Applying the cotangent functor gives the following structure in the symplectic category:
\begin{center}
\begin{tikzpicture}[thick]
	\node (UL) at (0,1) {$T^*G$};
	\node (UM) at (2,1) {$T^*B$};
	\node (UR) at (4,1) {$T^*H$};
	\node (ML) at (0,.5) {$\downdownarrows$};
	\node (LL) at (0,0) {$T^*G_0$};
	\node (MR) at (4,.5) {$\downdownarrows$};
	\node (LR) at (4,0) {$T^*H_0.$};

	\tikzset{font=\scriptsize};
	\draw[->>] (UM) to node [above] {$T^*J_G$} (LL);
	\draw[->>] (UM) to node [above] {$T^*J_H$} (LR);
\end{tikzpicture}
\end{center}
According to Proposition~\ref{prop:lifted-cot-action}, $T^*B$ has Hamiltonian actions of the symplectic grouopids $T^*G \rightrightarrows A^*G$ and $T^*H \rightrightarrows A^*H$, where $AG$ and $AH$ are the Lie algebroids of $G$ and $H$ respectively, obtained by taking the cotangent lift of the $G$ and $H$-actions on $B$. In addition, one can check that the ``moment maps'' $T^*J_G$ and $T^*J_H$ are then equivariant for these actions and the lifted actions of $T^*G \rightrightarrows A^*G$ and $T^*H \rightrightarrows A^*H$ on $T^*G_0$ and $T^*H_0$ respectively. Hence, $T^*B$ can in some sense be viewed as a ``bimodule'' between the symplectic hopfoids $T^*G \rightrightarrows T^*G_0$ and $T^*H \rightrightarrows T^*H_0$. It is unclear, however, the sense in which this ``bimodule'' is then left or right principal.

In addition, the construction above only uses part of the symplectic hopfoid structure---namely, it only uses one of the symplectic groupoid structures on $T^*G$ and $T^*H$. Since these hopfoids in fact encode the full symplectic double groupoid structures on $T^*G$ and $T^*H$, one should expect that this full structure is needed when defining the correct notion of ``bimodule'' between symplectic hopfoids. Indeed, perhaps the correct definition should come from first determining the correct notion of ``bibundle'' between symplectic double groupoids, and then translating this into the language of symplectic hopfoids.

We note that the notion of an action of double Lie groupoid on a Lie groupoid has already been studied (see for example \cite{BM}), and so it should be simple to define the notion of a Hamiltonian action of a symplectic double groupoid on a symplectic groupoid by requiring that certain relations involved be lagrangian. As above, it would still remain to study the sense in which such actions can be left or right principal. One approach to this would be to take the ``categorical'' formulations of left and right principality given in \cite{B}.

Coming from another direction, we may instead try to talk first about the notion of a ``canonical relation'' between symplectic orbifolds. This concept is not so clear from the point of view of orbifolds themselves, but it becomes a natural one to consider from the point of view of Conjecture~\ref{conj:main}. Based on this, a canonical relation between symplectic orbifolds should be something which ``lifts'' to a ``canonical relation'' between the corresponding symplectic hopfoids (or more generally, simplicial symplectic manifolds). Hence again we return to:

\begin{ques}
What is the correct notion of a morphism (or of a bimodule) between symplectic hopfoids, or more generally between simplicial symplectic manifolds satisfying a Kan-like condition?
\end{ques}

\section{Quantization of Symplectic Hopfoids}
After having clarified the sense in which symplectic hopfoids may be used to understand symplectic (and cotangent in particular) stacks, the next natural problem is that of quantization. Given that such a structure is described completely in terms of canonical relations, there is a natural way to proceed: we quantize the total and base space via geometric quantization, and the canonical relations between them, to produce some kind of ``quantum'' algebraic structure. Indeed, this is the one of the main motivations for expressing constructions in symplectic geometry in terms of lagrangian submanifolds.

The question then becomes: what algebraic structure arises from the quantization of a symplectic hopfoid? In the simplest case, that of a symplectic hopfoid over a point, which is the same as a Hopf algebra object in the symplectic category, the result is a quantum group or more generally a Hopf algebra. This point of view has seen some success in \cite{SS}; see also \cite{W6}. Following this line of reasoning, symplectic hopfoids should quantize to quantum groupoids.

As a motivating example, consider the symplectic hopfoid $T^*G \rightrightarrows T^*M$ obtained as the cotangent lift of a Lie groupoid $G \rightrightarrows M$. Quantizing these canonical relations produces the following structure.\footnote{Here we consider the quantization which assigns to a cotangent bundle the Hilbert space of $L^2$ functions on its zero section, and which assigns to a cotangent lift $T^*f$ the linear map given by integrating against a delta function supported on the graph of $f$.} First, the total and base space quantize (after making some choices) to $L^2(G)$ and $L^2(M)$ respectively. Now, the target and source relations $T^*\ell$ and $T^*r$ quantize to the linear maps
\[
L, R: L^2(G) \to L^2(M)
\]
which are given by integrating along the target and source fibers of $G \rightrightarrows M$ respectively.

\begin{rmk}
Of course, such linear maps may not actually be defined on all of $L^2(G)$, so we may have to consider a smaller domain or perhaps impose some properness assumptions. Also, it may be that these spaces should be enlarged to include delta functions or more general distributions, in which case care must be taken when precisely defining these spaces and the maps below. We ignore these issues for now.
\end{rmk}

The product $T^*m: T^*G \times T^*G \rightrightarrows T^*G$ quantizes to the convolution map
\[
L^2(G) \otimes L^2(G) \to L^2(G),
\]
and the dual of the coproduct $(T^*\Delta)^t: T^*G \times T^*G \to T^*G$ quantizes to the usual product of functions. In the end, all of these maps fit together to form a Hopf-algebroid like object
\[
L^2(G) \rightrightarrows L^2(M).
\]
In particular, quantizing the symplectic hopfoid $T^*G \rightrightarrows pt$ obtained from a Lie group $G$ produces the standard Hopf algebra structure on $C^\infty(G)$.

\begin{rmk}
We may also first reverse---by taking transposes---the directions of all structure morphisms of $T^*G \rightrightarrows T^*M$ to get an object of the form $T^*M \rightrightarrows T^*G$, and then attempt to quantize this. This would produce an algebraic object of the form
\[
L^2(M) \rightrightarrows L^2(G),
\]
where the maps are now going in the directions they normally would in an honest Hopf algebroid. Indeed, this now more closely resembles the standard Hopf algebroid structure on $C^\infty(M) \rightrightarrows C^\infty(G)$ obtained by taking the induced map on functions of the structure maps of $G \rightrightarrows M$.
\end{rmk}

\begin{rmk}
The above was only a vague description of a possible quantization procedure. Defining such a procedure rigorously would be an important step in any future work. In particular, seeing that we may have to use not only functions but possibly distributions as well---in which case the linear maps produced may be only partially-defined---perhaps the correct setting for such a construction is that of \emph{rigged Hilbert spaces}.
\end{rmk}

After having such a quantization procedure in hand, the problem then turns to the sense in which these quantizations describe quantizations of the ``symplectic stacks'' described by these hopfoids. In particular, in the setting of Chapter~\ref{chap:stacks}, it would be interesting to know if this procedure produces a new method for quantizing symplectic orbifolds.

To summarize, we have:

\begin{ques}
What results from the process of quantizing a symplectic hopfoid, and how can this be used to quantize symplectic and cotangent stacks?
\end{ques}


\bibliography{thesis}

\bibliographystyle{plain}

\appendix
\chapter{Groupoids and Stacks}\label{appendix}
In this appendix we collect some necessary background on Lie and symplectic groupoids. We also review some basics of the theory of differentiable stacks, focusing mainly on the case of orbifolds.

\section{Lie Groupoids}
As mentioned in the introduction, a groupoid is a small category in which every morphism is invertible, and a Lie groupoid is such an object with compatible smooth structures. Concretely:

\begin{defn}\label{grpd}
A \emph{Lie groupoid} is a groupoid object in $\Man$ such that the source and target maps are surjective submersions. In concrete terms, a Lie groupoid $G \rightrightarrows M$ is given by the following data:
\begin{itemize}
\item smooth manifolds $G$ and $M$ called the total and base space respectively; we allow the possibility that $G$ is non-Hausdorff,
\item surjective submersions $\ell,r: G \to M$ called the target and source respectively,
\item a closed embedding $e: M \to G$ called the unit embedding,
\item a smooth map $m: G \times_M G \to G$ called the product, where the fiber product is taken over $r$ and $\ell$, and
\item a diffeomorphism $i: G \to G$ called the inverse
\end{itemize}
satisfying the following requirements:
\begin{enumerate}[(i)]
\item $\ell \circ e = id_M = r \circ e$, $\ell \circ m = \ell \circ pr_1$, and $r \circ m = r \circ pr_2$,
\item (associativity) the diagram
\begin{center}
\begin{tikzpicture}[>=angle 90]
	\node (UL) at (0,1) {$G \times_M G \times_M G$};
	\node (UR) at (5,1) {$G \times_M G$};
	\node (LL) at (0,-1) {$G \times_M G$};
	\node (LR) at (5,-1) {$G$};

	\tikzset{font=\scriptsize};
	\draw[->] (UL) to node [above] {$id \times m$} (UR);
	\draw[->] (UL) to node [left] {$m \times id$} (LL);
	\draw[->] (UR) to node [right] {$m$} (LR);
	\draw[->] (LL) to node [above] {$m$} (LR);
\end{tikzpicture}
\end{center}
commutes,
\item (left and right unit) the diagrams
\begin{center}
\begin{tikzpicture}[>=angle 90]
	\node (UL) at (0,1) {$G \times G$};
	\node (U) at (2.5,1) {$M \times G$};
	\node (UR) at (5,1) {$G \times_M G$};
	\node (LL) at (0,-1) {$G$};
	\node (LR) at (5,-1) {$G$};
	
	\node at (6.5,0) {$\text{and}$};
	
	\node (UL2) at (8,1) {$G \times G$};
	\node (U2) at (10.5,1) {$G \times M$};
	\node (UR2) at (13,1) {$G \times_M G$};
	\node (LL2) at (8,-1) {$G$};
	\node (LR2) at (13,-1) {$G$};

	\tikzset{font=\scriptsize};
	\draw[->] (UL) to node [above] {$\ell \times id$} (U);
	\draw[->] (U) to node [above] {$e \times id$} (UR);
	\draw[->] (LL) to node [left] {$\Delta$} (UL);
	\draw[->] (UR) to node [right] {$m$} (LR);
	\draw[->] (LL) to node [above] {$id$} (LR);
	
	\draw[->] (UL2) to node [above] {$id \times r$} (U2);
	\draw[->] (U2) to node [above] {$id \times e$} (UR2);
	\draw[->] (LL2) to node [left] {$\Delta$} (UL2);
	\draw[->] (UR2) to node [right] {$m$} (LR2);
	\draw[->] (LL2) to node [above] {$id$} (LR2);
\end{tikzpicture}
\end{center}
commute, and
\item (left and right inverse) the diagrams
\begin{center}
\begin{tikzpicture}[>=angle 90]
	\node (UL) at (0,1) {$G \times G$};
	\node (UR) at (4,1) {$G \times_M G$};
	\node (LL) at (0,-1) {$G$};
	\node (LR) at (4,-1) {$G$};
	\node (L) at (2,-1) {$M$};
	
	\node at (6,0) {$\text{and}$};
	
	\node (UL2) at (8,1) {$G \times G$};
	\node (UR2) at (12,1) {$G \times_M G$};
	\node (LL2) at (8,-1) {$G$};
	\node (LR2) at (12,-1) {$G$};
	\node (L2) at (10,-1) {$M$};

	\tikzset{font=\scriptsize};
	\draw[->] (UL) to node [above] {$i \times id$} (UR);
	\draw[->] (LL) to node [left] {$\Delta$} (UL);
	\draw[->] (UR) to node [right] {$m$} (LR);
	\draw[->] (LL) to node [above] {$r$} (L);
	\draw[->] (L) to node [above] {$e$} (LR);
	
	\draw[->] (UL2) to node [above] {$id \times i$} (UR2);
	\draw[->] (LL2) to node [left] {$\Delta$} (UL2);
	\draw[->] (UR2) to node [right] {$m$} (LR2);
	\draw[->] (LL2) to node [above] {$\ell$} (L2);
	\draw[->] (L2) to node [above] {$e$} (LR2);
\end{tikzpicture}
\end{center}
commute,
\end{enumerate}
where $\Delta: G \to G \times G$ is the standard diagonal embedding $g \mapsto (g,g)$.
\end{defn}

When $r(g) = \ell(h)$, so that $(g,h) \in G \times_M G$, we say that $g$ and $h$ are \emph{composable} and often denote their product simply by concatenation: $gh := m(g,h)$, or as a composition: $g \circ h := m(g,h)$.\footnote{One often thinks of elements of $G$ as arrows with endpoints denoted by elements of $M$. Under our conventions, arrows go from right to left} Similarly, we often denote the inverse $i(g)$ of $g$ by $g^{-1}$. The conditions in (i) simply say that
\[
\ell(e(p)) = p = r(e(p)),\ \ell(gh) = \ell(g), \text{ and } r(gh) = r(h)
\]
for all $p \in M$ and all composable $g$ and $h$ in $G$.

For each $p \in M$, the \emph{isotropy group} of $G$ at $p$ is the group
\[
G_p := \{g \in G\ |\ \ell(g) = p = r(g)\},
\]
and the \emph{orbit} of $G$ through $p$ is $\O_p = \{\ell(g) \in M\ |\ r(g) = p\} = \ell(r^{-1}(p))$. It is a basic fact that isotropy groups are Lie groups and that orbits are injectively immersed submanifolds of $M$.

\begin{ex}
For any smooth manifold, the trivial groupoid $M \rightrightarrows M$ has all structure maps equal to the identity.\footnote{The set of composable arrows in this case is the diagonal of $M \times M$, which we identify with $M$ itself.} The pair groupoid $M \times M \rightrightarrows M$ has target and source respectively given by
\[
r(x,y) = y \text{ and } \ell(x,y) = x,
\]
unit given by $e(x) = (x,x)$, product given by $(x,y)(y,z) = (x,z)$, and inverse $i(x,y) = (y,x)$.
\end{ex}

\begin{ex}
Suppose that a Lie group $G$ acts smoothly on a manifold $M$. The \emph{action groupoid} $G \ltimes M \rightrightarrows M$ has target and source
\[
\ell(g,p) = gp \text{ and } r(g,p) = p,
\]
identity $e(p) = (1,p)$, where $1$ is the identity element of $G$, product
\[
(g,hp)(h,p) = (gh,p),
\]
and inverse $i(g,p) = (g^{-1},gp)$. The isotropy groups and orbits in this case are the usual isotropy groups and orbits of the $G$-action on $M$.
\end{ex}

\begin{ex}
A Lie groupoid $G \rightrightarrows pt$ over a point is nothing but a Lie group.
\end{ex}

\begin{defn}
A Lie groupoid $G \rightrightarrows M$ is
\begin{itemize}
\item \emph{\'etale} if $\ell$ and $r$ are \'etale maps,
\item \emph{proper} if the map $\ell \times r: G \to M \times M$ is proper, and
\item \emph{regular} if its orbits are all of the same dimension.
\end{itemize}
\end{defn}

In particular, an action groupoid $G \ltimes M \rightrightarrows M$ is proper if and only if the $G$-action on $M$ is proper. We also note that for a proper Lie groupoid, the isotropy groups are all compact and the orbits are actually \emph{embedded} submanifolds of the total space.

The following notion is often a convenient tool for transporting the groupoid structure from one point to another:

\begin{defn}
A \emph{bisection} of a Lie groupoid $G \rightrightarrows M$ is a submanifold $B$ of $G$ so that the restrictions $\ell|_B, r|_B: B \to M$ are both diffeomorphisms. A submanifold $B \subset G$ is a \emph{local bisection} if $\ell|_B$ and $r|_B$ are diffeomorphisms onto open subsets $U$ and $V$ of $M$ respectively.
\end{defn}

\begin{prop}
Any Lie groupoid admits local bisections.
\end{prop}

The obvious notion of a morphism between Lie groupoids is the following; we later give a ``better'' notion of morphism:

\begin{defn}
A \emph{homomorphism} $(F,f): G \to H$ from a Lie groupoid $G \rightrightarrows M$ to a Lie groupoid $H \rightrightarrows N$ is a pair of smooth maps
\begin{center}
\begin{tikzpicture}[>=angle 90]
	\node (U1) at (0,.5) {$G$};
	\node (U2) at (2,.5) {$H$};
	\node (L1) at (0,-.5) {$M$};
	\node (L2) at (2,-.5) {$N$};
	\node (L) at (0,0) {$\downdownarrows$};
	\node (R) at (2,0) {$\downdownarrows$};

	\tikzset{font=\scriptsize};
	\draw[->] (U1) to node [above] {$F$} (U2);
	\draw[->] (L1) to node [above] {$f$} (L2);
\end{tikzpicture}
\end{center}
making the diagram above commute when we use either the targets of $G$ and $H$ or the sources of $G$ and $H$, and such that
\[
F(gh) = F(g)F(h) \text{ for all composable } g,h \in G
\]
and
\[
F(e_G(p)) = e_H(f(p)) \text{ for all } p \in M
\]
where $e_G$ and $e_H$ are the units of $G$ and $H$ respectively. We usually denote a homomorphism simply by map $F$ on total spaces ($f$ is then the restriction of $F$ to the unit submanifolds) and can summarize the above conditions by saying that $F$ preserves targets, sources, products, and units.
\end{defn}

The infinitesimal object corresponding to a Lie groupoid is known as a Lie algebroid:

\begin{defn}
A \emph{Lie algebroid} over a smooth manifold $M$ is a vector bundle $A \to M$ equipped with
\begin{itemize}
\item an $\R$-linear Lie bracket $[\cdot\,,\cdot]: \Gamma(A) \otimes \Gamma(A) \to \Gamma(A)$ on smooth sections of $A$, and
\item a bundle map $\rho: A \to TM$ called the anchor,
\end{itemize}
so that the following Leibniz identity holds:
\[
[a,fb] = f[a,b] + (\rho(a)\cdot f)b \text{ for any $a,b \in \Gamma(A)$ and $f \in C^\infty(M)$}.
\]
\end{defn}

\begin{rmk}
A Lie algebroid should be thought of us a generalization of both a Lie algebra and of a tangent bundle; in particular, the anchor gives a way to differentiate functions using sections of the Lie algebroid.
\end{rmk}

To associate a Lie algebroid to a Lie groupoid, we use a procedure similar to the one which constructs a Lie algebra from a Lie group. The only issue is that, since the groupoid multiplication is only partially-defined, the notion of a vector field being right-invariant only makes sense for vector fields $X$ which are tangent to the $r$-fibers:

\begin{defn}
A vector field $X$ on the total space of a groupoid $G \rightrightarrows M$ is \emph{right-invariant} if it is tangent to the $r$-fibers and is $R_g$-related to itself for each $g \in G$, where
\[
R_g: r^{-1}(\ell(g)) \to r^{-1}(r(g))
\]
is right groupoid multiplication by $g$.
\end{defn}

Any such right-invariant vector field is then completely determined by its values along the unit embedding, so we have an identification
\[
\{\text{right-invariant vector fields on $G$}\} \cong \Gamma(\ker dr|_M).
\]
One can also check that $[X,Y]$ is right-invariant if $X$ and $Y$ are, so we have an induced Lie bracket on right-invariant vector fields.

\begin{defn}
Let $G \rightrightarrows M$ be a Lie groupoid. The \emph{Lie algebroid of $G$} is the vector bundle $A := \ker dr|_M$ equipped with the Lie bracket given by the above identification of sections with right-invariant vector fields on $G$ and with anchor $\rho := d\ell|_{A}: A \to TM$.
\end{defn}

Note that we could just as easily have used \emph{left-invariant} vector fields and defined a Lie algebroid structure on $\ker d\ell|_M$ instead; the result is anti-isomorphic (as a Lie algebroid\footnote{The notion of a morphism between Lie algebroids is easiest to express in the so-called \emph{supergeometric} formulation; we will not need this notion here and refer to \cite{V} for details.}) to the Lie algebroid we have defined above.

It is an easy exercise to show that the image of the anchor $\rho(A)$ is an involutive  distribution on $M$. Thus, it integrates to a (possibly singular) foliation. The following is then an easy consequence of the definition of $A$, $\rho$, and the definition of an orbit of $G$:

\begin{prop}
The foliation induced by the distribution $\rho(A)$ is the foliation on $M$ given by the orbits of $G$ if the source fibers of $G$ are connected.
\end{prop}

\begin{ex}The Lie algebroid of a trivial groupoid $M \rightrightarrows M$ is the zero bundle over $M$. The Lie algebroid of a pair groupoid $M \times M \rightrightarrows M$ is the tangent bundle $TM$ with the usual Lie bracket of vector fields and anchor $TM \to TM$ equal to the identity.
\end{ex}

\begin{ex}The Lie algebroid of an action groupoid $G \ltimes M \rightrightarrows M$ is the trivial bundle $\g \ltimes M \to M$ with bracket induced by the one on $\g$ and anchor $\g \times M \to TM$ given by the induced infinitesimal action of $\g$ on $M$.
\end{ex}

\begin{ex}The Lie algebroid of a Lie group $G \rightrightarrows pt$ is the Lie algebra $\g \to pt$ of $G$. In general, a Lie algebroid over a point is simply a Lie algebra.
\end{ex}

As opposed to the situation for (finite-dimensional) Lie algebras, a Lie algebroid need not be the Lie algebroid of a Lie groupoid. The problem of determining when a Lie algebroid comes from a Lie groupoid is referred to as the \emph{integration problem}; it was completely solved by Crainic and Fernandes in \cite{CF}.

Just as groups can act on manifolds, so too can groupoids act:

\begin{defn}
An \emph{action} of a Lie groupoid $G \rightrightarrows M$ on a smooth manifold $N$ consists of a smoth map $J: N \to M$ (called the \emph{moment map} of the action) and a smooth map $\tau: G \times_M N \to N$\footnote{The fiber product is taken over $J$ and the source $r$ of $G$, so the action of an element in $G$ should be viewed as an arrow moving an element of $N$ on the right to an element of $N$ on the left.} such that
\begin{enumerate}[(i)]
\item $J(\tau(g,n)) = \ell(g)$ for any $g \in G$ and $n \in J^{-1}(r(g))$,
\item $\tau(e(p),n) = n$ for any $p \in M$ and $n \in J^{-1}(p)$,
\item $\tau(g,\tau(h,n)) = \tau(gh,n)$ for any composable $g,h \in G$ and $n \in J^{-1}(r(h))$.
\end{enumerate}
The second condition says that units act as identities and the third expresses a compatibility between the action $\tau$ and the groupoid product. The above is actually a \emph{left} action---the notion of a right action $N \times_M G \to N$ is similarly defined.
\end{defn}

\begin{ex}
Any groupoid $G \rightrightarrows M$ acts on its base via the target $\ell$, where the moment map $M \to M$ is the identity. Similarly, $G$ acts on its base on the right via $r$.
\end{ex}

\begin{ex}
An action of a trivial groupoid $M \rightrightarrows M$ on $N$ is nothing more than a smooth map $N \to M$.
\end{ex}

\begin{ex}
An action of a action groupoid $G \ltimes M \rightrightarrows M$ on a manifold $N$ consists of an action of $G$ on $N$ so that the moment map $J: N \to M$ is equivariant with respect to this action and the action of $G$ on $M$ given by the target map.
\end{ex}

The notion of an action groupoid for a group action can easily be extended to the case of an action of a general groupoid $G \rightrightarrows M$ on a manifold $N$; the resulting groupoid is denoted $G \ltimes_M N \rightrightarrows N$.

\begin{defn}
Let $G \rightrightarrows M$ and $H \rightrightarrows N$ be Lie groupoids. A $(G,H)$-\emph{bibundle} (or \emph{bimodule}) is a smooth manifold $B$ equipped with commuting left and right actions of $G$ and $H$ respectively. A bibundle $B$ is \emph{right-principal} if the moment map $J_G: B \to M$ for the $G$-action is a surjective submersion and  the $H$-action is free and transitive on its fibers. A \emph{left-principal} bibundle is similarly defined. A bibundle is \emph{biprincipal} if it is both left and right-principal.
\end{defn}

\begin{ex}
Let $F: G \to H$ be a homomorphism from a Lie groupoid $G \rightrightarrows M$ to a Lie groupoid $H \rightrightarrows N$. Then we can naturally associate to $F$ a right-principal bibundle $\hat F$ as follows. First, the bibundle itself is the fiber product
\[
\hat F := M \times_N H = \{(m,h) \in M \times H\ |\ f(m) = \ell_H(h)\},
\]
where $f$ is the restriction of $F$ to the units and $\ell_H$ is the target map of $H$. The moment maps of the two required actions are projection onto $M$ and the composition of projection onto $H$ followed by $\ell_H$. The left $G$-action on $\hat F$ is induced by the map $F$ and the right $H$-action is induced by the source $r_H$.

Any right-principal bibundle whose left moment map $J_G: B \to M$ admits a smooth section is of this form.
\end{ex}

\begin{ex}
A right-principal bibundle between a trivial groupoid $M \rightrightarrows M$ and a groupoid $H \rightrightarrows N$ is nothing but a \emph{right-principal $H$-bundle} over $M$.
\end{ex}

It is possible to define a natural \emph{composition} of right-principal bibundles, and it can then be checked that the process which associates to a homomorphism a bibundle preserves compositions up to isomorphism. In this sense, right-principal bibundles may be viewed as generalizations of groupoid homomorphisms, and thus form the morphisms of a category whose objects are Lie groupoids.\footnote{Actually, isomorphism classes of right-principal bibundles form the morphisms of a category; or, said another way, right-principal bibundles themselves form a \emph{2-category}.} We refer to \cite{B} for further details.

\section{Symplectic Groupoids}

\begin{defn}
A \emph{symplectic groupoid} is a Lie groupoid $S \rightrightarrows P$ where $S$ is a symplectic manifold and such that the graph of the groupoid product:
\[
graph(m) := \{(g,h,gh) \in S \times S \times S\}
\]
is a lagrangian submanifold of $\overline{S}\times\overline{S}\times S$. This condition is equivalent to the requirement that
\[
m^*\omega = pr_1^*\omega + pr_2^*\omega
\]
where $pr_1,pr_2: S \times_M S \to S$ are the two projections.
\end{defn}

\begin{thrm}[Coste-Dazord-Weinstein \cite{CDW}, see also \cite{DZ}]\label{thrm:symp-grpd}
Let $S \rightrightarrows P$ be a symplectic groupoid. Then the unit embedding $P \to S$ is Lagrangian, the target and source fibers are symplectically orthogonal to each other, and there exists a unique Poisson structure on $P$ which makes the source $r$ a Poisson map and the target $\ell$ an anti-Poisson map.
\end{thrm}

The Lie algebroid of a symplectic groupoid is isomorphic to the cotangent bundle $T^*P$ equipped with the standard Lie algebroid structure induced by the Poisson structure $\pi$ on $P$. We recall the details: the anchor $T^*P \to TP$ is given by the map
\[
\xi \mapsto \pi(\xi,\cdot),
\]
and the bracket is determined by the requirement that
\[
[df,dg] = d(\pi(df,dg)).
\]

\begin{defn}
A Poisson manifold $P$ is said to be \emph{integrable} if there exists a symplectic groupoid $S \rightrightarrows P$ so that the given Poisson structure on $P$ is the one induced by $S$.
\end{defn}

\begin{thrm}
The integrability of a Poisson manifold $P$ is equivalent to the integrability of the cotangent Lie algebroid $T^*P$.
\end{thrm}

\begin{ex}
Let $G \rightrightarrows M$ be a Lie groupoid. Then there is a symplectic groupoid structure on the cotangent bundle $T^*G$. The base of this groupoid is $A^*$, the dual bundle of the Lie algebroid $A$ of $G$. The unit $A^* \hookrightarrow T^*G$ comes from the identification of $A^*$ with the conormal bundle $N^*M \subset T^*G$. Under this identification, the source $\wt r$ is determined by the requirement that
\[
\wt r(g,\xi)\big|_{\ker d\ell_{e(r(g))}} = (dL_g)_{e(r(g))}^*\left(\xi\big|_{\ker d\ell_g}\right)
\]
where $L_g: \ell^{-1}(r(g)) \to \ell^{-1}(\ell(g))$ is left multiplication by $g$; this determines the source when considering the splitting
\[
TG|_M = TM \oplus \ker (d\ell)|_M.
\]
Similarly, the target is determined by the requirement that
\[
\wt \ell(g,\xi)\big|_{\ker dr_{e(\ell(g))}} = (dR_g)_{e(\ell(g))}^*\left(\xi\big|_{\ker dr_g}\right)
\]
where $R_g: r^{-1}(\ell(g)) \to r^{-1}(r(g))$ is right multiplication by $g$. From these descriptions, it is straightforward to write down target and source maps $T^*G \to A^*$ without the requirement that $A^*$ be identified with $N^*M$.

The product $\wt m$ is determined by the requirement that
\[
\langle \wt m(\xi,\eta), uv \rangle = \langle \xi, u \rangle + \langle \eta, v \rangle
\]
for composable $(g,\xi), (h,\eta) \in T^*G$, composable $(g,u)$ and $(h,v)$ in the tangent groupoid, and where $uv$ denotes their product in the tangent groupoid. We give another description of this symplectic groupoid in Chapter~\ref{chap:symp}.
\end{ex}

As a consequence of Theorem~\ref{thrm:symp-grpd}, $A^*$ naturally inherits a Poisson structure. Indeed, this Poisson structure is present for any Lie algebroid $A$, integrable or not: the bracket of fiberwise-linear functions on $A^*$ is induced by the bracket and anchor of $A$, and this then extends uniquely to all smooth functions on $A^*$.\footnote{In fact, having such a fiberwise-linear Poisson structure on $A^*$ is equivalent to having a Lie algebroid structure on $A$.}

\begin{thrm}
The integrability of a Lie algebroid $A$ is equivalent to the integrability of the induced Poisson structure on $A^*$.
\end{thrm}

As a result of the symplectic structure, we can consider a special type of bisection:

\begin{defn}
A \emph{(local) lagrangian bisection} of a symplectic groupoid $S \rightrightarrows P$ is a (local) bisection which is lagrangian as a submanifold of $S$.
\end{defn}

One can show that local lagrangian bisections always exist. Such bisections may  be used to give an explicit expression (see for example \cite{X}) for the groupoid product on $T^*G$ in the example above.

\begin{defn}
An action of a symplectic groupoid $S \rightrightarrows P$ on a symplectic manifold $Q$ is said to be \emph{Hamiltonian} if the graph of the action map $S \times_P Q \to Q$ is a lagrangian submanifold of $\overline{S} \times \overline{Q} \times Q$.
\end{defn}

\begin{ex}
A Hamiltonian action of a Lie group $G$ on a symplectic manifold $Q$ is the same as a Hamiltonian action of the symplectic groupoid $T^*G \rightrightarrows \g^*$ on $Q$. The moment map $Q \to \g^*$ for this latter action is the moment map for the Hamiltonian action of $G$. We return to this example in Chapter~\ref{chap:symp}.
\end{ex}

A symplectic groupoid is a special case of more general \emph{Poisson groupoids}:

\begin{defn}
A \emph{Poisson groupoid} is a Lie groupoid $P \rightrightarrows M$ where $P$ is a Poisson manifold and such that the graph of the groupoid product is a coisotropic submanifold of $\overline{P} \times \overline{P} \times P$.
\end{defn}

In addition to symplectic groupoids, Poisson-Lie groups are also a special case---namely the case where $M = pt$. The following generalizes the relation between Poisson-Lie groups and Lie bialgebras:

\begin{thrm}[\cite{W2},\cite{MX}\footnote{Specifically, it was Weinstein who showed in \cite{W2} that $A^*$ has a Lie algebroid structure, and Mackenzie-Xu who then introduced Lie bialgebroids in \cite{MX}.}]
Suppose that $P \rightrightarrows M$ is a Poisson groupoid with Lie algebroid $A$. Then the dual bundle $A^*$ naturally inherits a Lie algebroid structure making $(A,A^*)$ a Lie bialgebroid.
\end{thrm}

We refer to the above cited references for further details on Poisson groupoids and Lie bialgebroids.

\section{Stacks and Orbifolds}
For our purposes, the following characterization of differentiable stacks will be sufficient. For the standard definition in terms of categories fibered in groupoids, we refer to \cite{BX}.

\begin{defn}
Two Lie groupoids $G \rightrightarrows M$ and $H \rightrightarrows N$ are said to be \emph{Morita equivalent} if there exists a biprincipal bibundle between them; such a bibundle is said to be a \emph{Morita equivalence}. A \emph{differentiable stack} is a Morita equivalence class $\X$ of Lie groupoids. A groupoid in such an equivalence class is called a \emph{presentation} of $\X$. We also use the notation $M//G$ for the stack presented by the groupoid $G \rightrightarrows M$.
\end{defn}

\begin{rmk}
Note that Morita equivalence is simply the notion of ``isomorphic'' in the category whose morphisms are right-principal bibundles. Such an equivalence induces an equivalence of categories between the categories of right-principal $G$-bundles on a manifold $Q$ and right-principal $H$-bundles on $Q$, mimicking the usual algebraic definition of Morita equivalence.
\end{rmk}

\begin{ex}
Suppose that $G$ acts properly and freely on a manifold $M$. Then the action groupoid $G \ltimes M \rightrightarrows M$ is Morita equivalent to the trivial groupoid $M/G \rightrightarrows M/G$ of the smooth quotient $M/G$. This equivalence is implemented by $M$ itself with the obvious bibundle structure. Thus the stack presented by this action groupoid---usually denoted by $[M/G]$---is simply the manifold $M/G$ itself.

For a non-free action, the quotient $M/G$ is no longer smooth and one works with the \emph{quotient stack} $[M/G]$ directly.
\end{ex}

\begin{ex}
A \emph{transitive} groupoid---meaning one which has only one orbit---is Morita equivalent to any of its isotropy groups.
\end{ex}

It is a basic fact that two Morita equivalent groupoids have homeomorphic orbit spaces and isomorphic isotropy groups.

\begin{defn}
Let $\X$ be a differentiable stack with groupoid presentation $G \rightrightarrows M$. The \emph{tangent stack} $T\X$ of $\X$ is the stack presented by the tangent groupoid $TG \rightrightarrows TM$.
\end{defn}

The study of vector fields and flows on differentiable stacks was begun in \cite{H}. A subtle point is that the tangent stack $T\X$ is in general \emph{not} a vector bundle over $\X$, but rather a ``2-vector bundle''. We refer to the previous cited references for further details.

The definition of a cotangent stack $T^*\X$ is more subtle, with no good candidate in general. One instance in which such a construction is known is in the case where $\X$ is an \emph{orbifold}.

\begin{defn}
A differentiable stack $\X$ is said to be \emph{proper} if there exists a presentation of $\X$ by a proper groupoid; $\X$ is said to be \emph{\'etale} if there exists a presentation by an \'etale groupoid.
\end{defn}

According to \cite{CM} (see also \cite{LM}), a groupoid $G \rightrightarrows M$ is Morita equivalent to an \'etale groupoid if and only if the anchor $A \to TM$ of its Lie algebroid is injective, and this is true if and only if all isotropy groups of $G$ are discrete.

\begin{defn}\label{orbifold}
An \emph{orbifold} is a proper, \'etale differentiable stack---that is, a stack presented by a proper Lie groupoid with finite isotropy groups.\footnote{Following an analogous construction in algebraic geometry, orbifolds are also called \emph{Deligne-Mumford} stacks.}
\end{defn}

We note that if $\X$ is an orbifold, then $T\X$ is actually a \emph{vector bundle} over $\X$ in a natural way. For the notion of a differential form on an orbifold, we use the following characterization from \cite{LM}:

\begin{defn}
A \emph{differential k-form} on an orbifold $\X$ is a $k$-form $\omega = \Omega^k(M)$ on the base of a groupoid presentation $G \rightrightarrows M$ of $\X$ such that $\ell^*\omega = r^*\omega$. This latter condition is equivalent to the form being invariant under the natural action of $G \rightrightarrows M$ on its base.
\end{defn}

The following proposition justifies this definition by implying that it defines a structure on the orbifold itself:

\begin{prop}\label{prop:diff-form}
If a groupoid presentation of an orbifold admits a form as in the definition above, then any Morita equivalent groupoid inherits such a form as well.
\end{prop}

The notion of a differential form on a general differentiable stack is subtle and requires more care; indeed, the condition that $\ell^*\omega = r^*\omega$ is no longer necessarily equivalent to $\omega$ being $G$-invariant. Instead, differential forms on general stacks are usually defined via sheaves; in the orbifold case this more general definition is equivalent to the one given above.

Any orbifold has a cotangent orbifold. Usually, one starts with the sheaf of differential forms on an orbifold $\X$, and shows that such objects are ``sections'' of a vector bundle over $X$.\footnote{In particular, this is what no longer holds for a general differentiable stack.} The cotangent stack $T^*\X$ is then defined to be the total space of this vector bundle. We will use the following characterization from \cite{LM}:

\begin{prop}\label{prop:cot-orb}
Let $G \rightrightarrows M$ be an arbitrary presentation of an orbifold $\X$. Then the cotangent stack $T^*\X$ of $\X$ is presented by a groupoid of the form
\[
N^*\F_\ell \cap N^*\F_r \rightrightarrows N^*\F
\]
where $\F_\ell$ is the foliation on $G$ induced by $\ell$-fibers, $\F_r$ is the foliation induced by $r$-fibers, and $\F$ is the foliation on $M$ given by orbits.
\end{prop}

We describe this groupoid structure in more detail in Chapter~\ref{chap:stacks}, where we obtain it via a new construction.

\end{document}